\documentclass[3p,number,preprint,times]{elsarticle}
\usepackage{graphicx}
\usepackage{subcaption}
\usepackage{siunitx} 
\usepackage{multirow}
\usepackage{xcolor}
\usepackage{longtable}
\usepackage{hhline}
\usepackage[linesnumbered,ruled]{algorithm2e}
\usepackage{amsmath}
\usepackage{lineno}
\usepackage{parskip}
\usepackage{adjustbox}

\graphicspath{ {./images/} }

\journal{TR-A -- Special Issue: Integration of passenger and freight transport}

\biboptions{authoryear}

\begin{document}
	
\begin{abstract}
        The continuous increase in urban deliveries and the ongoing urbanization of large cities require the development of efficient and sustainable transportation solutions. This study investigates the impact of modular vehicle concepts and the consolidation of different demand types in the route planning on the efficiency of the urban freight and passenger transportation system. Modularity is achieved by connecting multiple vehicles together to form a platoon. The consolidation of different demand types is realized by simultaneously consider passenger and freight demand in the optimization algorithm. The considered vehicles are specific for each demand type by can be connected freely, hence it is possible to transport different demand types in the same platoon. The cost terms in the problem formulation are comprised of travel time costs, travel distance costs, fleet size costs, and cost considering unserved requests. The modular vehicle operations are modeled in a novel pickup and delivery problem which is solved using CPLEX and Adaptive Large Neighborhood Search (ALNS). In an extensive scenario study, the potentials of the new modular vehicle type are explored for different spatial and temporal demand distributions. A parameter study on vehicle capacity, vehicle range and cost saving assumptions is performed to study their influence on the efficiency. The experiments carried out indicate a general cost savings of 48\% due to modularity and an additional 9\% due to consolidation. The reduction mainly stems from reduced operating costs and reduced trip duration, while the same number of requests can be served in all cases. Empty vehicle kilometers are reduced by more than 60\% by consolidation and modularity. An additional case study in Stockholm is analyzed which highlights the applicability of the modular transport system. The proposed model and optimization framework can be used by companies and policy makers to identify required fleet sizes, optimal vehicle routes and cost savings due to different types of operation and vehicle technology. 
        
		\fbox{\parbox{0.9\textwidth}{Highlights:
				\begin{itemize}
					\item extension of the existing pickup and delivery problem with modular and multi-purpose operations
					\item extensive scenario and parameter study for different spatial and temporal demand pattern
					\item total cost savings by over 50\% 
					\item reduction of empty vehicle kilometers by over 60\% 
					\item reduction of trip duration by more than 60\%
				\end{itemize}
			}
		}
\end{abstract}
	
\begin{keyword}
		Public transportation \sep
		Freight transportation \sep
		Modular vehicles \sep
		Heuristic optimization 
\end{keyword}

\begin{frontmatter}
		\title{Modular Vehicle Routing for Combined Passenger and Freight Transport 
		}

		\author{Jonas Hatzenb{\"u}hler\fnref{myfootnote}\corref{mycorrespondingauthor}}
		\ead{jonashat@kth.se}
		\author{Erik Jenelius\fnref{myfootnote}}
		\author{Gy\H{o}z\H{o} Gid\'ofalvi \fnref{secondfootnote}}
		\author{Oded Cats \fnref{myfootnote,thirdfootnote}}

		\cortext[mycorrespondingauthor]{Corresponding author:}
	
		\fntext[myfootnote]{Division of Transportation Planning, KTH Royal Institute of Technology, Stockholm, Sweden}
		\fntext[secondfootnote]{Division of Geoinformatics, Integrated Transport Research Lab (ITRL), KTH Royal Institute of Technology, Stockholm, Sweden}
		\fntext[thirdfootnote]{Department of Transport and Planning, Delft University of Technology, Delft, Netherlands}
		
		\date{Received: April, 2022}
\end{frontmatter}

\section{Introduction}
\label{sec:introduction}

In recent years the volume of e-commerce and the number of urban deliveries have steadily increased \citep{EcommerceStatisticsIndividuals}. This ongoing trend, which has been accelerated during the pandemic crisis, directly affects the number of vehicles needed to meet transportation demand, as well as the number of delivery trips made by these vehicles. To be able to serve this increasing demand without requiring more and larger vehicles, an efficient delivery system is required. Furthermore, the potential reduction of vehicles and a reduction in distance traveled, due to a new operation system, can have a positive impact on the sustainable development goals set by the European Commission \citep{sachsSustainableDevelopmentReport2021}.

As a complement to electrifying vehicle fleets \citep{aggoune-mtalaaSolvingNewUrban2015, mellobandeiraElectricVehiclesLast2019}, autonomous systems \citep{mouradIntegratingAutonomousDelivery2020, munsterAutonomousVehicleConcept2020, schlentherPotentialPrivateAutonomous2020}, and on-demand transportation services \citep{ronaldExploringCoModalityUsing2016, winterPerformanceAnalysisFleet2018, wenValueDemandInformation2019}, transportation concepts including collaboration and modularity have recently been studied to evaluate their impact on the efficiency of urban delivery systems. The term collaboration is typically used to describe the consolidation of the same type of demand, e.g., freight, between multiple operators. Such collaboration has been studied in the context of vehicle routing problems by \cite{ganstererCollaborativeVehicleRouting2018, wangCollaborationTransportationResource2018, cleophasCollaborativeUrbanTransportation2019,
ferrellHorizontalCollaborationOpportunities2020, vargasSustainableCollaborativeLogistics2020, losValueInformationSharing2020}. Past studies have used well-established game theory concepts to demonstrate potential reductions in vehicle kilometers traveled if multiple freight operators share the demand for delivery (see \cite{lozanoCooperativeGameTheory2013, guajardo_review_2016, wangCompensationProfitDistribution2018}). Vehicle routing problems that incorporate collaboration concepts focus mainly on the delivery of goods, whereas most urban trips serve the movements of people. As is highlighted by \cite{losValueInformationSharing2020}, urban freight and passenger transport system could improve their efficiency if both types of demand are consolidated into a joint vehicle routing framework. 

A related promising transportation concept are modular vehicles. \cite{peiVehicleDispatchingModular2021} describe a passenger transportation system consisting of multiple medium-sized vehicle pods that can be coupled to form variable-length platoons. In their model, the authors allow passengers to transfer between pods to further increase the flexibility of the system. The authors formulate a mixed-integer nonlinear model which is solved on small test scenarios. The authors conclude that the modular system reduces the total costs by 6.63\%, i.e. the operation costs by 25.16\% and the passenger waiting time cost by 6.16\%, compared with fixed-capacity shuttle buses. \cite{lin_autonomous_2022} describe a bimodal system which enables the integrated transportation of passenger and last-mile logistics. The vehicles are driven autonomously and can couple and de-couple individually along their route to dynamically form platoons. The authors discuss challenges and potentials of such systems and conclude that modularity and integrated transport solutions can result in more efficient, and flexible transport systems.

\begin{figure*}[ht]
	\centering
	\begin{subfigure}[b]{0.49\textwidth}
		\centering
		\includegraphics[width=\textwidth]{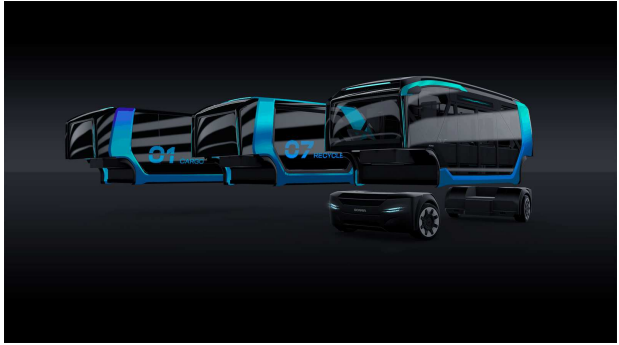}
		\caption[]%
		{Vehicle concept NXT by \cite{scaniaNXTConceptVehicle2020}}
		\label{fig:nxt}
	\end{subfigure}
	\hfill    
	\begin{subfigure}[b]{0.49\textwidth}
		\centering
		\includegraphics[width=\textwidth]{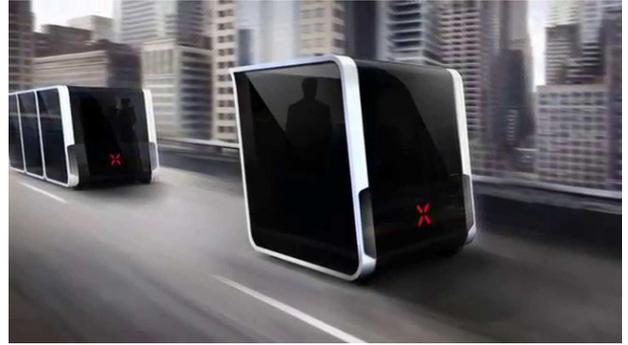}
		\caption[]%
		{Modular vehicle concept NEXT by \cite{nextHOME}}
		\label{fig:next}
	\end{subfigure}
	\caption{Illustrations of the modular vehicle concept}
	\label{fig:vehicle_concepts}
\end{figure*}

In this work consolidation and modularity are jointly modeled in a novel pickup and delivery problem. The necessary vehicle concepts for this novel operation system has been illustrated by several vehicle manufacturers. In Figure \ref{fig:vehicle_concepts} two example illustrations are shown, Figure \ref{fig:nxt} shows a vehicle concept announced by \cite{scaniaNXTConceptVehicle2020} and Figure \ref{fig:next} shows a vehicle concept designed by \cite{nextHOME}. The core innovation of the Scania NXT vehicle concept is its ability to transport different demand types (i.e. freight or passengers) using different vehicle modules. This is enabled through the exchange of modules at dedicated places. However, the total demand that can be transported in one trip is limited to the module capacity. In the NEXT vehicle concepts all modules are able to individually drive and connect to platoons, allowing for a more flexible capacity along one trip. Additionally, each platoon can be formed by multiple types of modules. Hence, this concept allows for consolidation and modularity, where the consolidation stems from the combined consideration of different demand/flow types (i.e. freight and passenger) and the modularity stems from the versatile module configuration of each platoon. 

This work focuses on the impacts of the modular multi-purpose vehicle technology when simultaneously optimizing the passenger and freight transportation system by determining optimal fleet size, vehicle routes, platoon configuration and arrival/departure times for each request. The impacts are evaluated in terms of travel distance, travel time, vehicle fleet size, and unserved demand. In order to strengthen the generality of the results an extensive array of experiments has been devised, where the scenarios differ in the temporal and spatial demand distribution. Additionally, a sensitivity analysis is performed in regard to the capacity, range, and cost parameters of the vehicles. 

The main contributions of this work can be summarized as follows. First, we present a mixed-integer program of the novel modular multi-purpose pickup and delivery problem (MMP-PDP). Second, an Adaptive Large Neighborhood Search (ALNS) algorithm is presented and used to efficiently solve the MMP-PDP. The ALNS is validated using CPLEX and smaller demo scenarios. Third, this paper extends the aforementioned studies on collaborative and modular vehicle routing problems (VRP) by performing an extensive parameter study and scenario analysis.

In the remainder of the paper, a literature review of the relevant and related works is first presented in Section \ref{sec:literature_review}. Second, the formulation of the optimization problem is provided and the ALNS algorithm described in Section \ref{sec:method}. The experimental design and parameter study are detailed in Section \ref{sec:experimental_design}. In Section \ref{sec:results} the ALNS algorithm is validated and the results are reported and discussed. The paper closes with a conclusion including a discussion of study limitations and an outlook on future research directions in Section \ref{sec:conclusion}.

\section{Literature Review}
\label{sec:literature_review}

Our work can be categorized as a specialized version of the widely studied vehicle routing problem (VRP). The general VRP describes the routing of a single vehicle serving a set of requests, where all trips start and end at a single depot. Since the problem was first described in \cite{dantzigTruckDispatchingProblem1959} many alternative formulations and extensions to the original problem have been developed. In this section we present the most relevant recent trends and variations for the proposed model. Following a brief discussion of the Pick and Delivery Problem (PDP), routing problems with modular vehicle, transfer and transshipment operations, and concepts of platooning are elaborated. The section highlights the key features of each variant, mentions key publications, and their results, while underlining the differences in relation to the model proposed in this work. The section closes with a brief discussion on multi-purpose transportation concepts and solution approaches for the respective routing problems. 

The PDP as introduced by \cite{PDP1988} and \cite{savelsberghGeneralPickupDelivery1995} is characterized by the possibility to model dedicated pickup and delivery processes for a request. This means that requests are not only delivered from a depot to their destinations but also could be picked up along the route. Adding this feature allows modelling of on-demand services and better captures the operations of passenger delivery systems. \cite{tothVehicleRoutingProblems2014, kocThirtyYearsHeterogeneous2016, dundarReviewSustainableUrban2021} provide comprehensive overviews of the various versions of VRP and PDP beyond the scope of this article. In addition to the definition of pickup and delivery nodes, the PDP typically also includes time window definitions, capacity constraints, range constraints, and multiple depot formulations. All of these features are shared with the proposed model in this work. The main distinctive feature of the model proposed in this work is the modularity of combining multiple types of demand and the possibility to form a platoon of vehicles along the route, both of which have been separately studied in previous works, but the combination of those has not hitherto been investigated. 

In a related study by \cite{hatz_multi_2022}, we propose a PDP using single multi-purpose vehicles to serve multiple demand types sequentially. In the model of that paper, modules can be exchanged during a vehicle route multiple times at dedicated depots, whereas the model in this paper forms platoon configurations at each depot before departure. The configuration of platoons cannot be changed en-route. The different vehicle operations require problem-specific destroy and repair operators to guarantee feasible solutions and efficient convergence. The purpose of a vehicle can be changed by replacing a designated passenger- or freight-specific container on the vehicle at dedicated swapping locations. The model is applied to several scenarios in Stockholm, Sweden. In \cite{hatz_multi_2022} we show that the operation of multi-purpose vehicles on consolidated demand can lead to a cost reduction of 13\%, while vehicle trip durations can be reduced by 33\% on average. In contrast to the model proposed in \cite{hatz_multi_2022}, this study considers the simultaneous transportation of multiple demand types. Additionally, the vehicle operations considered in this study allow for vehicle platoons of variable length and with different configurations.

A related research topic is the concept of vehicle platooning. \cite{bhoopalamPlanningTruckPlatoons2018} categorize different platooning strategies and create an extensive list of works for each platooning strategy. The authors distinguish between opportunistic, real-time, and scheduled planning, of which the last concept is closest related to the form of platooning considered in this work. In addition, the authors review work related to fixed and flexible routing strategies. The main motivation to form platoons in the listed works is expected cost savings due to reduced operational costs, reduced emissions, and improved traffic flow. \cite{gongCooperativePlatoonControl2018} develop a cooperative control for mixed platoons of connected, autonomous and human-driven vehicles. By controlling the speed of each vehicle in the platoon, traffic oscillations can be minimized, and traffic flow smoothness is improved. To investigate the operational cost benefit of platooning \cite{liangFuelsavingPotentialsPlatooning2014} estimate the benefits of opportunistic platoon formation and find that on long-haul trips in Europe, a spontaneous platooning rate of 1.2\% can be achieved, which in turn leads to fuel savings of 0.07\%. This result indicates that a certain level of coordinated platooning is necessary. In \cite{sokolovPlatoonFormationMaximization2017} and \cite{nourmohammadzadehFuelEfficientPlatooningHeavy2016}, the authors study scheduled platooning operations for long-haul and metropolitan areas. The works utilize heuristic and exact solvers in order to estimate the cost savings on smaller and large scenarios. \cite{nourmohammadzadehFuelEfficientPlatooningHeavy2016} identify a potential fuel reduction of up to 5\% for large scenarios. 

In addition to the freight transportation sector also the public transport sector can benefit from the formation of platoons along certain route segments. In \cite{zhangEfficiencySemiAutonomousFully2019} the authors investigate the efficiency of semi-autonomous and fully autonomous buses on trunk-and-branch networks. The authors formulate an analytical model that integrates passenger related costs (e.g. waiting time, travel time) into the total cost formulation. The authors conclude that the platooning along network corridors carries the potential to reduce the waiting time and operating costs. The theoretical cost reduction is sensitive to the demand level and the operating speed of the vehicles. For scenarios with higher speeds and higher demand along the corridor higher cost savings can be achieved. 


In the works of \cite{massonAdaptiveLargeNeighborhood2013, wolfingerLargeNeighborhoodSearch2021, dakicDesignOptimalFlexible2021, fuPickupDeliveryProblem2022, gongTransferbasedCustomizedModular2021, liuImprovingFlexrouteTransit2021, zhangModularTransitUsing2020, chen_designing_2021, chen_continuous_2022} the authors propose multiple approaches on how to integrate transfers and transshipments into modular passenger routing problems. Different approaches have been proposed to analyze this concept, such as continuous approximation, heuristics, and exact algorithms. These studies have in common to simplify the problem by either considering single-line networks or limiting the potential transfer / transshipment locations. The general conclusion of these works is that modular transportation systems can improve conventional systems, while gains are mainly present for passenger travel, since shorter travel times, reduced waiting times, and fewer transfers do not greatly affect the cost of freight transportation. 

In contrast to the works mentioned above, our work focuses on the combination of passenger and freight transportation using the same transportation system. In this study, the platoon configuration is determined before departure from the depot and remains unchanged until its arrival at the destination depot. Since complex in-vehicle transfers and transshipments are not modeled, our model allows us to study larger and more complex transportation scenarios. 

A problem similar to the modular vehicle problems mentioned above including transfers and transshipment is the truck-and-trailer routing problem (see \cite{derigsTruckTrailerRouting2013}). Here, the fleet consists of two types of vehicles, trucks and trailers. Each request has to be served with the right vehicle configuration, that is, truck or truck and trailer. In comparison with the M-VRP, the truck and trailer problem considers only three possible vehicle configurations: (1) A truck in combination with a trailer, (2) a truck which has temporarily detached its trailer, and (3) a truck without a trailer. In the M-VRP no limitations on the possible vehicle configurations are imposed. In addition, in the proposed model all vehicles can operate within the entire network, whereas in the truck and trailer problem certain vehicle configurations are prohibited to enter specific areas of the network. However, in the truck and trailer problem vehicles can change configuration during their route, e.g. de/reattach a trailer, while in the proposed model no en-route reconfigurations are possible.

Table \ref{tab:literature_comparison} gives an overview over the most relevant works as discussed in this section. The different papers are categorized based on the problem formulation, the vehicle operations, solution approaches, transportation system modelled (e.g. collaboration, vehicle platoons) and the scenario/case study analysed. 

\begin{table}
	\caption{Literature overview. Pas.: Passenger, MINLP: Mixed-integer non-linear program, GA: Genetic algorithm, ALNS: Adaptive large neighborhood search, DP: Dynamic Programming, PD: Pickup and Delivery, C: Collaborative logistics, S: Vehicle sharing, P: Vehicle platoons} 
	\label{tab:literature_comparison}
    \begin{adjustbox}{angle=90}
    \begin{tabular}{p{3.1cm}|p{2.5cm}|p{1.1cm}|p{3.7cm}|p{1.6cm}|p{2.1cm}|p{1cm}|p{0.3cm}|p{0.2cm}|p{0.2cm}|p{0.2cm}}
    \textbf{Authors}                                                          & \textbf{Operations}                          & \textbf{Demand}  & \textbf{Objective}                                                                & \textbf{Solution}                          & \textbf{Scenario}                               & \textbf{Type} & \textbf{PD} & \textbf{C} & \textbf{S} & \textbf{P}  \\ 
    \hline
    \cite{wangCollaborationTransportationResource2018}       & Multi-depots and shared vehicles & Freight               & Transportation, maintenance and depot costs          & NSGA-II~                                            & Chongqing city,      China        & MILP                & x                            & x                                & x                        & ~                          \\
    \cite{peiVehicleDispatchingModular2021}                  & Modular platoon formation                            & Pas.             & Vehicle operations and pas.
      trip times                                              & Gurobi                                              & Guangzhou, China & MINLP               & x                            & ~                                & ~                        & x                          \\
    \cite{hatz_multi_2022}                                 & Sequential multi-purpose
      vehicles                  & Pas. and freight & Vehicle operations, customer
      travel time costs, and unserved demand                      & ALNS and CPLEX                                      & Synthetic created real size
      problems          & MILP                & x                            & x                                & x                        & ~                          \\
    \cite{sokolovPlatoonFormationMaximization2017}           & Platoon formation                                    & Freight               & Fuel consumption                                                                           & POLARIS / GAMS                                    & Synthetic created                               & MINLP               & ~                            & ~                                & ~                        & x                          \\
    \cite{nourmohammadzadehFuelEfficientPlatooningHeavy2016} & Platoon
      formation                                  & Freight               & Fuel
      consumption                                                                         & LINDO
      and GA                                      & Synthetic
      created                             & ILP                 & ~                            & ~                                & ~                        & x                          \\
    \cite{massonAdaptiveLargeNeighborhood2013}               & PD with
      transfers                 & Freight               & Travel distance                                                                            & ALNS                                                & Synthetic and real instances         & MILP                & x                            & ~                                & ~                        & x                          \\
    \cite{wolfingerLargeNeighborhoodSearch2021}              & PD with
      transshipments            & Freight               & Travel and
      transshipment costs                                                           & ALNS                                                & Benchmark and case study                 & MILP                & x                            & ~                                & ~                        & ~                          \\
    \cite{fuPickupDeliveryProblem2022}                       & PD with
      transfers                 & Freight               & Travel distance, and
      waiting/transfer time & Ruin \& repair & Synthetic and real instances         & MILP                & x                            & ~                                & ~                        & x                          \\
    \cite{gongTransferbasedCustomizedModular2021}            & Modular bus systems with
      transfers                 & Pas.             & Operational
      cost, pas. travel cost, and unserved customers                          & PSO, CPLEX                                          & Chengdu, China                                  & MINLP               & x                            & ~                                & x                        & x                          \\
    \cite{liuImprovingFlexrouteTransit2021}                  & Flexroute modular transit                            & Pas.             & Pas. waiting time, travel
      time, and operation costs                                 & Heuristics and DP                 & Beijing, China                      & MILP                & x                            & ~                                & x                        & x                          \\
    \cite{zhangModularTransitUsing2020}                      & Modular transit                                      & Pas.             & Maximize the number of served
      trip requests                                              & Commercial solver                                   & Benchmark instances                             & ILP                 & x                            & ~                                & x                        & x                          \\
    \cite{chenDesigningCorridorSystems2021}                  & Modular autonomous vehicles                          & Pas.             & Pas. waiting time and total
      operational cost                                        & Branch and bound                                    & Benchmarks and Beijing, China                    & MILP                & x                            & ~                                & x                        & x                          \\
    \cite{ozturkOptimizationModelFreight2018}                & Freight delivery using urban
      rail                  & Freight               & Minimization of inventory levels
      and waiting time                                        & Heuristics and DP                & Syntetic generated                              & MILP                & x                            & ~                                & x                        & ~                          \\
    \cite{liUrbanRailService2021}                            & Collaborative urban rail                             & Pas. and freight & Revenue considering operation
      cost, loading cost, and inventory cost                     & Iterative based heuristics
      algorithm              & Two single line studies
      cases         & MILP                & ~                            & x                                & x                        & ~                          \\
    This work                                                                 & Modular vehicle platoons                             & Pas. and Freight & Vehicle operations, customer travel time costs, and unserved demand                        & ALNS and CPLEX                                      & Synthetic created real size problems            & MILP                & x                            &         x                         & x                        & x                         
    \end{tabular}
    \end{adjustbox}
    \end{table}

During the last decade, several simultaneous consolidation approaches as described above (i.e. transportation of passenger and freight request in the same platoon) have been tested and demoed in pilot projects. In \citet{steadieseifiMultimodalFreightTransportation2014}, \citet{cochraneMovingFreightPublic2017}, \citet{ozturkOptimizationModelFreight2018}, \citet{behiriUrbanFreightTransport2018a} and \citet{liUrbanRailService2021}, the authors provide several examples of consolidation operations. These transport systems have been mostly limited to specific routes and purposes, i.e., grocery delivery using a dedicated tram line from the warehouse to the central store or waste transportation from dedicated collection points in urban areas to the recycling station outside the city center. Typically, these applications were rail bound and dependent on political support and financial subsidies; hence none of these pilots has remained operational. A recent summary of sustainable collaboration examples and problems can be found in \cite{alouiSystematicLiteratureReview2021}. Despite unsuccessful long-term demonstrations \citet{losValueInformationSharing2020} highlight the large potential for integration and collaboration in modern transportation systems. Especially, the utilization of vehicles for multiple purposes is expected to yield high efficiency gains in an urban environment. 

The mentioned literature, the past pilot projects, and recent upcoming vehicle concepts all point to the high need and potential for efficient models to evaluate the impact of modular vehicle systems for multi-purpose transportation systems. Since the practical operation of modular platoons that can perform en-route reconfiguration is technically not yet feasible, we focus on determining the platoon configuration before departure. This reduces the complexity of the problem and allows for the study of larger and more complex transportation scenarios. Additionally, we close the identified research gap in modular vehicle systems for multi-purpose transportation systems by developing a combinatorial mixed linear integer problem which models the consolidated pick up and delivery of passenger and freight requests in an urban environment. In addition to the optimal route, our model also identifies the optimal platoon configuration to maximize vehicle utilization. Additionally, we propose a cost formulation which accounts for both, travel time and travel distance related costs, by that we account for the perspective of passengers and freight.
\section{Methodology}
\label{sec:method}

The transportation system modeled in this article is a variation of the pickup and delivery problem and consists of a number of available vehicle modules. Each module is manually driven and has a certain capacity and range. Modules can be of two types: the first module type can only transport passengers, and the second module type can only transport freight. Multiple modules can form a platoon, which allows the simultaneous transportation of passenger and freight. The two request types are spatially separated, and transfers between modules within a platoon are not possible. A platoon configuration is created at the beginning of a route and this configuration remains unchanged throughout the entire route. If multiple modules form a platoon, then only one driver is required, resulting in a reduction in operational costs compared to operating each module individually. The objective of this operation is to utilize the flexible capacity configurations and the reduction in operational costs to create a more efficient transportation system. The number of modules in each platoon and their type are optimized to serve the demand in the most efficient way. In the problem considered, not every request needs to be served, but unserved requests are penalized in the objective function.

This section focuses on highlighting problem-specific challenges to first formulate the proposed pickup and delivery problem and then develop an efficient solution algorithm. The section closes with a brief description of the used performance indicator and the performed optimization parameter tuning approach. 

\subsection{Proof of concept}
\begin{figure*}[ht]
	\centering
	\begin{subfigure}[b]{0.3\textwidth}
		\centering
		\includegraphics[width=\textwidth]{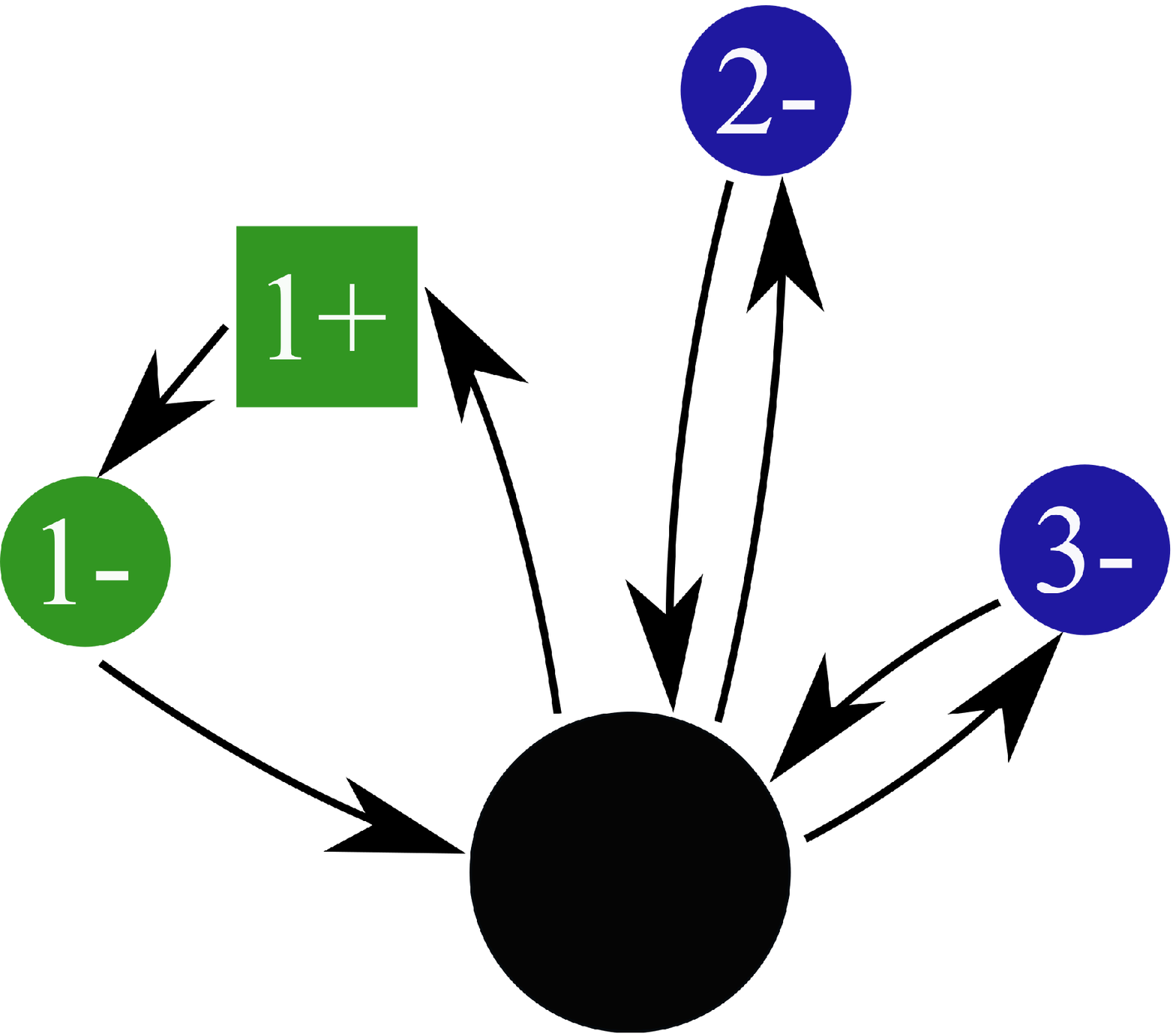}
		\caption[]%
		{Solution assuming conventional vehicle operations}
		\label{fig:poc_conventional}
	\end{subfigure}
	\hfill    
	\begin{subfigure}[b]{0.3\textwidth}
		\centering
		\includegraphics[width=\textwidth]{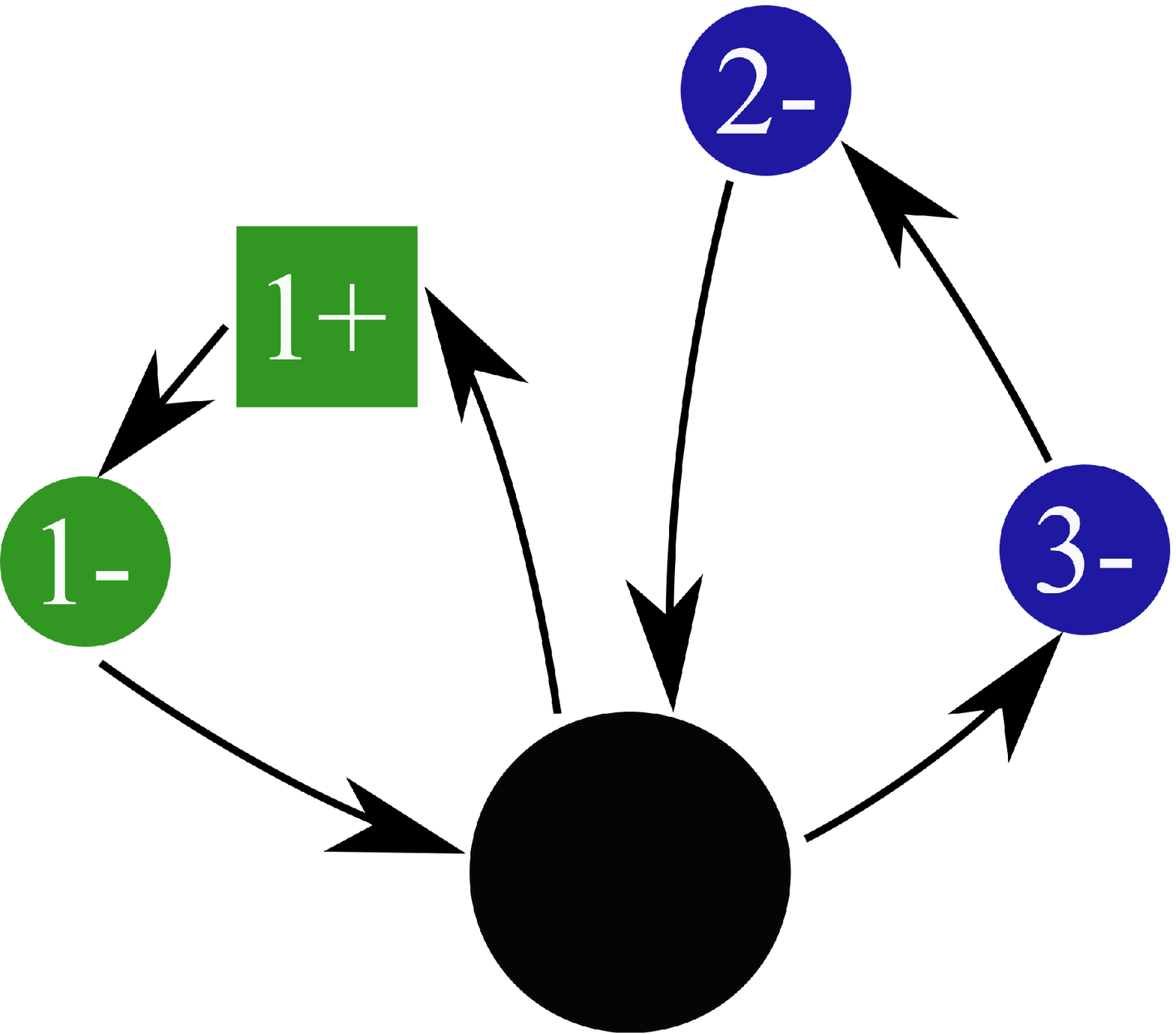}
		\caption[]%
		{Solution assuming modular vehicle operations and separated demand types}
		\label{fig:poc_modular_separated}
	\end{subfigure}
	\hfill    
	\begin{subfigure}[b]{0.3\textwidth}
		\centering
		\includegraphics[width=\textwidth]{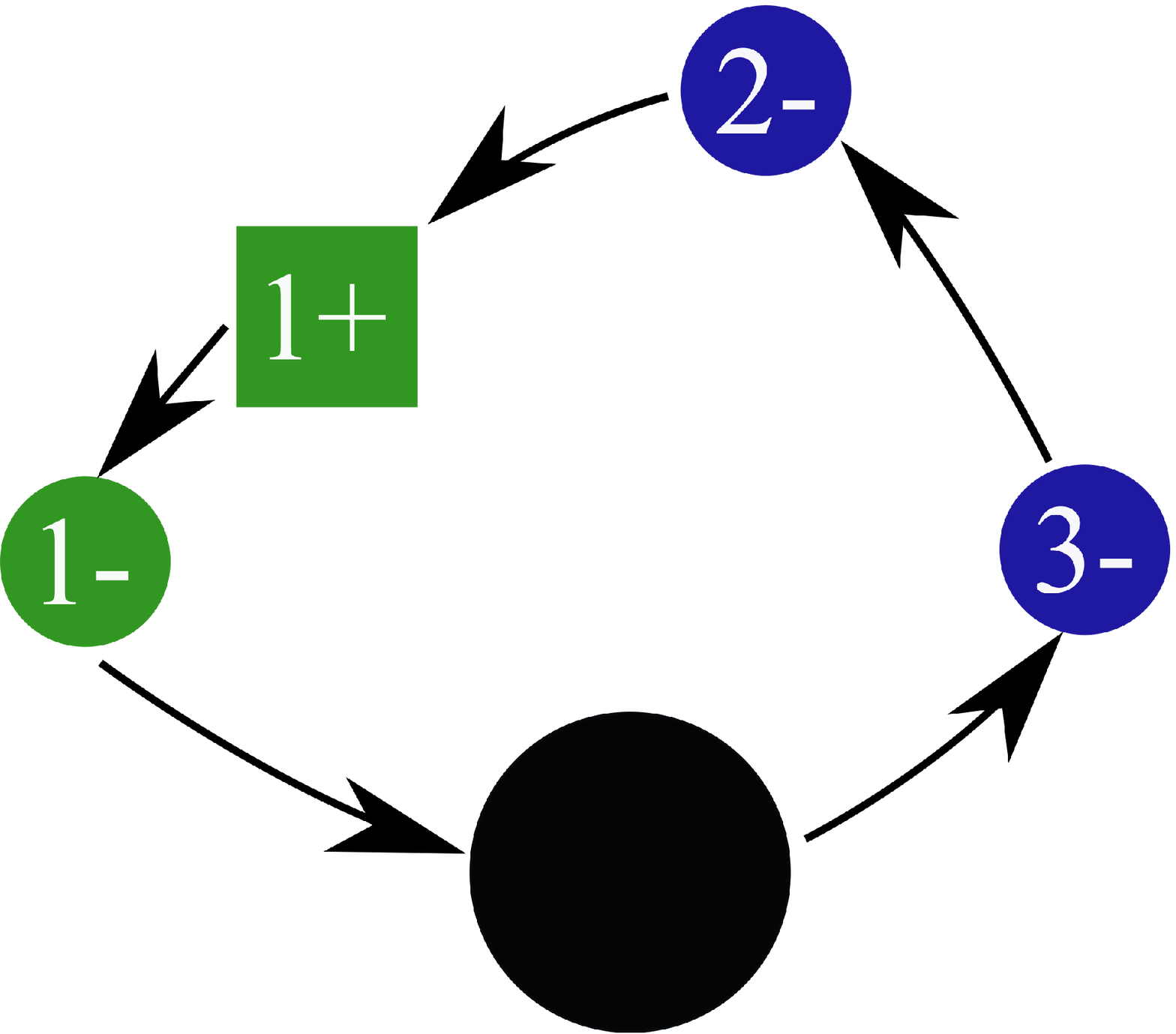}
		\caption[]%
		{Solution assuming modular vehicle operations and consolidation}
		\label{fig:poc_modular_coll}
	\end{subfigure}
	\caption{Proof of concept for different vehicle operations and consolidation}
	\label{fig:poc}
\end{figure*}

The potential benefits of the proposed transportation system are visualized in Figure \ref{fig:poc}. The scenario consists of one depot (black marker), one passenger request (green marker), and two freight requests (blue markers). The passenger request has a pickup node (square) and a drop off node (circle), while the freight requests are picked up at the depot and dropped off at the round blue nodes. The numerical values in Figure \ref{fig:poc} represent the request identifier including '+' for pickup nodes and '-' for drop off nodes. For each request a demand of 1 is assumed and the capacity of each module is also 1. For simplicity, no time windows are defined for any of the nodes. In Figure \ref{fig:poc_conventional} the optimal solution for a conventional transportation system is shown. One vehicle serves the request of the passenger and returns directly to the depot. Two vehicles are required to serve the freight requests and return to the same depot. The solution in Figure \ref{fig:poc_modular_separated} illustrates the benefits of modularity. Here only two trips are performed in the optimal solution. One trip serves the passenger demand, while the other two modules form a platoon and sequentially serve the two freight requests. Figure \ref{fig:poc_modular_coll} shows the additional benefits of consolidating the two demand types and forming a platoon, which further reduces the operational costs. Instead of simultaneously operating two vehicles with two drivers,  only one driver and one route are needed. In this example, both the total vehicle kilometers and the empty kilometers can be reduced by modularity and consolidation, which in turn reduces emissions.

The effect of flexible capacity over time due to vehicle modularity is not displayed in this simple example. If the temporal distribution of requests form peak and off-peak periods, the modular vehicle concept is able to utilize the vehicle capacity more efficiently.

\subsection{Problem formulation}

The mixed-linear program is formulated with the parameters, variables, sets, and decision variables from Table \ref{nomenclature}. Using a directed graph $G = (N,A)$ with a set of nodes $N$ and a set of arcs $A$, the proposed modular multi-purpose pickup and delivery problem (MMP-PDP) can be formulated as follows. The nodes $N$ are the union of all depot nodes ($N^{orig} \cup N^{dest}$) and nodes associated with customer requests ($N^r$), i.e., $N = N^{orig} \cup N^r \cup N^{dest}$. Here, $N^{orig}$ is the set of nodes that correspond to the depot locations. These nodes are the starting point of each trip. $N^{orig}$ is defined as $N^{orig} = \{0,...,n_d\}$ where $n_d$ is the number of depots. Nodes representing the destination of each trip are defined as $N^{dest} = \{n_d + 2\cdot n_r + 1, ...,2 \cdot n_d + 2\cdot n_r\}$, where $n_r$ is the number of requests. Note that the number of origin nodes and destination nodes is the same; hence, each origin node has one corresponding destination node with the same location. This pair of nodes represents a physical depot. Customer requests ($N^r$) are the set of nodes $N^r = \{n_d + 1, ..., n_d + 2\cdot n_r\}$. Each request has a dedicated pickup node (set $N^+$) and drop off node (set $N^-$). A request can then be modeled as a pair of nodes $(i, i+n_r)$ with $i \in N^+$ and $i+n_r \in N^-$. The arcs ($A$) of the fully connected graph contain the travel time information ($t_{i,j}$) between the nodes $i \in N$ and $j \in N$. The travel time is calculated as $t_{i,j} = d_{i,j} / v + t^s_i$, where $d_{i,j}$ is the travel distance between the nodes $i$ and $j$, $t^s_i$ is the service time at node $i \in N$, and $v$ is the constant speed of the vehicle of all platoons in the network. The service time at a node represents the duration for loading or unloading. It is independent of the level of demand and the type of request.

\begin{longtable}{c|l}
	\caption{Nomenclature and parameter values for the modular multi-purpose pickup and delivery problem (MMP-PDP)} 
	\label{nomenclature} \\
	Notation & Description \\ \hline
    $i$ & node index \\
    $j$ & node index \\
    $l$ & platoon index \\
    $p$ & platoon length index \\
    $m$ & module index \\
    $k$ & module type index \\ \hline
    $N^r$ & set of all pickup and drop off nodes for all types of request \\
    $N$ & union of $N^r$ and depot nodes \\
    $N^+$ & set of all passenger and freight pickup nodes  \\
    $N^-$ & set of all passenger and freight drop off nodes  \\
    $N^{orig}$ & set of all origin nodes (depots) \\ 
    $N^{dest}$ & set of all destination nodes (depots) \\ 
    $P$ & set of all available platoon lengths \\
    $M$ & set of all available modules \\
    $L$ & set of available platoons \\
    $W$ & set of distance related cost reduction values \\
    $U$ & set of fleet related cost reduction values \\
    $K$ & set of available module types \\ \hline
    $a_i$ & lower time-window limit at node $i \in N$ \\
    $b_i$ & upper time-window limit at node $i \in N$ \\
    $d_{i,j}$ & distance between node $i \in N$ and $j \in N$ \\
    $t_{i,j}$ & travel time between node $i \in N$ and $j \in N$ \\
    $t^s_{i}$ & service time at node $i \in N$ \\
    $q_i$ & demand at node $i \in N$ \\ \hline
    $n_r$ & number of requests \\
    $n_d$ & number of depots \\
    $H$ & large positive number  \\
    $R$ & maximum range for a platoon \\
    $Z_{max}$ & maximum number of modules in one platoon \\
    $Z_k$ & maximum number of modules per type $k \in K$ \\
    $Q_k$ & module capacity for module type $k \in K$ \\ \hline
    $\alpha_1$ & travel distance related cost parameter \\
    $\alpha_2$ & fleet size related cost parameter \\
    $\alpha_3$ & trip duration related cost parameter \\
    $\alpha_4$ & unserved requests related cost parameter \\ \hline
    $x_{i,j,l}$ & binary variable is 1 if a platoon $l \in L$ is driving between node $i$ and $j \in N$  \\
    $y_{l,m,k}$ & binary variable is 1 if platoon $l \in L$ consists of modules $m \in M$ of type $k \in K$ \\
    $e_{l,p}$ & binary variable is 1 if platoon $l \in L$ consists of $p \in P$ modules\\
    $f_{i,j,l,p}$ & binary variable is 1 if a platoon $l \in L$ of length $p \in P $ is driving between node $i$ and $j \in N$   \\
    $s_{i}$ & continuous variable describing the arrival time at node $i \in N$ \\
    $g_{l,i}$ & continuous variable describing the load for platoon $l \in L$ at node $i \in N$
\end{longtable}

The set of available module types ($K$) is defined as $K = \{0,1\}$, for passenger modules and freight modules respectively. The set can be freely extended to accommodate additional types of modules. The set of available platoons is defined as $L = \{1, ..., l_{max}\}$, where $l_{max}$ is the maximum number of platoons in a solution. Hence, we predefine the maximum number of platoons, but a feasible solution is not constrained to form $l_{max}$ platoons. A single platoon can carry out a single trip. Each platoon can contain a certain number of modules ($m \in M$) of each type, with $M = \{0, ... Z_{max}\}$. The zero is needed to allow for platoons with only one module type. The set of platoon lengths ($P$) contains all the length values that a platoon can have. In the problem formulation, the length of the platoon is given by the index $p$ and the number of modules of type ($k$) per platoon ($l$) is given by the index $m$. The maximum trip distance per platoon is $R$, hence the range is assumed to be independent of the platoon length. 

For each node $i \in N$, demand $q_i$, service duration $t^s_i$ and time windows ($a_i$, $b_i$) are defined. For depot nodes the demand and service duration are set to zero and the time windows are defined so that each depot is active for the entire planning period. Demand at the pickup and drop off nodes represents the total of items or passengers entering or leaving the modules. The demand at each drop off node is the negative value of the corresponding pickup node demand, that is $q_i + q_{i+n_r} = 0$.

In contrast to the work by \cite{hatz_multi_2022} the vehicles operating costs is formulated with respect to the platoon length. In the proposed model, two assumptions for the operation of modular vehicles are made. First, the scalar $W_p$ with $p \in P$ accounts for the reductions in distance-related costs. Due to the formation of platoons and the potential reduction of energy consumption a reduction of distance related operating costs is assumed. This goes in line with the results reported in \cite{nourmohammadzadehFuelEfficientPlatooningHeavy2016}. $W_p$ is defined as follows.

\begin{equation}
    W_p = \frac{1 + 0.95 \cdot \left(p - 1\right)}{ p} \qquad \forall p \in P \label{eq:distance}
\end{equation}

Second, it is assumed that operating costs related to fleet size, including vehicle capital costs, driver costs, and infrastructure costs (see \cite{militaoOptimalFleetSize2021}) are reduced due to the operation of modular vehicles. This cost reduction factor $U_p$ is calculated as

\begin{equation}
    U_p = \frac{1 + \left(1-\eta\right) \cdot \left(p - 1\right)}{p} \qquad \forall p \in P \label{eq:fleet_size}
\end{equation}

where $\eta$ is the train incentive parameter that models the reduction in operational costs due to forming a platoon. This cost reduction is motivated by two factors. First, driver costs can be significantly reduced if modules are operated in a platoon configuration compared to operated separately. Second, the creation of platoons may result in fewer trips, which in turn simplifies the fleet management processes. The numerical value for the parameter $\eta$ is subject to the conducted parameter study.

To model the proposed vehicle operations, six decision variables are required. The first decision variables are all binary and model the vehicle routing, platoon characteristics, and request assignment. $x_{i,j,l}$ is a binary variable that is 1 if the platoon $l \in L$ drives between nodes $i$ and $j \in N$ and is 0 otherwise. $y_{l,m,k}$ is a binary variable which is 1 if platoon $l \in L$ consists of a number of modules $m \in M$ of type $k \in K$ and 0 otherwise. $e_{l,p}$ is a binary variable which is 1 if platoon $l \in L$ consists of $p \in P$ modules and 0 otherwise. This variable is needed in order to compute the fleet size related cost term. $f_{i,j,l,p}$ is a binary variable which is 1 if platoon $l \in L$ of length $p \in P$ is driving between node $i$ and $j \in N$ and 0 otherwise. This variable is needed to linearize the distance-related cost term. The remaining two decision variables, model the arrival time $s_{i}$ at node $i \in N$ and the load $g_{l,i}$ for platoon $l \in L$ at node $i \in N$, respectively.

The objective function consists of four terms (see Equation \eqref{eq:objective}). The first term computes the cost for the total distance travelled, the second term computes the fleet size dependent cost, the third term accounts for travel time costs and the fourth cost term computes the penalty costs for unserved requests. The cost terms are scaled with the cost parameters $\alpha_{1}$, $\alpha_{2}$, $\alpha_{3}$, and $\alpha_{4}$.  

Travel time costs are computed by subtracting the depot starting time from the depot arrival times for each trip. Note that the problem formulation permits waiting at each node, therefore the travel time is not equal to the traveled distance as computed in term 1. This cost term represents the perspective of the user about using the proposed transportation system. The cost for unserved requests is computed by subtracting the served pickup nodes from the total number of requests. This value is multiplied by $\alpha_4$, which represents the average operational costs to serve an average node. 

\begin{align}
\begin{split}
	\mathop{min}_{x, s, e, f} 
	    & \quad \alpha_1 \cdot \sum_{i,j \in N, l \in L, p \in P} d_{i,j} \cdot W_p \cdot f_{i,j,l,p} \quad + 
        \quad \alpha_2 \cdot \sum_{l \in L, p \in P} p \cdot U_p \cdot e_{l,p}  \quad + \\
        & \quad \alpha_3  \cdot \sum_{i \in N^{orig}} \left( s_{i+2\cdot n_r+n_d} - s_{i}\right) \quad +  
        \quad \alpha_4 \cdot \left(n_r - \sum_{i \in N^+, j \in N, l \in L} x_{i,j,l}\right) \label{eq:objective} 
\end{split}
\end{align}

The problem formulation for the MMP-PDP is completed with the following constraints.
\begingroup
\allowdisplaybreaks
\begin{alignat}{2}
    \sum_{j \in N, l \in L} x_{i,j,l} &\leq 1 && \qquad \forall i \in N^+ \cup N^{orig}, \label{const:1} \\
    \sum_{j \in N} x_{i,j,l} &= \sum_{j \in N} x_{i+n_r,j,l} && \qquad \forall i \in N^+, \forall l \in L, \label{const:2} \\
    \sum_{i \in N^{orig}, j \in N} x_{i,j,l} &\leq 1 && \qquad \forall l \in L, \label{const:3} \\
    \sum_{j \in N} x_{i,j,l} &= \sum_{j \in N} x_{j,i+2 \cdot n_r+n_d,l} && \qquad \forall i \in N^{orig}, \forall l \in L, \label{const:4} \\
    \sum_{j \in N} x_{j,i,l} &= \sum_{j \in N} x_{i,j,l} && \qquad \forall i \in N^r, \forall l \in L, \label{const:5} \\
    x_{i,i,l} &= 0 && \qquad \forall i \in N, \forall l \in L, \label{const:6} \\
   	\sum_{j \in N} x_{i,j,l} &= 0 && \qquad \forall i \in N^{dest}, \forall l \in L, \label{const:7} \\
    s_{i} + t_{i,j} - H \cdot (1 - x_{i,j,l}) &\leq s_{j} && \qquad \forall i,j \in N, \forall l \in L, \label{const:8} \\
    s_{i} + t_{i,i+n_r} - H \cdot (1 - \sum_{j \in N} x_{i,j,l}) &\leq s_{i+n_r} && \qquad \forall i \in N^+, \forall l \in L, \label{const:9} \\
    a_{i} &\leq s_{i} \leq b_{i} && \qquad \forall i \in N, \label{const:10} \\
   	s_{i} &\leq H \cdot \sum_{j \in N, l \in L} x_{i,j,l} && \qquad \forall i \in N^{orig} \cup N^r, \label{const:11} \\
    \sum_{i,j \in N} d_{i,j} \cdot x_{i,j,l} &\leq R && \qquad \forall l \in L, \label{const:12} \\
    \sum_{m \in M, k \in K} m \cdot y_{l,m,k} &\leq Z_{max} && \qquad \forall l \in L, \label{const:13} \\
	\sum_{l \in L, m \in M} m \cdot y_{l,m,k} &\leq Z_k && \qquad \forall k \in K, \label{const:14} \\
    g_{l, i} &\leq Q_k \cdot \sum_{m \in M} m \cdot y_{l, m, k} && \qquad  \forall i \in N^r_k, \forall l \in L, \forall k \in K, \label{const:15} \\
    g_{l, j} + H \cdot \left(1 - x_{i, j, l}\right) &\geq g_{l, i} + q_{j} && \qquad  \forall i \in N, \forall j \in N^r_k, \forall l \in L, \forall k \in K, \label{const:16} \\
    g_{l, i} &= 0 && \qquad  \forall i \in N^{orig}, \forall l \in L, \label{const:17} \\
   	\sum_{m \in M} y_{l,m,k} &\leq 1 && \qquad \forall l \in L, \forall k \in K, \label{const:18} \\
    \sum_{m \in M, k \in K} m \cdot y_{l, m, k} &= \sum_{p \in P} p \cdot e_{l, p} && \qquad \forall l \in L, \label{const:19}\\
    f_{i, j, l, p} &\leq x_{i, j, l} && \qquad \forall i,j \in N, \forall l \in L, \forall p \in P, \label{const:20}\\
    f_{i, j, l, p} &\leq e_{l, p} && \qquad \forall i,j \in N, \forall l \in L, \forall p \in P, \label{const:21}\\
    x_{i, j, l} + e_{l, p} - 1 &\leq f_{i, j, l, p} && \qquad \forall i,j \in N, \forall l \in L, \forall p \in P, \label{const:22}\\
    x_{i,j,l} &= \{ 0 , 1 \} && \qquad \forall i,j \in N, \forall l \in L,\label{const:23} \\
    y_{l,m,k} &= \{ 0 , 1 \} && \qquad \forall l \in L, \forall m \in M , \forall k \in K \label{const:24} \\
    e_{l,p} &= \{ 0 , 1 \} && \qquad \forall l \in L, \forall p \in P , \label{const:25} \\
    f_{i,j,l,p} &= \{ 0 , 1 \} && \qquad \forall i,j \in N, \forall l \in L , \forall p \in P \label{const:26} \\
    \{ s_i \in \mathbb{R} & \, | \, s_i \geq 0 \} && \qquad \forall i \in N, \label{const:27} \\
    \{ g_{l,i}\in \mathbb{R} & \, | \, g_{l, i} \geq 0 \} && \qquad \forall l \in L, \forall i \in N \label{const:28}
\end{alignat}
\endgroup

The combination of constraints \eqref{const:1} and \eqref{const:2} guarantee that each request node is visited at most once and that the nodes of a request pair are served by the same vehicle. The constraints \eqref{const:3} ensure that a trip starts at an origin node, while the constraints \eqref{const:4} force each trip to terminate at the corresponding destination node. The one-sided constraints \eqref{const:1} and \eqref{const:3} allow for unserved requests. Constraint \eqref{const:5} restricts a platoon that enters a request node to also leave the node. 
In constraint \eqref{const:6} looping within the same node is prohibited while constraint \eqref{const:7} manifests that destination nodes cannot be the origin of a route. The determination of arrival times is modeled in constraints \eqref{const:8}--\eqref{const:11}. In constraint \eqref{const:8} the arrival times for two nodes along a served trip segment $(i,j)$ are modeled. Note that the travel time between nodes $i$ and $j$ includes the service time $s_i$ for node $i$. The number $H$ is a large positive number that is defined as $H = \max \left( 0, \max\left(b_i + t_{i,j} - a_j  \, \forall i,j \in N \right) \right)$. Similarly, constraint \eqref{const:9} requires that the arrival times of the pickup nodes be lower than the arrival times of the corresponding drop off nodes. The definitions of the time window are considered in constraint \eqref{const:10}. Constraint \eqref{const:11} explicitly enforces that unused depots have an arrival time of 0, which implicitly sets the arrival time at the corresponding destination node to 0 as well. This constraint is required to correctly compute the travel time costs in the third cost term. 

The range of each trip is constrained by \eqref{const:12}. In constraints \eqref{const:13} and \eqref{const:14} the maximum number of modules per trip and the maximum total number of modules per type are restricted to $Z_{max}$ and $Z_k$, respectively. 

With constraints \eqref{const:15}-\eqref{const:17} the platoon capacity is modeled. Additionally, this group of constraints connects the two main decision variables $x$ and $y$ which connects the spatial definition and the platoon configuration of a trip. The constraint \eqref{const:15} ensures that the load in a platoon at any node is less than the available capacity of that platoon. If a platoon is visiting a node, we assume that all requests at that node are either picked up or remain unserved. Since a platoon cannot visit the same node twice the demand cannot be split in multiple parts. The available capacity is calculated by multiplying the module capacity $Q_k$ by the number of modules of the same type in the platoon $m$. Constraint \eqref{const:16} ensures that the load of the arrival node $j$ is larger or equal to the load of the departure node $i$ plus the demand at node $j$ if the segment $(i,j)$ is served. This effectively models the boarding and alighting processes at each served node. Note that the demand $q_i$ at drop off nodes ($i \in N^-$) is negative. Constraint \eqref{const:17} ensures that every platoon starts its trip empty.

The group of constraints \eqref{const:18}--\eqref{const:22} are necessary to achieve a linear objective function. First, the total length of a platoon is computed. In constraint \eqref{const:18} the number of modules per type and platoon are defined, each platoon can consists of $m$ modules for each type $k$. The total length of a platoon can be registered using constraint \eqref{const:19}. In this equation the total number of modules $p$ in a platoon $l$ is set to the sum of modules $m$ for each type $k$ in the same platoon. Second, the information about the platoon length is added to the platoon route by constraints \eqref{const:20}--\eqref{const:22}. Constraint \eqref{const:20} guarantees that the helper variable $f_{i,j,l,p} = 0$, if $x_{i,j,l} = 0$. Constraint \eqref{const:21} guarantees that $f_{i,j,l,p} = 0$, if $e_{l,p} = 0$. Constraint \eqref{const:22} guarantees that $f_{i,j,l,p} = 1$, if $x_{i,j,l} = 1$ and $e_{l,p} = 1$. In combination, these three constraints ensure the desired combination of $x$ and $e$ and allow the linear formulation of the first term in Equation \eqref{eq:objective}. 

The remaining constraints \eqref{const:22}--\eqref{const:28} define the decision variable domains for platoon routing ($x_{i,j,l}$), platoon configuration ($y_{l,m,k}$), platoon length ($e_{l,p}$), route characteristics ($f_{i,j,l,p}$), node arrival times ($s_i$) and platoon load ($g_{l,i}$), respectively.

\subsection{Adaptive Large Neighborhood Search}
The proposed linear program is a variation of the widely studied vehicle routing problems. The complexity of these combinatorial optimization problems has been proven to be NP-hard \citep{lenstraComplexityVehicleRouting22}, therefore, the proposed model is also NP-hard. In order to solve larger instances of such complex problems, several heuristic algorithms have been developed over the years. In \cite{ropkeAdaptiveLargeNeighborhood2006} the Adaptive Large Neighborhood Search (ALNS) algorithm is described and presented for the first time. The authors show that the ALNS can outperform the best-known solutions of VRP benchmark problems by 50\%. The authors show high robustness, repeatability, and convergence speed of the proposed algorithm. Since then the ALNS algorithm has been adapted and applied in many optimization problems in the field of VRP (see \cite{davidpisingerLargeNeighborhoodSearch2010},  \cite{massonAdaptiveLargeNeighborhood2013}, \cite{ghilasAdaptiveLargeNeighborhood2016}, \cite{liAdaptiveLargeNeighborhood2016}, and \cite{hatz_multi_2022}).  

Due to the good performance and versatility of the ALNS optimization algorithm, we use this algorithm to solve the proposed MMP-PDP. The basic algorithm is identical to the one proposed in \cite{ropkeAdaptiveLargeNeighborhood2006}. The implemented heuristics are adjusted to consider the modularity and consolidation in the model. In Algorithm \ref{alg:1} the general outline of the ALNS implemented is given. 

\begin{algorithm} 
	\SetKwComment{Comment}{$\triangleright$\ }{}
	\KwData{passenger/freight demand, depot positions, parameter settings}
	\KwResult{best solution ($x^*$) for the MMP-PDP}
	create a feasible solution ($x$), set $x^* := x$\; \label{alg:line1}
	\Repeat{maximum number of iterations, or objective variation threshold}{
		random wheel selection (see Equation \eqref{operator_propability}) for a destroy \& a repair operator using weights\;
		create a destroyed solution ($x_d$) using the chosen destroy operator on $x$\;
		create a candidate solution ($x'$) using the chosen repair operator on $x_d$\;
		\eIf{Objective($x'$) $<$ Objective($x^*$)}
		{ \label{alg:line6}
		    set $x^* := x'$\;
		    set $x := x'$\;
			set score for chosen destroy \& repair operator to $\sigma_1$\;
		}
		{
		    \eIf{$x'$ is accepted (Simulated Annealing)}
		    {
		        set $x := x'$\;
			    \eIf{Objective($x'$) $<$ Objective($x$)}
			    {
			        set score for chosen destroy \& repair operator to $\sigma_2$\;
			    }
			    {
			    	set score for chosen destroy \& repair operator to $\sigma_3$\;
			    }
		    }
		    {
		        set score for chosen destroy \& repair operator to $\sigma_4$\;
		    }
		} \label{alg:line21}
		update weights using new operator scores (see Equation \eqref{weight_udpate}) \label{alg:line22}
	}
	return $x^*$
	\caption{General outline of ALNS framework (see \cite{ropkeAdaptiveLargeNeighborhood2006})}
	\label{alg:1}
\end{algorithm}

First, an initial feasible solution $x$, which is also the current global solution $x^*$, is created using one of the repair operators (see Line \ref{alg:line1} in Algorithm \ref{alg:1}). Second, the main iterative loop is executed. In this loop destroy and repair operators are selected based on a random wheel selection process (see Equation \eqref{operator_propability}). Applying the selected destroy and repair heuristics sequentially to the current solution creates a new candidate solution $x'$. This is done at each iteration. The new candidate solution is always feasible. In Equation \eqref{operator_propability} $p_{i,j}$ is the probability to choose operator $i$ at iteration $j$.

\begin{equation}
\begin{aligned}
    p_{i,j} = \frac{w_{i,j}}{\sum_{k \in O} w_{k,j}} \label{operator_propability}
\end{aligned}
\end{equation}

In the first iteration, the weights are set to 1. The weight for each operator is updated at the end of each iteration (see Line \ref{alg:line22} in Algorithm \ref{alg:1}) based on their performance score. In Equation \eqref{weight_udpate} this updating step is formalized, where $s_{i,j}$ is the score ($\sigma_{1}$, $\sigma_{2}$, $\sigma_{3}$, or $\sigma_{4}$) of the operator $i$ at iteration $j$. The parameter $\{ \delta \in \mathbb{R} \, | \, 0 \leq \delta \leq 1 \}$ in Equation \eqref{weight_udpate} is used to determine the decay of the weights. If $\delta$ is large the adjustment rate is slow and the new score values are not changing the overall scores much. The opposite is true for small values of $\delta$.

\begin{equation}
\begin{aligned}
    w_{i,j+1} &= w_{i,j}\cdot \delta + \left( 1 - \delta \right) \cdot s_{i,j}   \label{weight_udpate}
\end{aligned}
\end{equation}

The score is computed based on the objective function value of the new candidate solution $x'$ (see Lines \ref{alg:line6}-\ref{alg:line21} in Algorithm \ref{alg:1}). If the candidate solution is the new best global solution, the current solution and the best global solutions are updated accordingly and the score is set to $\sigma_1$. If this is not the case, the candidate solution is evaluated using a simulated annealing (SA) approach. If the candidate solution is accepted and its objective value is better than the current solution, the current solution is replaced and the operator score is set to $\sigma_2$. If the candidate solution is accepted and its objective value is worst, then the current solution, the current solution is replaced with the candidate solution, and the operator score is the to $\sigma_3$. If the candidate solution is not accepted the operator score is set to $\sigma_4$.

The main iteration loop is repeated until one of the two termination criteria is met. (1) The algorithm ends when the maximum number of iterations ($\lambda$) are computed, or (2) the algorithm ends when the changed in objective value over a certain number of iterations ($\omega$) is below a threshold ($\epsilon$). The second criteria is only active after a minimum of iterations ($\lambda_{min}$) are computed. Equation \ref{variation} formalizes the second termination criteria.

\begin{equation}
\begin{aligned}
    \frac{\sum^{i-\omega}_{j = i - 2\cdot\omega}z_j}{\sum^{i}_{j = i - \omega} z_j} - 1 &\leq \epsilon \quad for, \lambda_{min} \leq i \leq \lambda, \label{variation}
\end{aligned}
\end{equation}

\subsubsection{Heuristic Operators}

In this section, the implemented operators are briefly discussed. The discussions focus on the adjustments made to the operators to facilitate the modularity of the vehicles. In a related study by \cite{hatz_multi_2022} the authors propose heuristic operators to solve a multi-purpose PDP with sequential transportation of different demand types. In this work the operators are designed to efficiently optimize the PDP for simultaneous  transportation of multiple demand types, resulting in different destroy and repair heuristics. For a more detailed description of the basic operators and the reasoning for their implementation, we refer to the original paper from \cite{ropkeAdaptiveLargeNeighborhood2006}.

\paragraph{Destroy Operators}

Several destroy operators are implemented to diversify the search process and facilitate the intensification process. The implemented operators are \emph{Random removal}, \emph{Module removal}, route removal, \emph{Shaw removal}, and \emph{Worst removal}. The operators differ mainly in two ways. First, the number of requests which are removed from a solution is different for each operator, resulting in small or large neighborhood variations. Second, the selection process which request should be removed differs. Operators like \emph{Worst removal} or \emph{Shaw removal} remove requests based on deterministic decision rules, where operators like \emph{Random removal}, \emph{Module removal} and route removal are based on pseudo-random decisions. All destroy operators remove node pairs from a solution, hence the depot origin node and depot destination node, as well as requests pickup and drop off nodes are removed. 

\begin{itemize}
	\item \emph{Random removal}: a random selection of served requests are removed from the solution. The removed requests are added to the list of unserved requests.
	\item \emph{Module removal}: a random selection of vehicle modules is removed from the solution. The reduced platoon capacity triggers the removal of requests, which are added to the list of unserved requests.
	\item \emph{Platoon removal}: a randomly selected platoon is removed from the solution. All included requests and depots are added to the list of unserved requests. 
	\item \emph{Shaw removal}: As proposed by \cite{shawNewLocalSearch1997a} and implemented in \cite{ropkeAdaptiveLargeNeighborhood2006} this heuristics first computes a relatedness index (see Equation \eqref{similarity}) between two requests. Using this index the operator then removes highly related requests from the solution. 
	\item \emph{Worst removal}: If the objective value of a solution is reduced by removing a specific request from the solution, this request is considered to have a high contribution to the cost. The removal operator removes a number of requests, which contribute most to the objective value.
\end{itemize}

The consideration of modular platoon configurations requires some special adjustments of the above-mentioned operators. All operators, except for \emph{Module removal} and \emph{Platoon removal}, solely remove requests and/or depots from a solution, hence the platoon configuration is not affected. However, they allow exploring the order of served requests and by that define the route and arrival times of the platoon. The number of removed requests ($N$) is computed using the random distribution in Equation \eqref{n_removal} at each iteration. The number of removed request ranges between a minimum number ($\iota$) and maximum number $n_r \cdot \xi$, where $n_r$ is the total number of requests in a scenario. If the number of served requests $n_{served}$ is smaller than the expression $n_r \cdot \xi$, $N$ ranges between $\iota$ and $n_{served}$.

\begin{equation}
    N \sim U\left(\iota, ..., min\left(n_{served}, n_r \cdot \xi \right)\right), \label{n_removal}
\end{equation}

The \emph{Module removal} and \emph{Platoon removal} operator affects the fleet size and number trips performed, as well as the request order, further diversifying the solution and therefore allowing for a simultaneous optimization of platoon route, fleet size, and platoon configuration. If a module is removed from a solution, a number of randomly selected requests of the corresponding type are removed until the capacity of the removed module is met. Therefore, the number of requests removed can vary between different iterations. 

The computation of the relatedness index in the \emph{Shaw removal} operator differs slightly from the implementation as proposed by \cite{shawNewLocalSearch1997a}. In Equation \eqref{similarity} the formulation as implemented in this work is shown. The relatedness ($R_{i,j}$) between two requests $i$ and $j$ and their pickup ($a_i, a_j$) and drop off ($b_i, b_j$) nodes is computed using the distance, travel time and request load information, terms one to three respectively. In Equation \eqref{similarity} the distance term is computed using the travel distance $d_{i,j}$ between two nodes $i$ and $j$. The travel time term uses the difference between node arrival times ($s_i$) and the load term is computed as the difference between the loads ($q_i$) of two nodes $i$ and $j$. Each term is weighted using the parameter $\phi$, $\chi$, and $\psi$ respectively. 

\begin{equation}
\begin{aligned}
        R_{i,j} &= \phi \cdot \left(d_{a_i,a_j} + d_{b_i,b_j} \right) + \chi \cdot \left( |s_{a_i}-s_{a_j}| + |s_{b_i} - s_{b_j}| \right) + \psi \cdot \left( | q_i - q_j | \right), \label{similarity}
\end{aligned}
\end{equation}

\paragraph{Repair Operators}

Each repair operator has two main steps. The first step is to randomize the list of unserved requests and depots. In the second step these unserved requests and depots are inserted into the previously destroyed solution. Each repair operator has a different heuristic how the insertion is computed. If a request or depot cannot be feasible inserted into a solution, this request or node remains unserved. At the end of each repair operator a feasible solution is created, which is the new candidate solution ($x'$) (compare Algorithm \ref{alg:1}). The goal of repair operators is to intensify the search process. 

\begin{itemize}
    \item \emph{First fit insert}: a request or depot is inserted into the first feasible location of the temporary solution.
    \item \emph{Inter route insert}: a request or depot is inserted into the best feasible location within \emph{one} platoon. The platoon in which the request is inserted is the same platoon from which it has previously been removed. For the first insertion of a request/depot a random platoon is chosen.
    \item \emph{Best insert}: a request or depot is inserted into the best feasible location of \emph{all} platoons. 
\end{itemize}

\begin{figure*}[ht]
	\centering
	\begin{subfigure}[b]{0.45\textwidth}
		\centering
		\includegraphics[width=\textwidth]{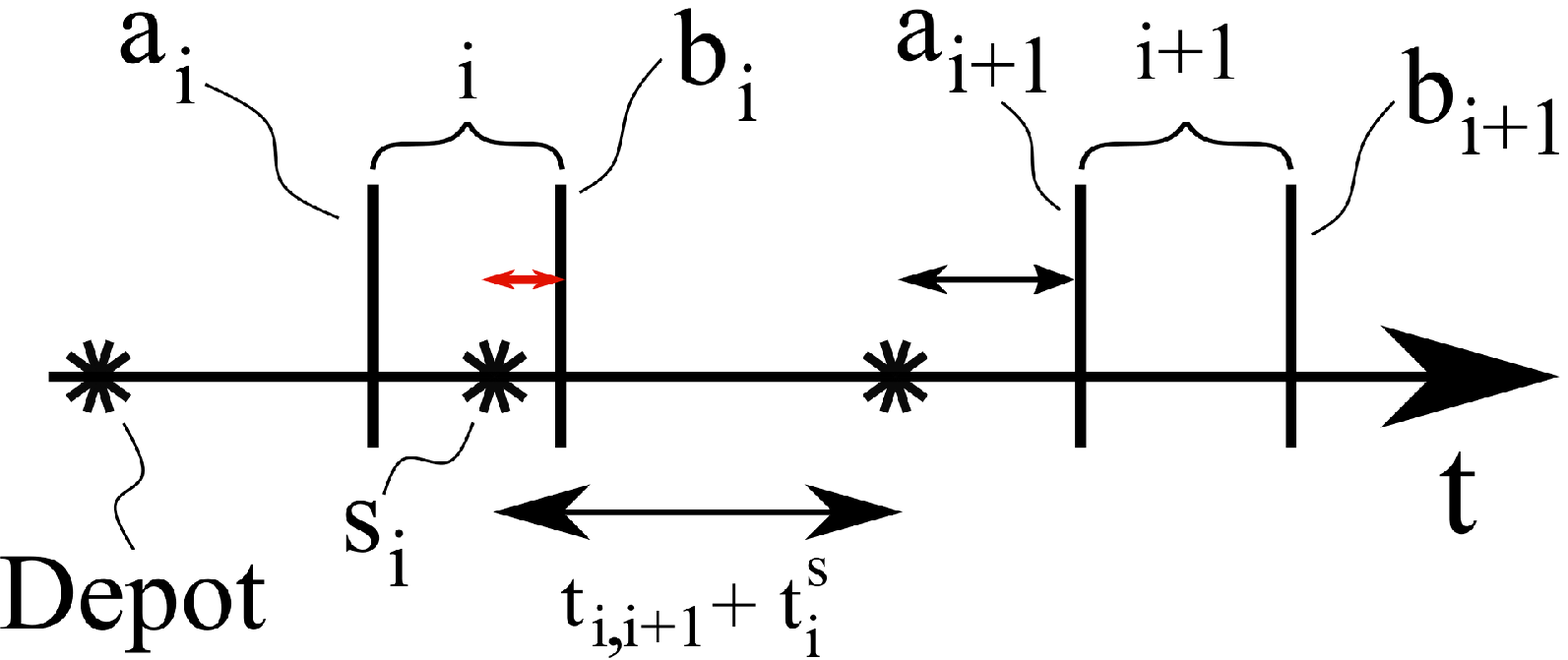}
		\caption[]%
		{Arrival times before waiting time adjustment}
		\label{fig:arrival_time_before}
	\end{subfigure}
	\hfill   
	\begin{subfigure}[b]{0.45\textwidth}
		\centering
		\includegraphics[width=\textwidth]{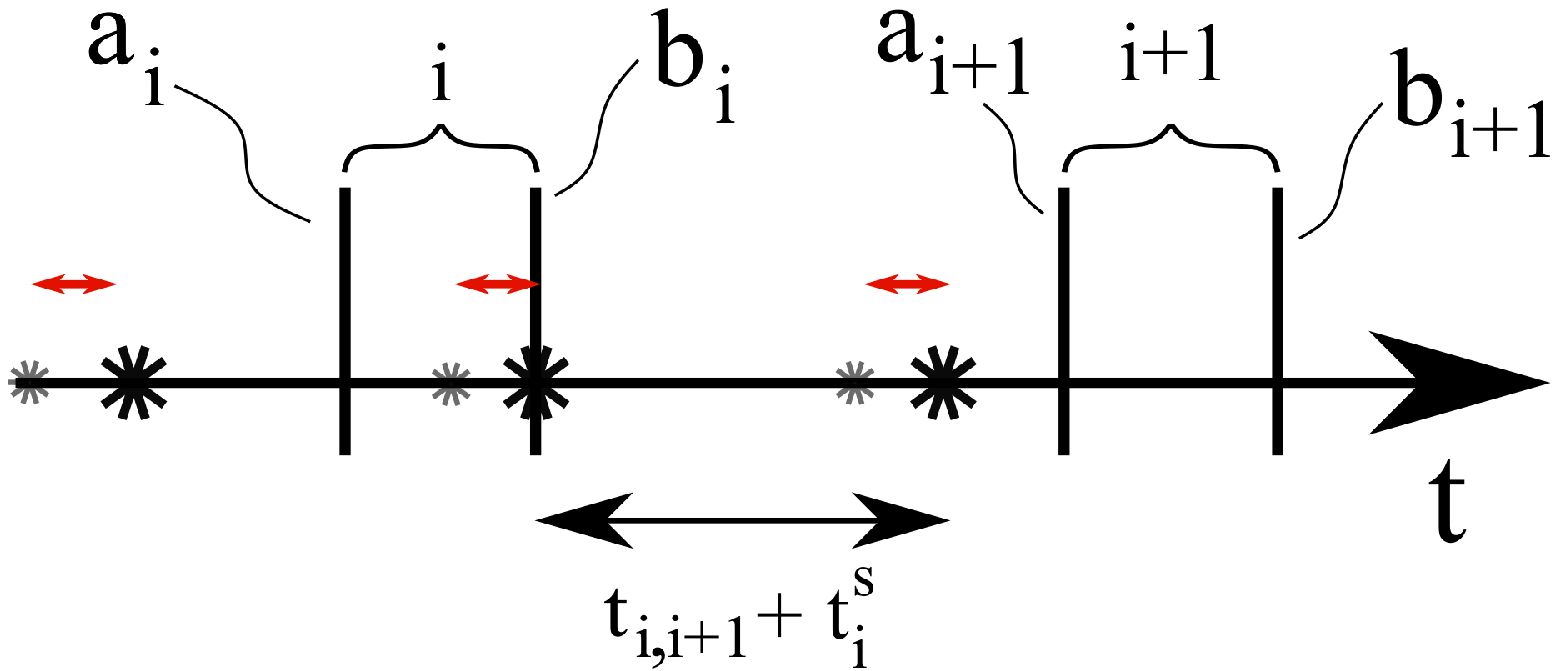}
		\caption[]%
		{Arrival times after waiting time adjustment}
		\label{fig:arrival_time_after}
	\end{subfigure}
	\caption{Minimizing waiting times to reduce trip duration}
	\label{fig:arrival_time}
\end{figure*}

\paragraph{Arrival time computation}

Due to the integration of route duration into the objective function the optimal solution for a given scenario is not only dependent on the sequence of visited nodes but also depends on the arrival time at each node. The route duration is computed as the time difference between the departure time at the origin depot and the arrival time at the destination depot. In order to guarantee the shortest route duration for a specific route and with that implicitly guaranteeing an optimal solution the arrival times at each node are determined in a two step approach. First, the arrival time of node $s_{i+1}$ is computed as the sum of arrival time at node $s_i$, the travel time $t_{i,i+1}$ and the service time $t^s_i$. Here $i$ is the index of nodes for a given route. If $i=0$ hence the beginning of a route at a depot the arrival time represents the departure time at the depot and is determined by subtracting $t_{i,i+1}$ and the service time $t^s_{i+1}$ from the lower time limit of node $i+1$ ($a_{i+1}$). If the resulting time $s_i$ is larger than $b_i$ , $s_i$ is adjusted to $s_i = b_i$ to respect the time window constraints. In a second step, the arrival time $s_{i+1}$ is evaluated. If $s_{i+1}$ is larger than $b_{i+1}$ the solution is not feasible. If $s_{i+1}$ is between $a_{i+1}$ and $b_{i+1}$ the arrival time is accepted, and the next node pair in the route is evaluated. If $s_{i+1}$ is smaller than $a_{i+1}$ this means the vehicle experiences waiting time. In such cases a recursive arrival time adjustment is performed. This calculation steps aims to minimize the duration of the route up to $i+1$ by potentially minimizing the waiting time. 

The concept of recursive arrival time adjustment is visualized in Figure \ref{fig:arrival_time}. For each arrival time the time gap to its upper time window is computed. These time gaps are compared with the waiting time and the arrival times of all nodes moved forward by the shortest time gap or the waiting time. This procedure guarantees feasible solutions while generating the shortest route duration. Both steps are performed for each node of a route until all nodes are evaluated.

\subsubsection{Parameter tuning}

For each ALNS specific parameter (see Table \ref{tab:ALNS_parameter}) a tuning procedure has been performed. The performed procedure is aligned with that described in \cite{ropkeAdaptiveLargeNeighborhood2006}, where the authors propose an iterative extensive search over a set of possible parameter values. For each parameter a set of 10 values was explored. A varied set of small and medium sized scenarios are solved using one parameter setting at a time. After all values for one parameter have been evaluated the parameter value resulting in the lowest average optimally gap with respect to the best available solution is chosen. This procedure is continued until all parameters are set. All the results in the remainder of this paper are computed with the parameter settings from Table \ref{tab:ALNS_parameter}.

\begin{table}
	\caption{ALNS parameter} 
	\label{tab:ALNS_parameter}
    \centering
    \begin{tabular}{l|l|l}
        Notation                        & Description  & Value  \\ 
        \hline
    $ \sigma_1 $        &    Operator weight for new best overall solution         &     7        \\
    $ \sigma_2 $        &     Operator weight for better candidate solution        &      2       \\
    $ \sigma_3 $        &   Operator weight for accepted solution        &       9      \\
    $ \sigma_4 $        &  Operator weight for rejected solution           &      1       \\
    $ \delta $           &      Operator decay per iteration       &       0.8       \\
    $ \lambda $          &    Iterations         &     10000        \\
    $ \lambda_{min}$ &      Min. iterations       &     5000         \\
    $ \omega $           &    Iterations look-back         &    1000          \\
    $ \epsilon $         &   Objective improvement threshold          &     0.001         \\
    $ T_{start} $                   &     Start temperature simulated annealing             &   90     \\
    $ T_{end}   $                   &     End temperature simulated annealing              &     0.0001   \\
    $ \nu $              &     Step size simulated annealing        &   0.9999           \\
    $ \phi $             &    Relatedness parameter for distance         &    9          \\
    $ \chi $             &     Relatedness parameter for time        &         4     \\
    $ \psi $             &  Relatedness parameter for load           &      9        \\
    $ \rho $             & Randomising selection parameter for \emph{Shaw removal} requests            &     6         \\
    $ \rho_{worst}$  &    Randomising selection parameter for \emph{Worst removal} requests         &     4         \\
    $ \xi $              &    Max. request removal factor         &     0.32         \\
    $ \iota$   &    Min. request removals         &      1       
    \end{tabular}
\end{table}

\section{Experimental Design}
\label{sec:experimental_design}

In this study two sets of experiments are performed. The first set evaluates the proposed transportation system with respect to the variation of temporal and spatial demand. The second set is a parameter analysis with respect to capacity, range, and vehicle operation costs. All the described scenarios and parameter studies are solved using the proposed ALNS. Each problem instance is solved as one ensemble run consisting of 10 individual optimization runs. This is to increase the generality of the presented results.

\subsection{Scenario definition}
The scenarios included in this study are created synthetically. The created scenarios aim at representing typical transportation scenarios in medium-sized to large European city. This entails an area size of \SI{12.25}{km^2} and a planning period of 6 am to 10 pm. Each scenario consists of 80 requests and 5 depots. Each request has demand $q \sim U(1,15)$ associated which follows a uniform distribution. Furthermore, the service time for each node ranges uniformly between \SI{1}{min} and \SI{5}{min}, i.e. $t^s \sim U(1,5)$. The scenarios differ in the spatial and temporal distribution of the requests.

\begin{figure*}[ht]
	\centering
	\begin{subfigure}[b]{0.49\textwidth}
		\centering
		\includegraphics[width=\textwidth]{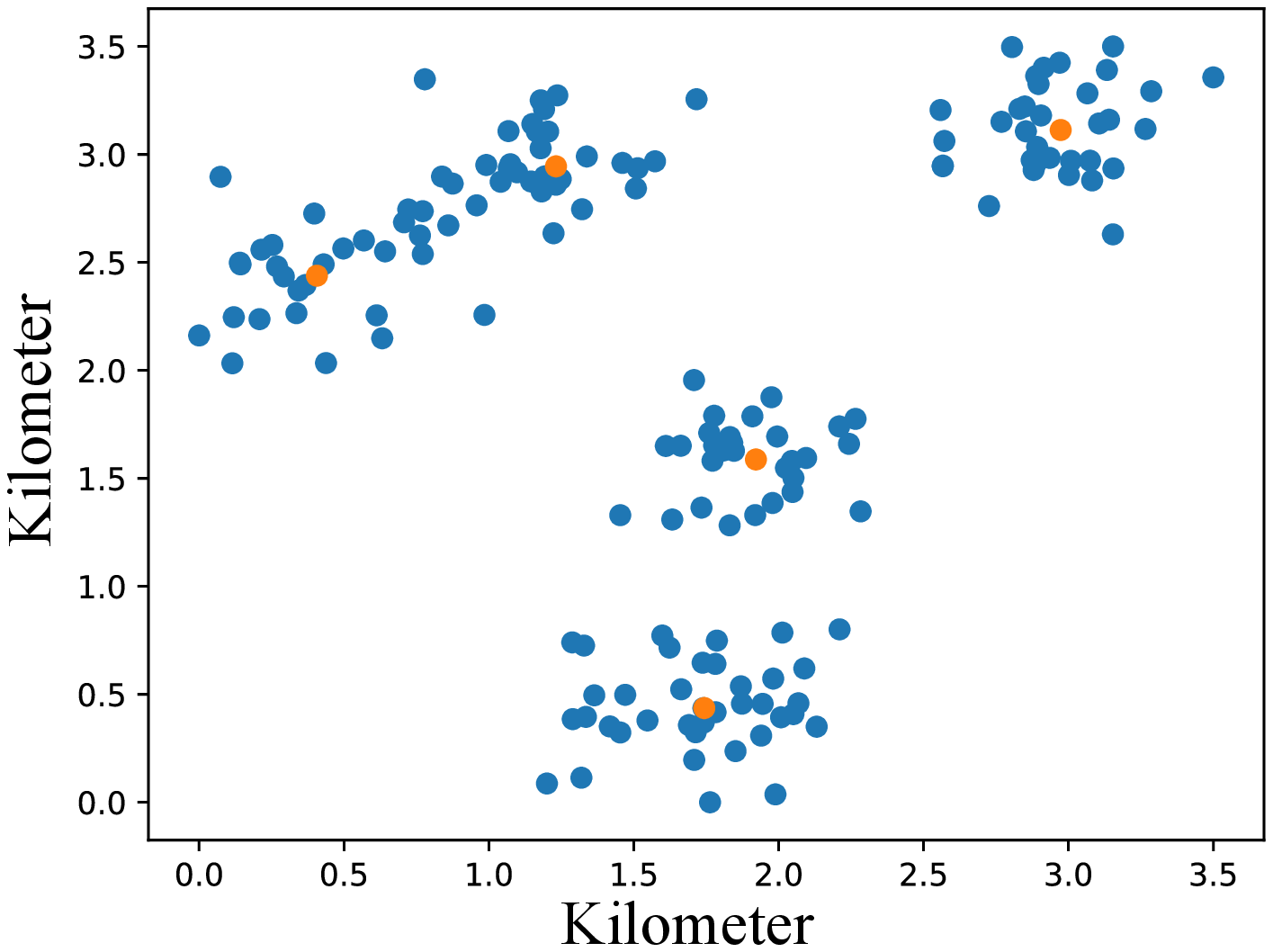}
		\caption[]%
		{clustered scenario}
		\label{fig:clustered_scenario}
	\end{subfigure}
	\hfill    
	\begin{subfigure}[b]{0.49\textwidth}
		\centering
		\includegraphics[width=\textwidth]{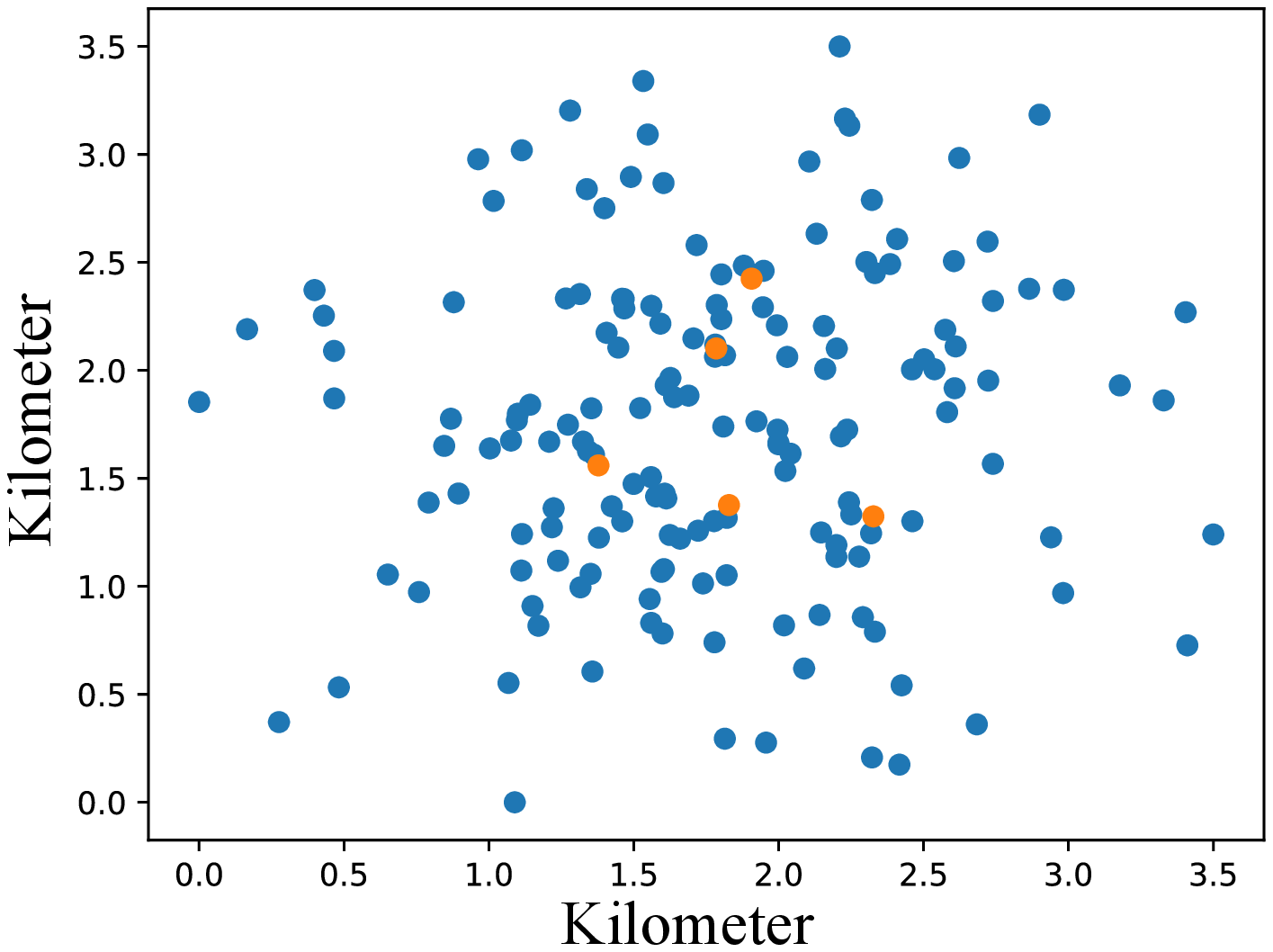}
		\caption[]%
		{distributed scenario}
		\label{fig:distributed_scenario}
	\end{subfigure}
	\caption{Illustration of spatial scenario variations}
	\label{fig:spatial_variation}
\end{figure*}

\paragraph{Spatial distribution}
The spatial variation is distinguished in clustered and distributed scenarios. The two scenarios allow us to study the effect of the demand pattern on the efficiency of the modular transport system. In the clustered scenarios the node locations are forming clusters following a Gaussian distribution using the depot locations as centers. The depots themselves are positioned following a Gaussian distribution on the area using the areas center as mean location. In Figure \ref{fig:clustered_scenario} an example clustered scenario is given. The x and y coordinates are given in kilometers. In distributed scenarios (see Figure \ref{fig:distributed_scenario}) both the depot and the node locations are spread over the area using a Gaussian distribution with the center of the area as the mean location. For passenger requests pickup and drop off nodes are spread over the case study area, whereas for freight requests the pickup location is at one of the depot locations and the drop off location is spread over the case study area. 

\begin{figure*}[ht]
	\centering
	\begin{subfigure}[b]{0.49\textwidth}
		\centering
		\includegraphics[width=\textwidth]{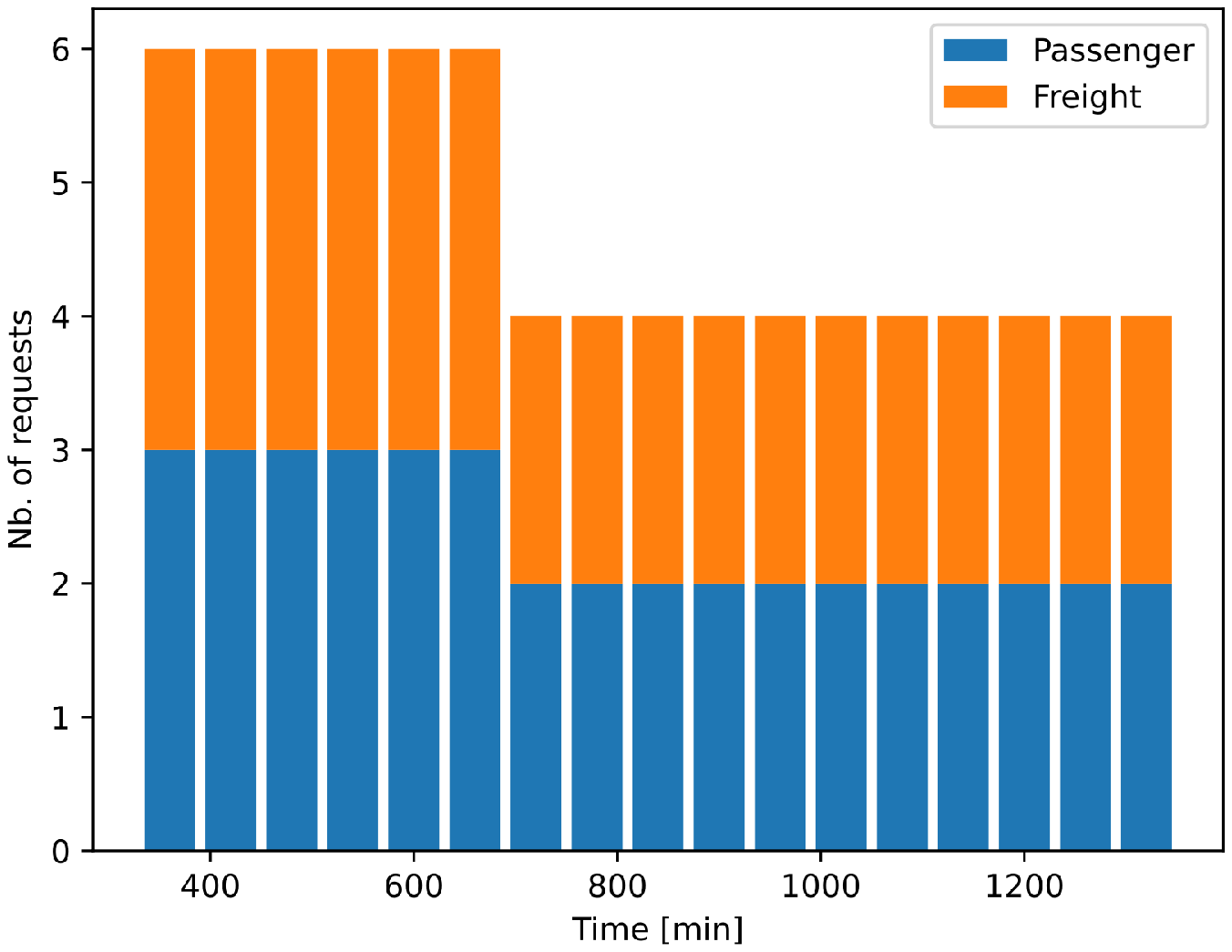}
		\caption[]%
		{even scenario}
		\label{fig:even_scenario}
	\end{subfigure}
	\hfill    
	\begin{subfigure}[b]{0.49\textwidth}
		\centering
		\includegraphics[width=\textwidth]{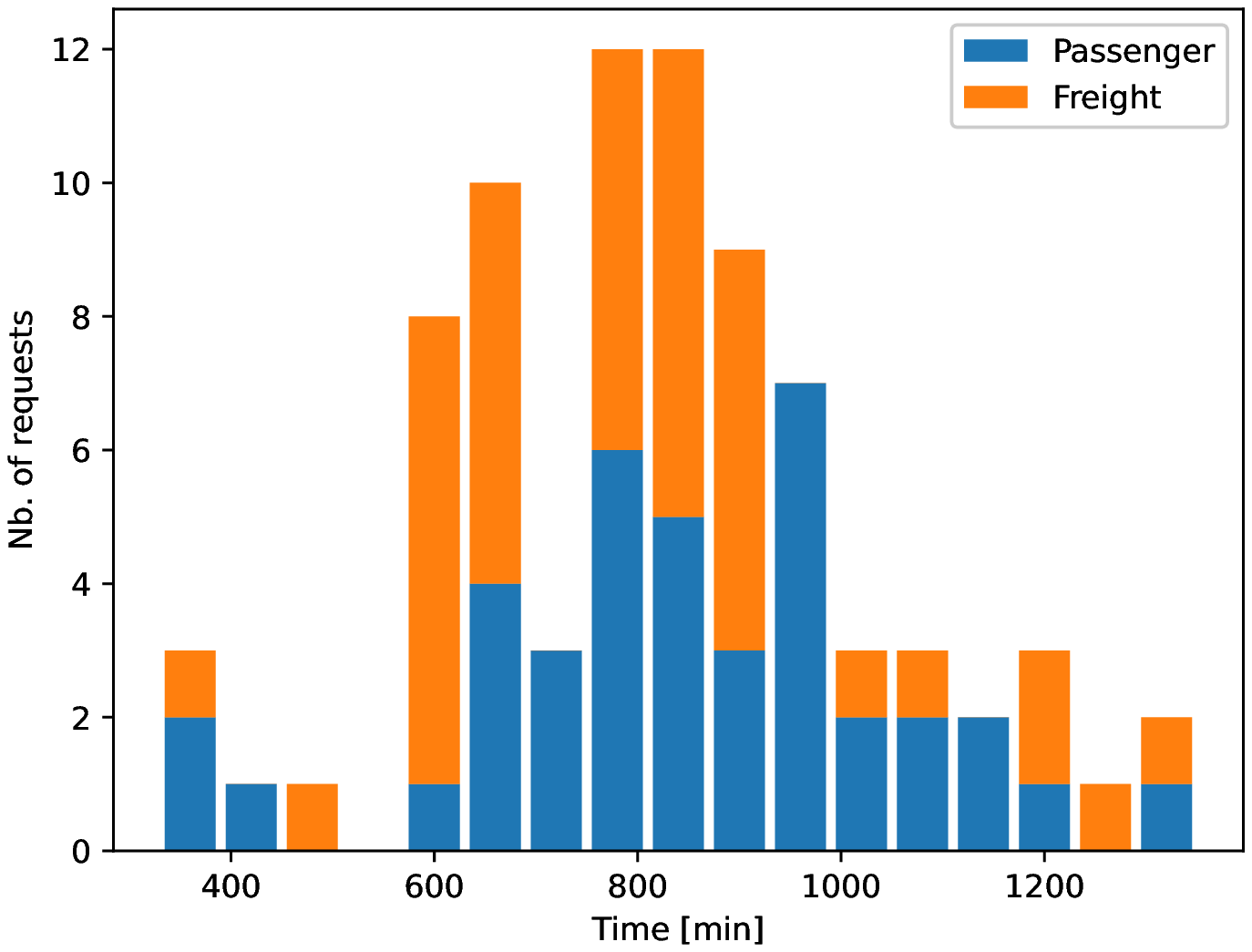}
		\caption[]%
		{peak scenario}
		\label{fig:normal_scenario}
	\end{subfigure}
	\caption{Illustration of temporal scenario variations}
	\label{fig:temporal_variation}
\end{figure*}

\paragraph{Temporal distribution}
The temporal demand variation is distinguished in evenly distributed and peak demand levels. The creation of even and peak demand levels allows to study the effect peak demand has on the proposed transportation system. In scenarios with even demand distribution, the requests are spread out evenly over the planning period. This is achieved by assigning evenly spread arrival time definitions for each passenger pickup node and freight drop off node. In Figure \ref{fig:even_scenario} an example distribution is shown. Here, the colors indicate the different demand types, the x-axis shows the time, and the y-axis shows the number of requests. 

 In the peak demand scenarios, the number of requests in each time period is computed so that it follows a demand curve with a peak during the middle of the time period. The total number of requests remains unchanged. For each passenger request the pickup and drop off time windows are designed based on the time distribution assumed. The pickup time windows ($i\in N^+_p$) are defined as $a_i \sim U(t - \Delta, t + \Delta)$ and $b_i = a_i + \Delta$, with $\Delta = \sim U(\SI{5}{min},\SI{20}{min})$. The drop off time windows for passenger requests ($j\in N^-_p$) are computed with $a_j \sim U(a_i - t_{i,j}, a_i + \SI{60}{min})$ and $b_j = a_j + \Delta$, with $\Delta = \sim U(\SI{5}{min},\SI{20}{min})$, the travel time $t_{i,j}$ between the corresponding pickup node $i$ and the drop off node $j$. For freight requests the pickup time windows range over the entire planning period, hence $a_i = \SI{360}{min} \quad \forall i \in N^+_f$ and $b_i = \SI{1320}{min} \quad \forall i \in N^+_f$. The drop off time windows are designed in the same way as the passenger pickup time windows, i.e. $a_i \sim U(t - \Delta, t + \Delta)$ and $b_i = a_i + \Delta$, with $\Delta = \sim U(\SI{5}{min},\SI{20}{min})$ and $t$ being the computed time based on the assumed temporal distribution. With these time window definitions, each request can be served by an individual platoon so that not all requests are infeasible by design. In Figure \ref{fig:normal_scenario} an example of the peak demand levels is shown. 

The conducted experiments contain a combination of the spatial and temporal distributions, hence there are four different scenario classes investigated, clustered-even, clustered-peak, distributed-even, and distributed-peak scenarios. In total, there are 5 different instances of each scenario type computed, hence the reported results are based on a total of 20 different scenarios.

\paragraph{Consolidation}
In order to study the impacts of consolidation of passenger and freight requests on the efficiency of the transportation system each scenario is solved once without consolidated optimization and once with consolidated optimization. If no consolidation is assumed the passenger and freight requests are solved separately and summed up. For conventional vehicle operation the solutions with and without consolidation are identical, since each vehicle module can only transport one type of demand. 

\newpage 

\subsection{Parameter settings}

The parameters used in the performed experiments are shown in Table \ref{tab:model_parameter}. If not specified further, these parameters are employed in all scenarios. The values for the number of requests, number of depots, service time and demand per request are set so that the computation time is within reasonable values and the complexity of the problems are comparable with real-sized problems.

\begin{table}
	\caption{Model parameter} 
	\label{tab:model_parameter}
    \centering
    \begin{tabular}{c|l|c}
        Notation                        & Description  & Value  \\ 
        \hline
        $n_r$ & number of requests & 80\\
    $n_d$ & number of depots & 5\\
    $R$ & maximum range for a platoon & \SI{200}{km}\\
    $l_{max}$ & maximum number of platoons & 20\\
    $Z_{max}$ & maximum number of modules in one platoon &10\\
    $Z_k$ & maximum number of modules per type $k \in K$ &10 per type\\
    $Q_k$ & module capacity for module type $k \in K$ &15 per type\\ 
    $\alpha_1$ & travel distance related cost parameter &\SI{0.096}{EUR/km}\\
    $\alpha_2$ & fleet size related cost parameter &\SI{19.37}{EUR/h}\\
    $\alpha_3$ & trip duration related cost parameter& \SI{6.9}{EUR/h}\\
    $\alpha_4$ & unserved requests related cost parameter& \SI{19.37}{EUR/h}\\     
    $\eta$ & operational cost reduction& 0.6    
    \end{tabular}
\end{table}

The values for module range ($R$), module capacity ($Q_k$) and maximum platoon length ($Z_{max}$) are in line with reported values from research reports and prototype vehicles \citep{nextHOME}. The number of modules per type ($Z_k$) is defined to provide sufficient supply to absorb the demand, hence it is feasible to serve every request with the provided combination of fleet size and module capacity. The maximum number of platoons is set to the sum of available modules per type to allow one trip per module. Hence it is feasible to operate the fleet of modules as conventional, individual vehicles. The cost parameters $\alpha_1$ and $\alpha_2$ are based on the values reported in \cite{militaoOptimalFleetSize2021}. In their work the authors estimate the distance-based and time-based operation costs of electric vehicles based on available data from the city of Munich, Germany. The authors propose a linear estimation model based on vehicle capacity to estimate the two cost parameters. On the basis of the published data, the models are $\alpha_1 = 0.003599 \cdot Q_k + 0.04162$ and $\alpha_2 = 0.1753 \cdot Q_k + 16.8$. The parameter $\alpha_2$ is scaled with the duration of the planning period. The parameter $\alpha_3$ scales travel time costs and is based on the value-of-time for public transport users in Stockholm, Sweden as reported by \cite{borjessonExperiencesSwedishValue2014}. In Equation \eqref{eq:objective} the parameter $\alpha_3$ corresponds to the value of time for passenger travel. Hence, it increases the overall cost of total travel time per route in the objective function. Since this value is also applied on freight requests, which can be assumed to have a lower value of time, this implies an overestimation of the total travel time costs. However, the different value of time for passenger and freight requests is implicitly considered using the time window definitions. In the scenario definition passenger time windows are set to be tighter compared to freight requests, indicating the higher importance of on-time operations for passengers. The parameter $\alpha_4$ penalizes the unserved requests. The value for this parameter is chosen to scale with the fleet size related cost parameter. The general reasoning behind the value is that it should not lead to a reduced objective value if a single module serves a single request. Hence, the penalty for keeping demand unserved equals that of $\alpha_2$ multiplied with the duration of the planning period and the average demand for the request. The value of $\eta$ is estimated to be 0.6 as it is reported in \cite{zhangEfficiencySemiAutonomousFully2019} who provide cost reduction estimations for bus platoon operations based on data published in project reports. 

The difference between the conventional vehicle operations and the proposed modular vehicle operations is modeled by adjusting the $Z_{max}$ parameter. For the case of conventional vehicle operations $Z_{max} = 1$ and for modular vehicle operations $Z_{max} = 10$. The other cost parameters remain unchanged to investigate the impact of the changed vehicle operations. When optimizing a scenario using conventional vehicle operations, the best solution for the corresponding scenario using modular vehicle operations with separated demand is used as the initial solution. That way a high quality initial solution is used which improves the optimization process by reducing computation time and achieving high robustness. 

\subsection{Parameter analysis}

A parameter analysis is performed in order to explore the influence different parameter settings have on the objective value and the benefits modular vehicles might have on the transportation system . The three parameters varied in this analysis are the vehicle capacity, vehicle range and the operational cost reduction parameter. For each of these parameters an array of discrete values is created. The vehicle capacity spans from 15 to 45, the vehicle range spans from \SI{50}{km} to \SI{250}{km} and the cost reduction parameter ($\eta$) spans from 0.2 to 1. All other parameters remain unchanged. The parameter study is performed on scenarios with distributed-peak demand characteristics. Each of the 5 created scenarios with this temporal and spatial demand definition is computed with all parameter settings. The reported values are the average values for all scenarios and its ensemble runs. In total $3 \cdot 5 \cdot 10 \cdot \left(5+5+4\right) = 2100$ experiments are evaluated for parameter analysis. 

\subsection{Key performance indicators}

In addition to the cost terms (route distance, route duration, fleet size, and unserved demand), the analysis includes the following key performance indicators.

\paragraph{Fill rate}
The fill rate is computed by adding the demand of all the pickup nodes on all the routes of a solution. This total demand units transported is divided by the total capacity of all modules of all platoons used in the solution.  Therefore, this indicator provides an understanding of the overall utilization of the supplied capacity. The higher the fill rate, the better the module capacity is used for a solution. 

\paragraph{Request kilometer}
The request kilometer is computed as the sum of all travelled distances over the requests (passenger and freight). The travelled distance is the distance from pickup node to drop off node of that request along the travelled route. 

\paragraph{Request minutes}
Similar to the request kilometer this indicator sums the travelled time from the pickup node to the drop off node of all passenger and freight requests. In this indicator, also potential waiting times along a requests travel path are included. 

\paragraph{Load per platoon distance}
This indicator is calculated by dividing the total number of passenger and freight served with the total distance of all platoons in that solution.

\paragraph{Empty vehicle kilometers}
The empty kilometers of a solution are computed as the sum of the length of all segments in each route where the platoon is driving without any load on board. This indicator also includes the empty kilometers at the beginning and end of a route before returning to the depot.

\section{Results}
\label{sec:results}

\subsection{ALNS performance analysis}
The validation of the implemented heuristic optimization algorithm is performed by solving small and medium-sized routing instances with the exact solver CPLEX and ALNS. The solutions of both solvers are then compared. The small- and medium-sized instances are created as described in Section \ref{sec:experimental_design}. In Table \ref{tab:ALNS_validation} the results of the benchmark tests are summarized. For each routing instance the number of nodes, number of depots, the objective values for CPLEX and and best objective value using ALNS, and the computation times for CPLEX and ALNS are given in Table \ref{tab:ALNS_validation}. The ALNS algorithm is computed 10 times for each routing instance. The standard deviation of the objective value and the mean ALNS computation time of all ensemble runs are reported in Table \ref{tab:ALNS_validation}. The column "Mean First Found" shows the mean iteration for all ensemble runs where the final objective value was evaluated for the first time, indicating the convergence speed of the ALNS. The reported computation times are recorded from the start of the ALNS algorithm until one convergence criteria is met. CPLEX is configured to have a maximum computation time of \SI{1}{h} per instance and the accepted optimally gap is set to 0.01\%.

\begin{table}
\centering
\caption{ALNS validation results }
	\label{tab:ALNS_validation}
\begin{tabular}{c||c|c|c|c|c|c|c|c|}
ID & Nodes & Depots & \begin{tabular}[c]{@{}c@{}}Obj. \\CPLEX\end{tabular} & \begin{tabular}[c]{@{}c@{}}Obj. \\ALNS\end{tabular} & \begin{tabular}[c]{@{}c@{}}Stdev. Obj. \\ALNS\end{tabular} & \begin{tabular}[c]{@{}c@{}}Time [sec]\\CPLEX\end{tabular} & \begin{tabular}[c]{@{}c@{}}Mean Time [sec]\\ALNS~~\end{tabular} & Mean First Found  \\ 
\hline
1  & 8     & 1      & 547.465                                              & 547.465                                             & 0                                                          & 10.3                                                      & 17.17                                                           & 7.3               \\
2  & 8     & 1      & 800.808                                              & 800.808                                             & 0                                                          & 15.6                                                      & 17.8                                                            & 10.7              \\
3  & 8     & 1      & 676.032                                              & 676.032                                             & 0                                                          & 16.6                                                      & 9.92                                                            & 8                 \\
4  & 8     & 1      & 674.499                                              & 674.451                                             & 0                                                          & 6                                                         & 14.38                                                           & 5.3               \\
5  & 8     & 1      & 797.586                                              & 797.586                                             & 0                                                          & 6                                                         & 12.48                                                           & 3.1               \\
6  & 8     & 2      & 641.608                                              & 641.608                                             & 0.66                                                       & 236.5                                                     & 15.77                                                           & 108               \\
7  & 8     & 2      & 572.472                                              & 572.472                                             & 0                                                          & 213.1                                                     & 12.62                                                           & 104               \\
8  & 8     & 2      & 765.994                                              & 765.994                                             & 0                                                          & 68.6                                                      & 7.91                                                            & 15.8              \\
9 & 8     & 2      & 1015.083                                             & 1015.083                                            & 0.06                                                       & 273.1                                                     & 9.45                                                            & 34.2              \\
10 & 8     & 2      & 768.393                                              & 768.393                                             & 0.03                                                       & 273                                                       & 7.31                                                            & 1271.6            \\
11 & 8     & 2      & 849.648                                              & 849.648                                             & 0                                                          & 15.5                                                      & 8.45                                                            & 8.4               \\
12 & 8     & 2      & 845.769                                              & 845.769                                             & 0                                                          & 156.5                                                     & 9.2                                                             & 8.8               \\
13 & 8     & 2      & 847.112                                              & 847.112                                             & 0                                                          & 38.7                                                      & 8.35                                                            & 8.3               \\
14 & 12    & 1      & 705.163                                              & 705.163                                             & 0                                                          & 149.4                                                     & 27.07                                                           & 125.7             \\
15 & 12    & 1      & 925.901                                              & 925.881                                             & 0                                                          & 23                                                        & 30.15                                                           & 11.6              \\
16 & 12    & 1      & 579.082                                              & 579.082                                             & 0                                                          & 186.6                                                     & 15.24                                                           & 29.3              \\
17 & 12    & 1      & 580.514                                              & 580.514                                             & 0                                                          & 121.8                                                     & 36.08                                                           & 113.6             \\
18 & 12    & 1      & 801.68                                               & 801.68                                              & 0                                                          & 21                                                        & 38.67                                                           & 12.3              \\
19 & 12    & 2      & 654.029                                              & 654.029                                             & 0.09                                                       & 1083.8                                                    & 32.53                                                           & 801.4             \\
20 & 12    & 2      & 704.316                                              & 580.303                                             & 0                                                          & 3600.2                                                    & 17.18                                                           & 83.9              \\
21 & 12    & 2      & 1020.164                                             & 1020.164                                            & 1.79                                                       & 3601.2                                                    & 12.19                                                           & 430.2             \\
22 & 12    & 2      & 967.319                                              & 967.319                                             & 0.14                                                       & 668.9                                                     & 24.29                                                           & 21.4              \\
23 & 12    & 2      & 649.579                                              & 649.579                                             & 0.23                                                       & 3604                                                      & 11.07                                                           & 784.4             \\
24 & 12    & 2      & 701.256                                              & 701.256                                             & 0                                                          & 3601.5                                                    & 13.74                                                           & 149.5             \\
25 & 12    & 2      & 907.161                                              & 907.161                                             & 0                                                          & 818.8                                                     & 12.23                                                           & 52.6              \\
26 & 12    & 2      & 781.22                                               & 781.22                                              & 0                                                          & 229.8                                                     & 12                                                              & 36.9              \\
27 & 12    & 2      & 1099.148                                             & 1099.148                                            & 38.99                                                      & 3600.2                                                    & 21.55                                                           & 408.6             \\
28 & 12    & 2      & 676.335                                              & 676.334                                             & 0                                                          & 260.1                                                     & 28.49                                                           & 60.1              \\
29 & 16    & 1      & 832.738                                              & 708.393                                             & 0                                                          & 3600.4                                                    & 40.15                                                           & 54.9              \\
30 & 16    & 1      & 708.483                                              & 708.482                                             & 0                                                          & 56.8                                                      & 55.51                                                           & 65.7             
\end{tabular}
\end{table}

The ALNS algorithm is able to compute the optimal solution determined by CPLEX in all instances. In four instances (ID 4, ID 15, ID 28, ID 30), the objective value determined by the ALNS is slightly lower than the CPLEX value. The difference between the two values lies within the predefined optimality gap in CPLEX of 0.01\%, hence these instances can be assumed as solved to optimality as well. Additionally, the ALNS ensemble runs for each instance have a small standard deviation, indicating a robust convergence. In 22 out of the 30 instances, all 10 ensemble runs per instance converge to the same optimal solution. Even for some smaller instances (ID 8, ID 12), ALNS outperforms CPLEX in computation time. This finding is emphasized by the medium sized instances of 12 and 16 nodes where the ALNS converges faster to all optimal values. The computation time of ALNS is more related to the size of the problem, whereas the computation time of CPLEX is fluctuating within one problem size (see instances 29 and 30). In the instances 20 and 29 ALNS outperforms the CPLEX algorithm in both computation time and objective value by a large margin. In these instances CPLEX does not converge to a global optimal solution within 1h. The last column in Table \ref{tab:ALNS_validation} shows the mean iteration over all ensemble runs in which the final optimal solution was evaluated first. This metric shows the fast convergence speed of the ALNS algorithm and motivates the chosen convergence criteria. Using this metric and the average computation time per iteration, the ALNS finds the optimal solution for all benchmark instances within \SI{0.5}{sec}. 

\subsection{Scenario analysis}
In this section the results for the different scenarios in terms of spatial and temporal demand distributions are presented. The numerical values are the average values over all computed scenarios and their ensemble runs. The analysis focuses on the effects of vehicle technology and consolidation on system performance. 

\begin{figure*}[ht]
	\centering
	\includegraphics[width=\textwidth]{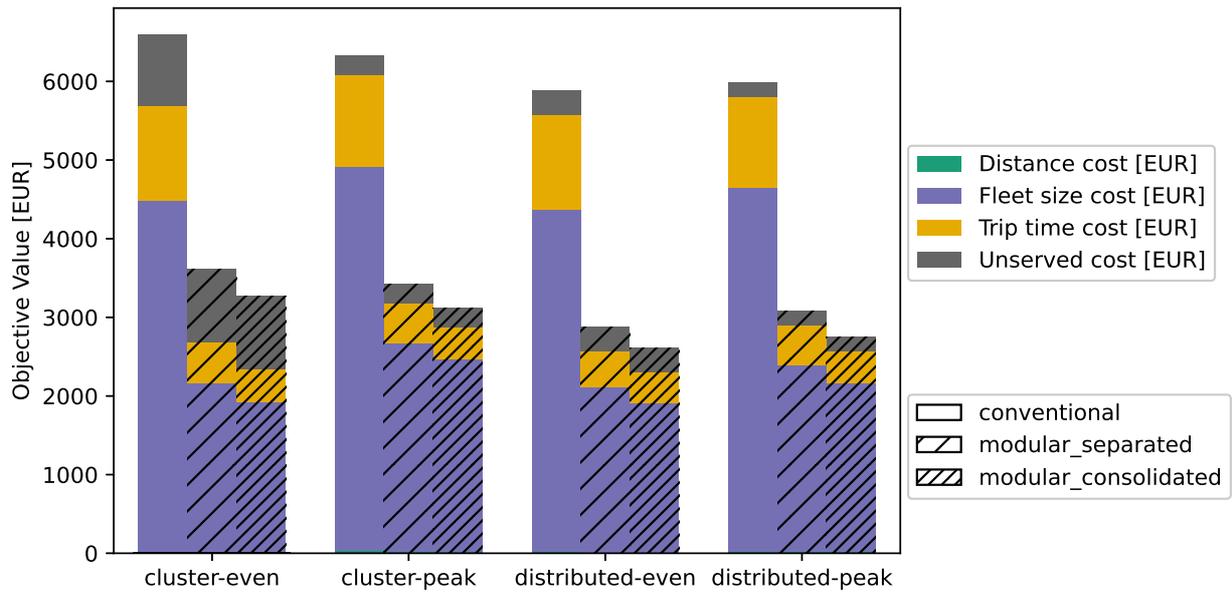}
	\caption{Objective function values per scenario}
	\label{fig:scenario_study_Objective_Value}
\end{figure*}

Figure \ref{fig:scenario_study_Objective_Value} shows the objective values with each color representing a cost term. It can be seen that the cost of the fleet size and the cost of the trip time dominate the objective value. The distance-related cost is negligible which can be explained by the parameter settings compared to the higher unit costs for modules. When analyzing the impacts of vehicle technology and consolidation, it can be seen that, due to a large reduction in fleet size and trip time costs, modular operations result in lower overall costs, while the unserved demand remains unchanged. Additional savings can be achieved by introducing consolidation between passenger and freight demand. These cost savings are independent of the spatial and temporal demand distribution. In addition to the effect of vehicle technology and consolidation, the results show that spatially distributed demand patterns lead to slightly lower total costs. In numerical values, the use of modular vehicles results in an average reduction of 48\%, and the consolidated optimization further reduces the objective value by 9\%. The fleet size costs, distance costs, and trip time costs can be reduced by 50\%, 15\%, and 58\% for modular vehicles and by further 9\%, 4\%, and 18\% for consolidation, respectively. The costs for unserved demand remain unchanged. The reason why some scenarios have unserved requests is the time windows for these requests. It would require one one-module platoon per request to serve these requests. Since the penalty cost for one unserved request is as high as the cost of one platoon with one module, serving a single request with one module is not beneficial for the overall system due to the additional travel time/distance. By design (fleet size costs, penalty for unserved demand and travel distance/time cost), serving a single request with a single module results in a higher objective value compared to not serving it. 

\begin{figure*}[ht]
	\centering
	\begin{subfigure}[b]{0.49\textwidth}
		\centering
    	\includegraphics[width=\textwidth]{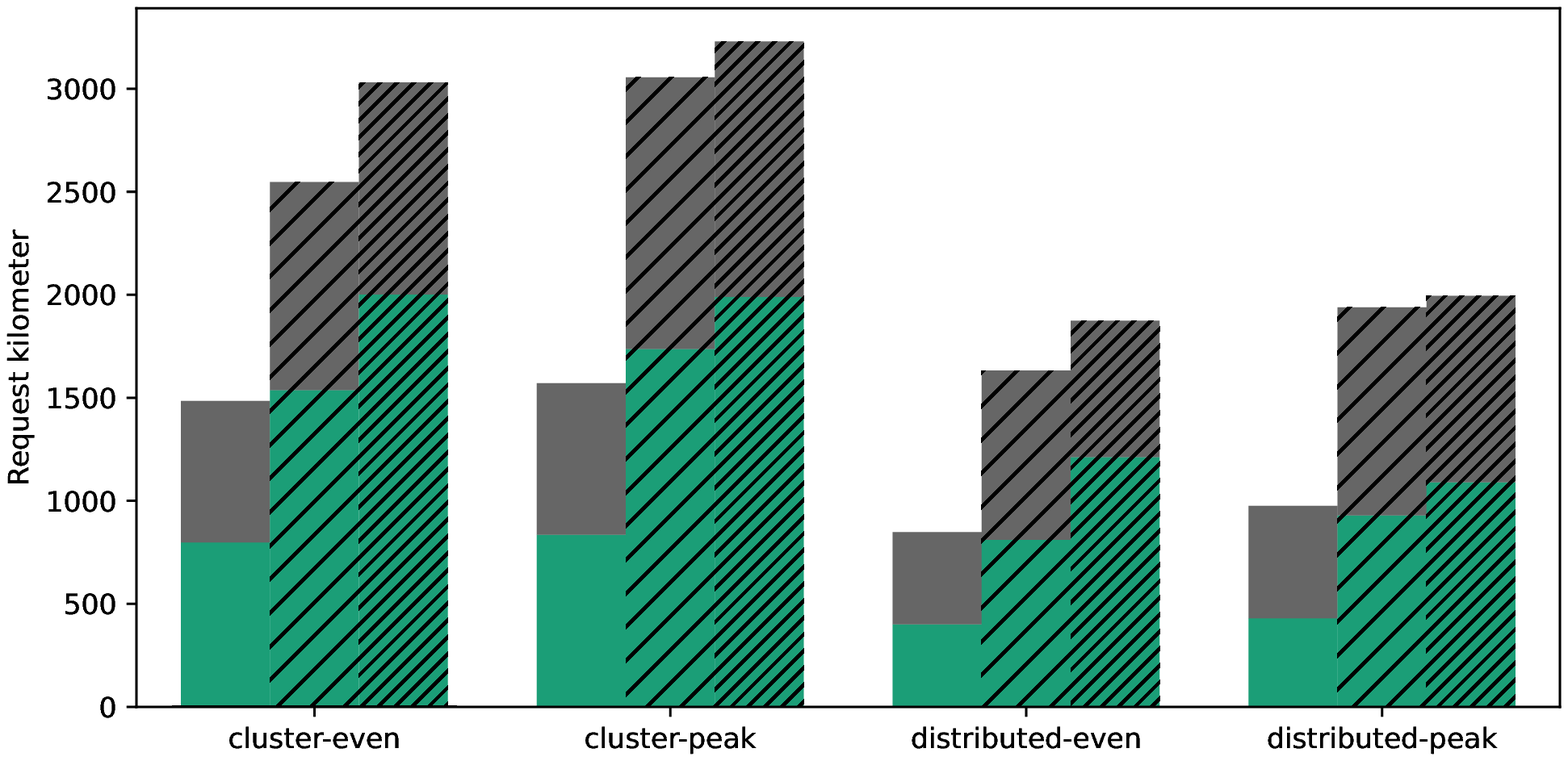}
    	\caption{Effect on Request kilometer}
    	\label{fig:scenario_study_Request_kilometer}
	\end{subfigure}
	\hfill    
	\begin{subfigure}[b]{0.49\textwidth}
		\centering
    	\includegraphics[width=\textwidth]{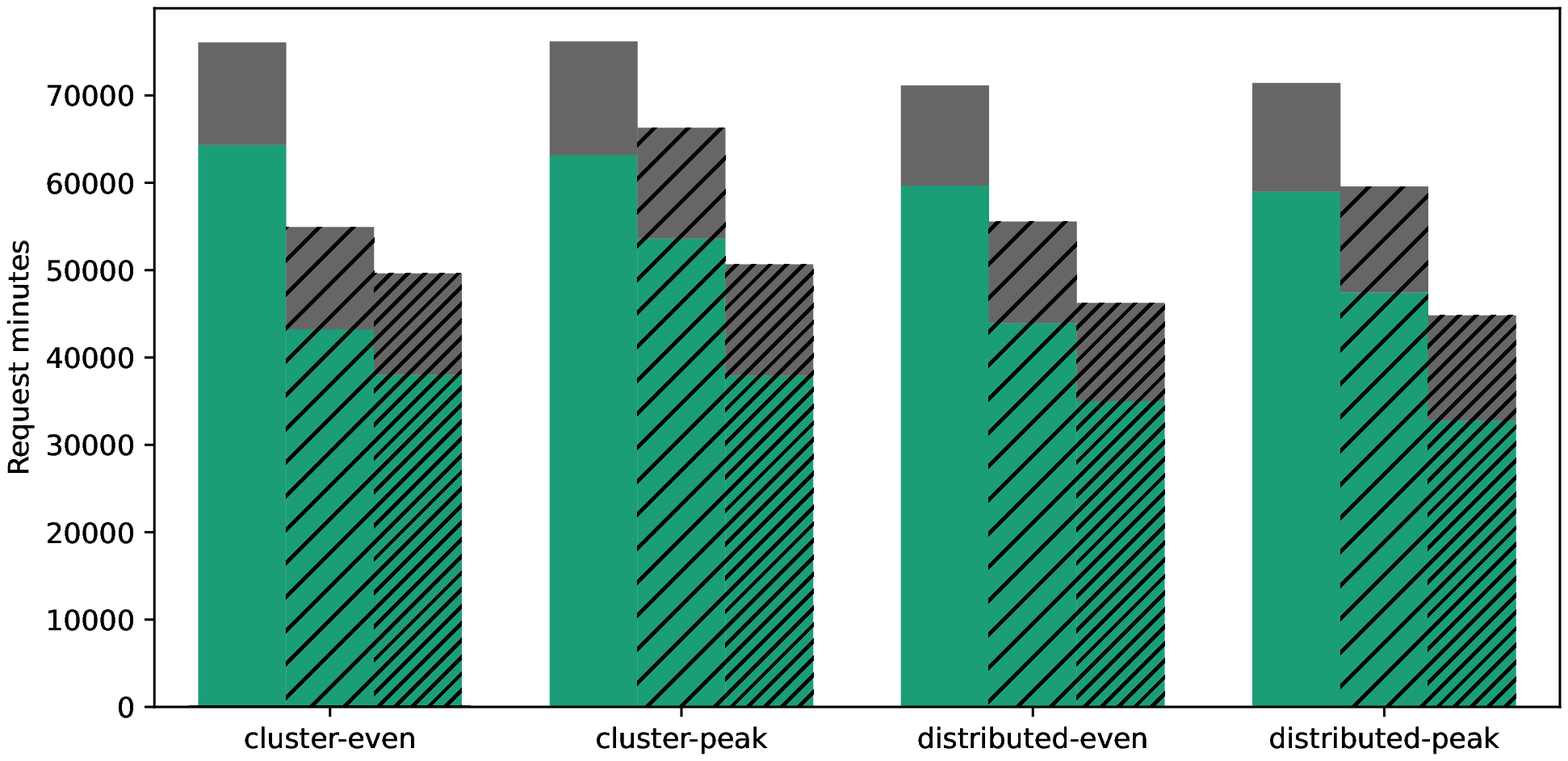}
    	\caption{Effect on Request minutes}
    	\label{fig:scenario_study_Request_minutes}
	\end{subfigure}
	\hspace*{0.5\fill}
	\begin{subfigure}[b]{0.49\textwidth}
		\centering
    	\includegraphics[width=\textwidth]{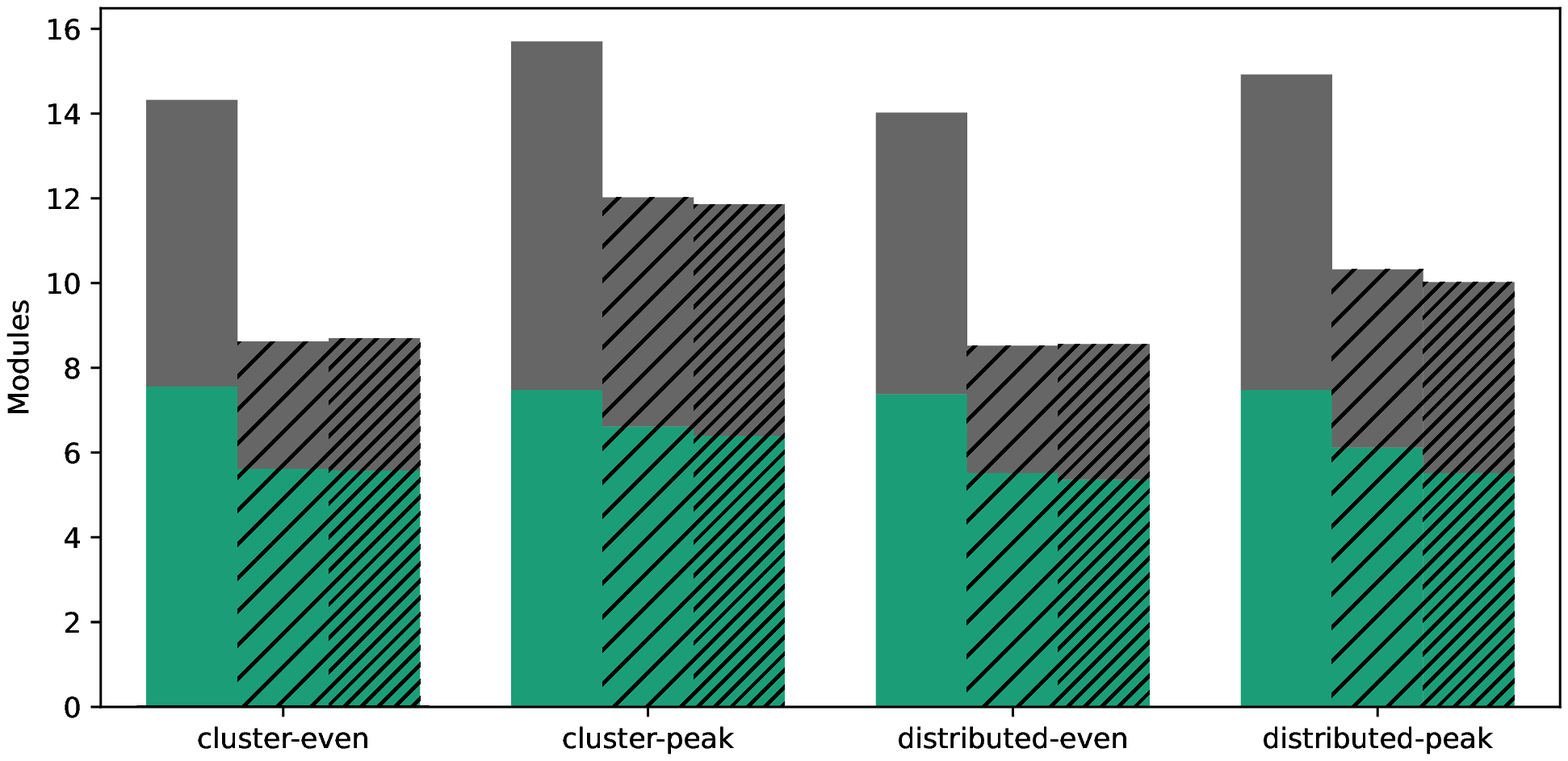}
    	\caption{Effect on Number of modules}
    	\label{fig:scenario_study_Modules}
	\end{subfigure}
	\caption{Effects of spatial and temporal demand distributions. The vehicle operations are indicated by different patterns, conventional (\protect\includegraphics[height=1.5ex]{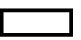}), modular without consolidation (\protect\includegraphics[height=1.5ex]{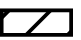}), and modular with consolidation (\protect\includegraphics[height=1.5ex]{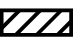}). The colors indicate the demand type, passenger demand (\protect\includegraphics[height=1.5ex]{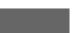}) and freight demand (\protect\includegraphics[height=1.5ex]{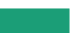}).}
	\label{fig:scenario_analysis_1}
\end{figure*}

\begin{figure*}[ht]
	\centering
	\begin{subfigure}[b]{0.49\textwidth}
		\centering
    	\includegraphics[width=\textwidth]{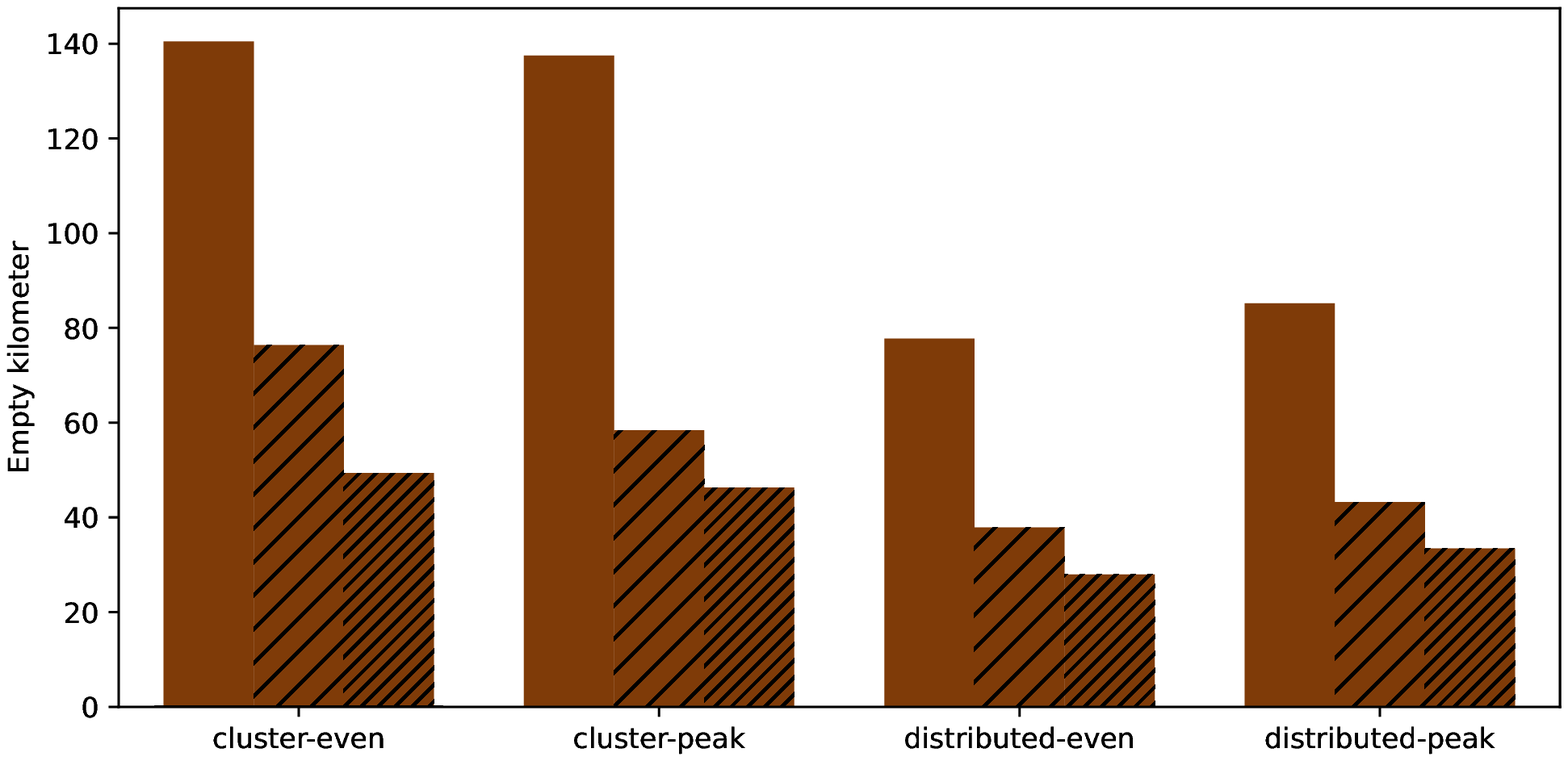}
    	\caption{Effect on Empty kilometer}
    	\label{fig:scenario_study_Empty_kilometer}
	\end{subfigure}
	\hfill
	\begin{subfigure}[b]{0.49\textwidth}
		\centering
    	\includegraphics[width=\textwidth]{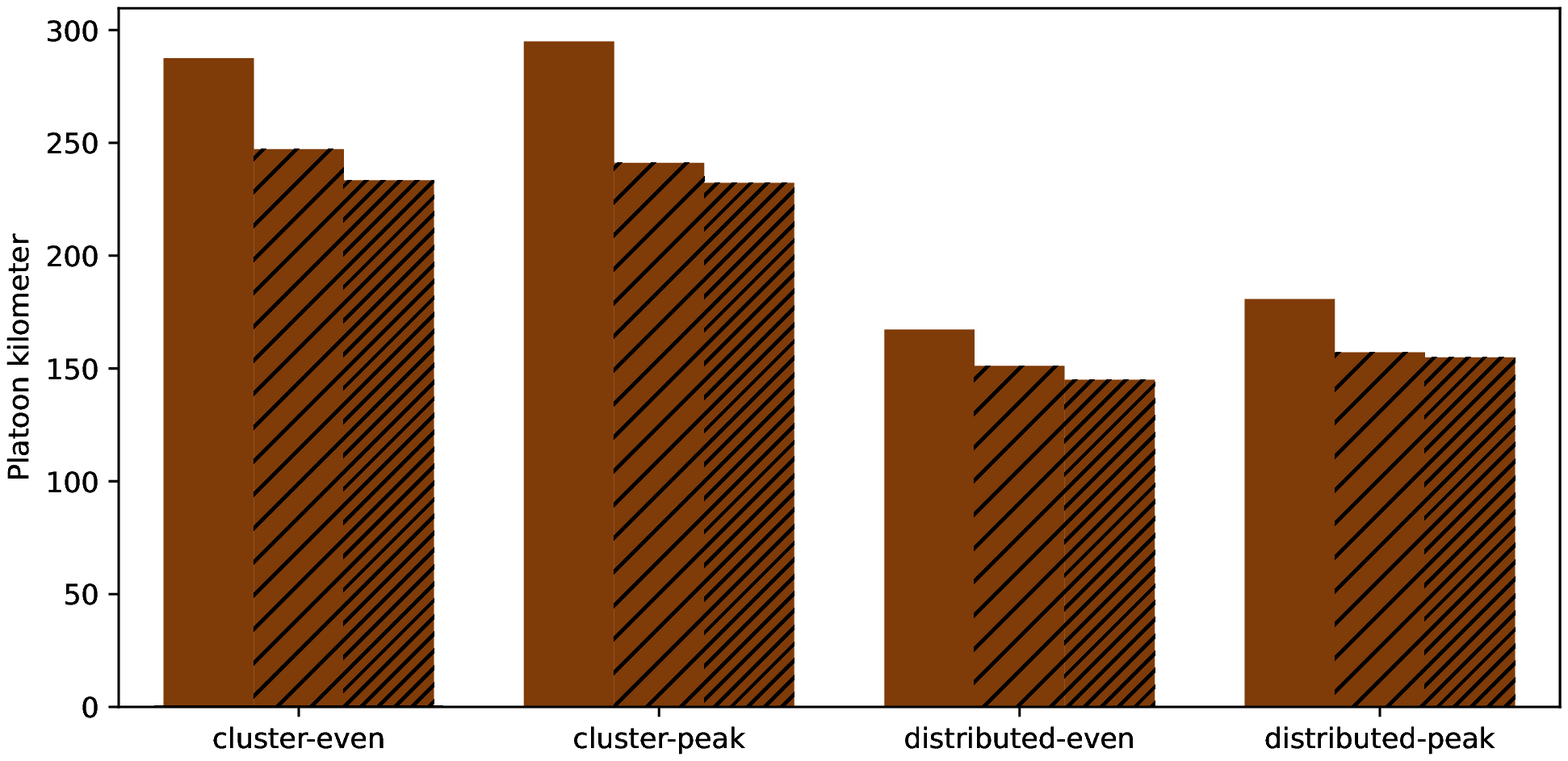}
    	\caption{Effect on Platoon kilometer}
    	\label{fig:scenario_study_Platoon_kilometer}
	\end{subfigure}
	\hspace*{0.5\fill}
	\begin{subfigure}[b]{0.49\textwidth}
		\centering
	    \includegraphics[width=\textwidth]{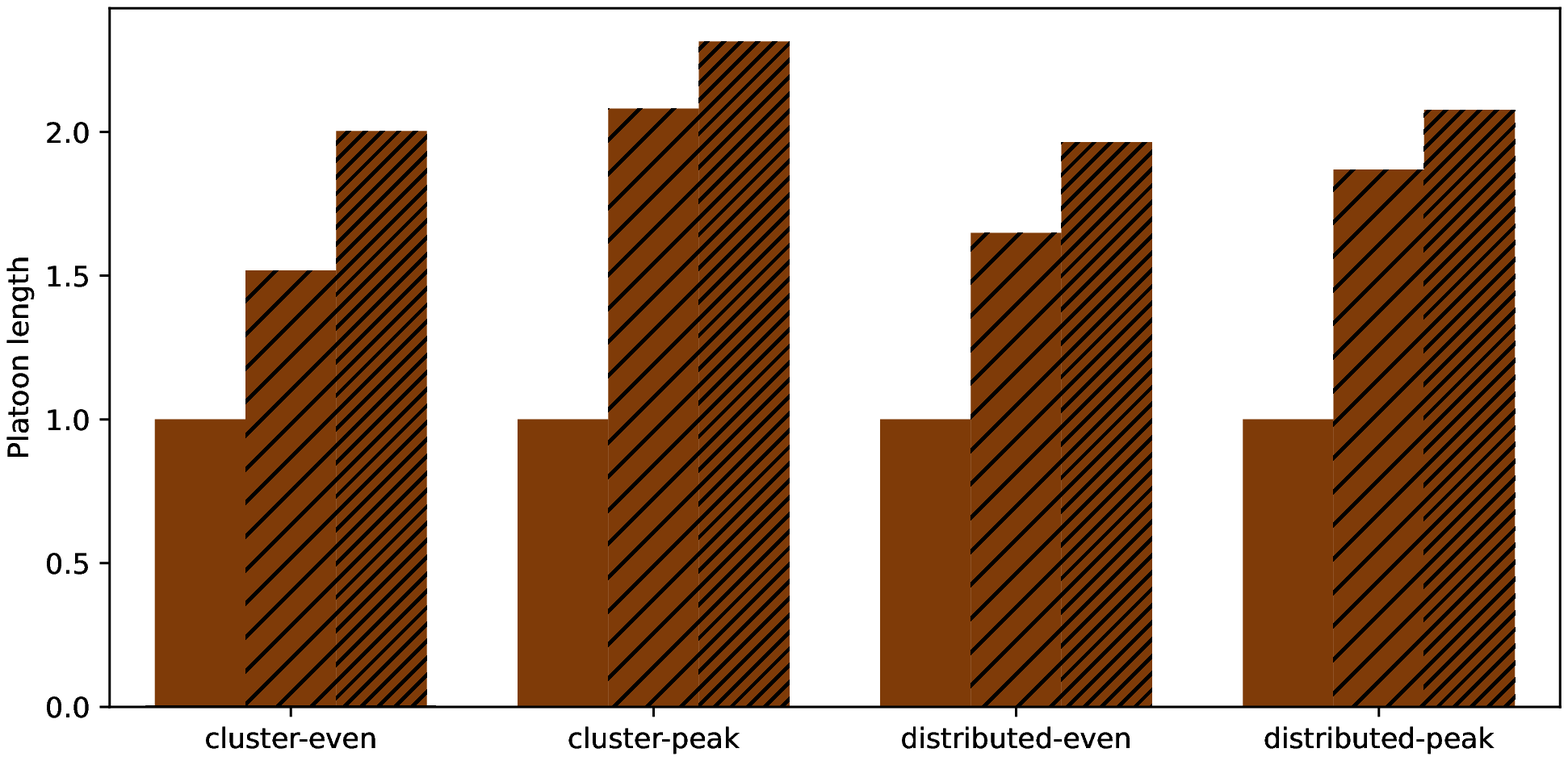}
	    \caption{Effect on Platoon length}
	    \label{fig:scenario_study_Platoon_length}
	\end{subfigure}
	\caption{Effects of spatial and temporal demand distributions. The vehicle operations are indicated by different patterns, conventional (\protect\includegraphics[height=1.5ex]{images/conventional.eps}), modular without consolidation (\protect\includegraphics[height=1.5ex]{images/modular_separated.eps}), and modular with consolidation (\protect\includegraphics[height=1.5ex]{images/modular_collaboration.eps}).}
	\label{fig:scenario_analysis_2}
\end{figure*}

The large reduction in the cost of fleet size is mainly due to the reduction in operational costs made possible by the formation of platoons. The creation of platoons and their average length can be seen in Figure \ref{fig:scenario_study_Platoon_length}, indicating that on average 1.8 modules form a platoon. The average platoon length increases to 2.1 if passenger and freight demand are jointly optimized. The large reduction in total fleet size cost is also due to a general reduction in the number of required modules for the different types of operation. Figure \ref{fig:scenario_study_Modules} shows the average number of modules used for each scenario. Additionally, in the conventional operation approximately the same number of modules is used for the transportation of passengers and freight, and the reduction in the number of modules in the other cases mainly stems from the reduction of required passenger modules. This observation is due to a change in the delivery order for the freight requests served. Due to the connection of multiple modules, more passengers can be jointly transported, allowing one to meet the time window constraints more efficiently. Such exploitation might result in longer travel distances which are compensated by the reduction in the required number of modules. The longer travel distances for the various types of operation can be seen in Figure \ref{fig:scenario_study_Request_kilometer}, where the kilometers traveled for each request demand unit increases with modularity and consolidation. It can also be seen that distributed demand patterns need fewer traveled kilometers to serve the demand compared to clustered demand patterns. This is expected since the average distance between nodes is lower for distributed scenarios. 

When looking at the empty vehicle kilometers for each scenario (see Figure \ref{fig:scenario_study_Empty_kilometer}) two clear trends can be seen. First, the introduction of modular vehicle types reduces empty kilometers by more than 50\% in all cases, while an additional 27\% can be saved by combining passenger and freight demand. This observation is another explanation for the reduction in the number of modules required. In distributed demand scenarios the empty kilometers are lower than in the clustered scenarios, hence a higher utilization of the supply capacity can be reached. The difference between even and peak demand patterns is not significant.

Figure \ref{fig:scenario_study_Platoon_kilometer} underscores the previous results. Although kilometers per request (see Figure \ref{fig:scenario_study_Request_kilometer}) increase by 90\% when modular vehicles are operated and increase by an additional 10\% if the demand is consolidated, the reduction of required modules leads to a general reduction of required platoon kilometers to serve the demand. This observation directly translates to a reduction of emissions and traffic in the area. This positive effect is most prominent for clustered scenarios and scenarios with peak demand. 

In contrast to the increase in request kilometers, request minutes are slightly reduced for the different vehicle operations. For modular operations without consolidation the request minutes are reduced by an average of 20\%, and the consideration of consolidation reduced the request minutes by another 19\%. Most of this reduction comes from a reduction in freight request minutes, which are reduced by 23\% for each operation mode respectively. This is again indicating a different delivery order for this type of demand. Passenger request minutes vary by less than 1\%, resulting in an unchanged level of service for users of such transportation systems. The discrepancy between freight and passenger request minutes can be explained by the design of the scenario. Freight requests have looser time window definition for the pickup nodes, meaning that longer waiting times between multiple requests, and there for longer in-vehicle times are feasible. The average in-vehicle time per freight request is approximately \SI{91}{min}, whereas for passenger requests the average in-vehicle time is \SI{17}{min} The passenger demand time windows are defined tightly and therefore do not allow for large pickup variations.

\subsection{Sensitivity analysis}

\begin{figure*}[ht]
	\centering
	\includegraphics[width=\textwidth]{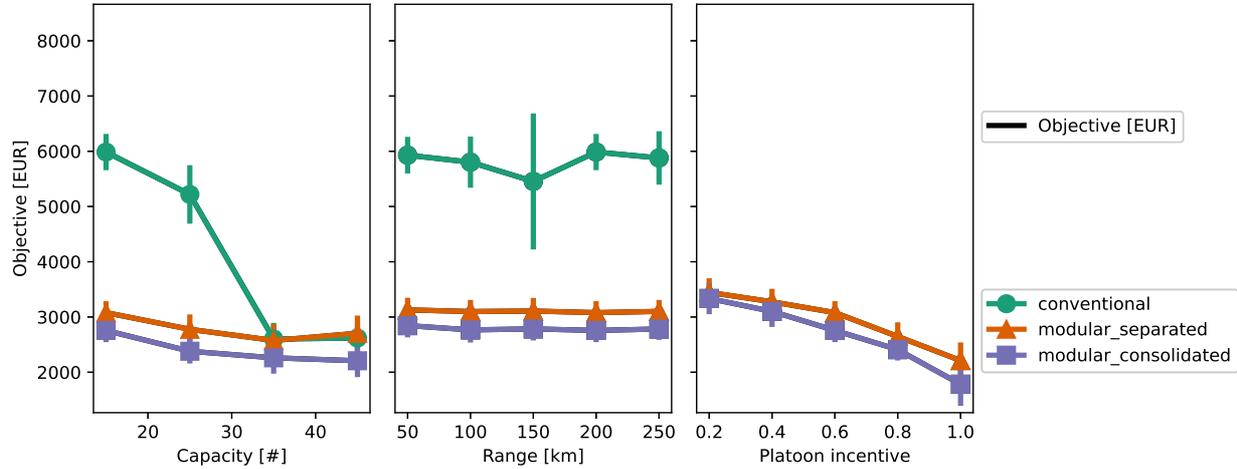}
	\caption{Effect of different parameter settings in the objective value.}
	\label{fig:parameter_study_Objective_}
\end{figure*}

\begin{figure*}[!ht]
	\centering
	\begin{subfigure}[b]{0.49\textwidth}
		\centering
    	\includegraphics[width=\textwidth]{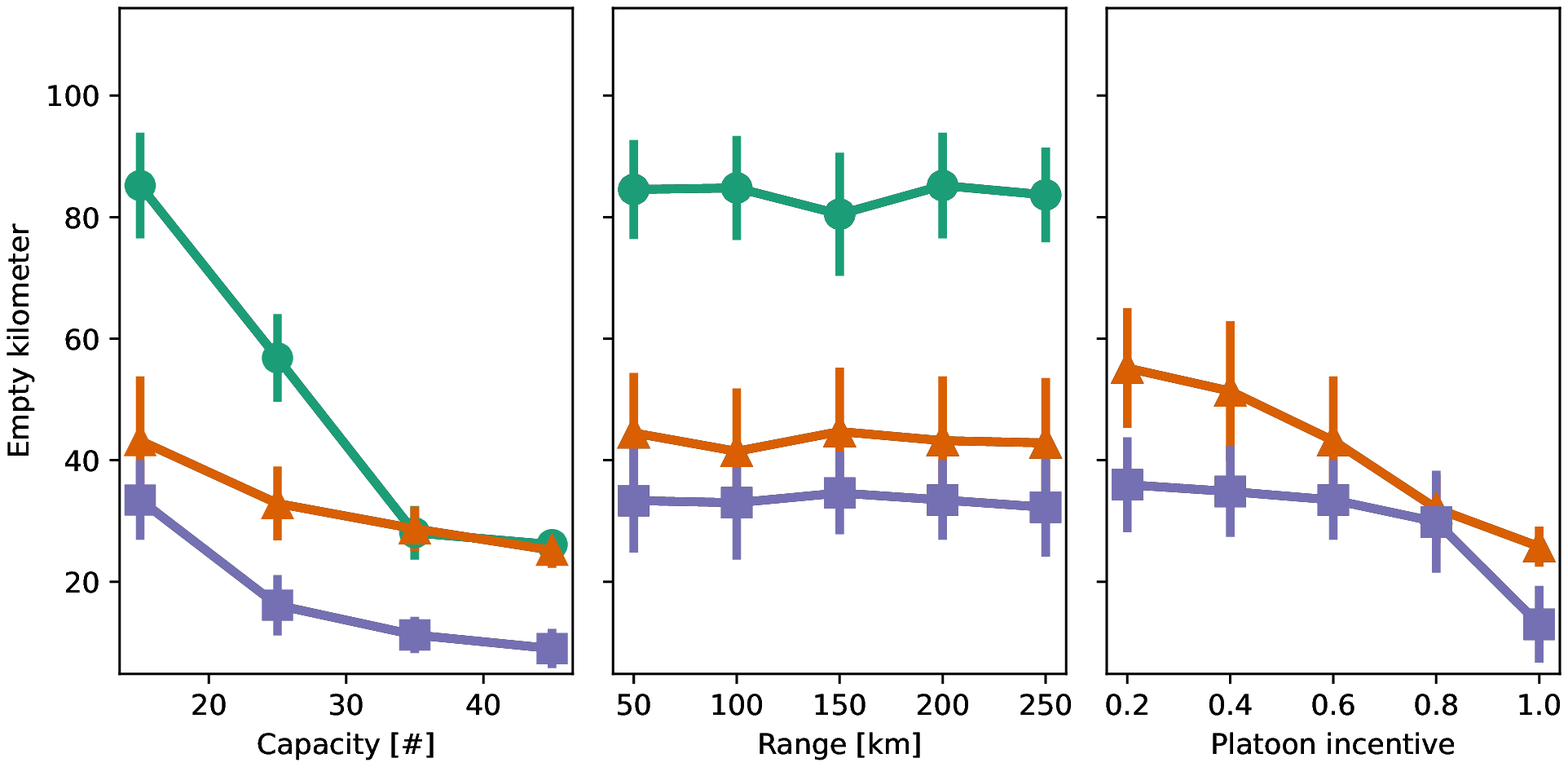}
    	\caption{Impact on Empty kilometer}
    	\label{fig:parameter_study_Empty_kilometer}
	\end{subfigure}
	\hfill    
	\begin{subfigure}[b]{0.49\textwidth}
		\centering
    	\includegraphics[width=\textwidth]{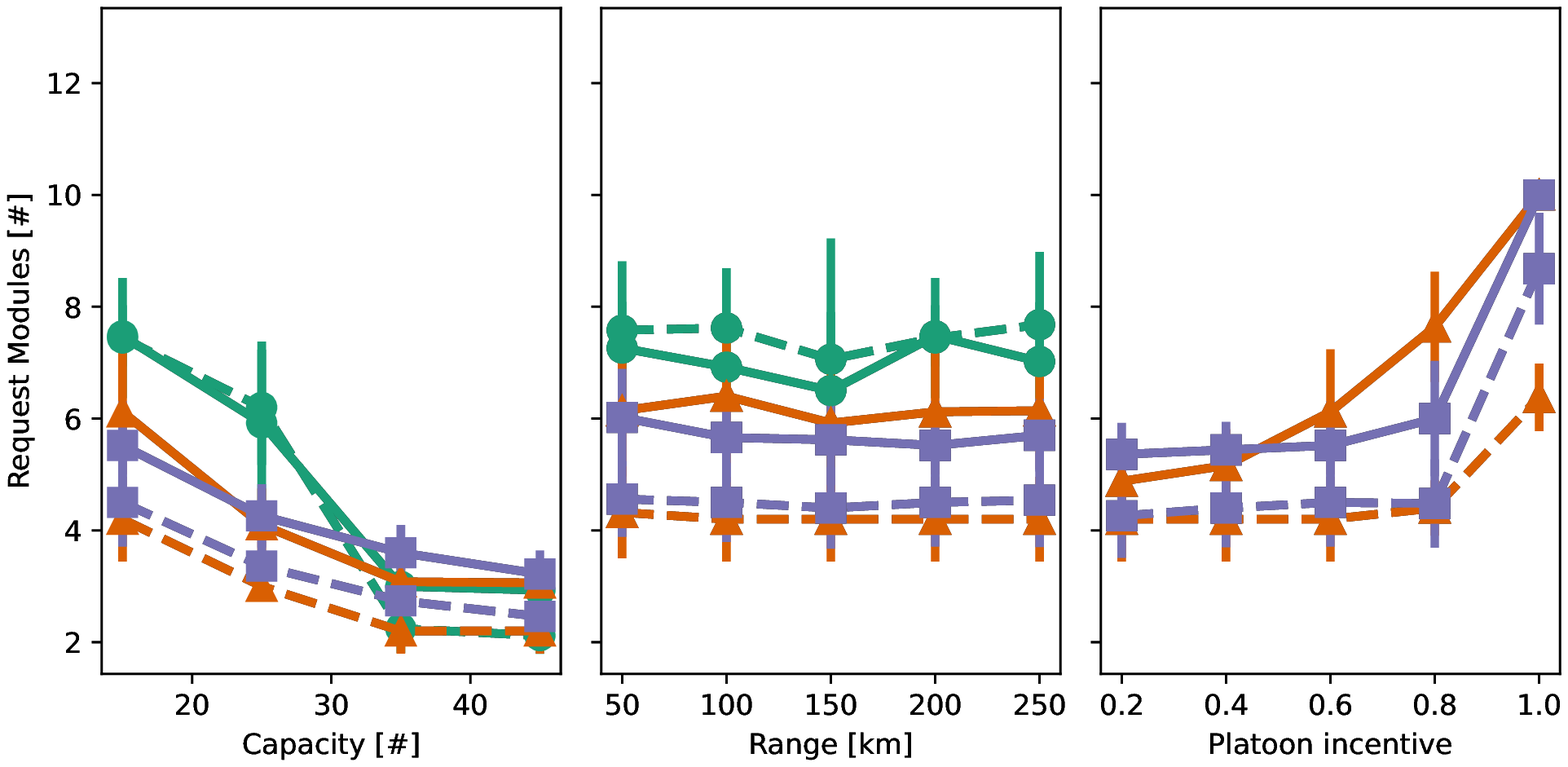}
    	\caption{Impact on Number of modules}
    	\label{fig:unit_parameter_study_Request_Modules_}
	\end{subfigure}
	\hspace{0pt}    
	\begin{subfigure}[b]{0.49\textwidth}
		\centering
    	\includegraphics[width=\textwidth]{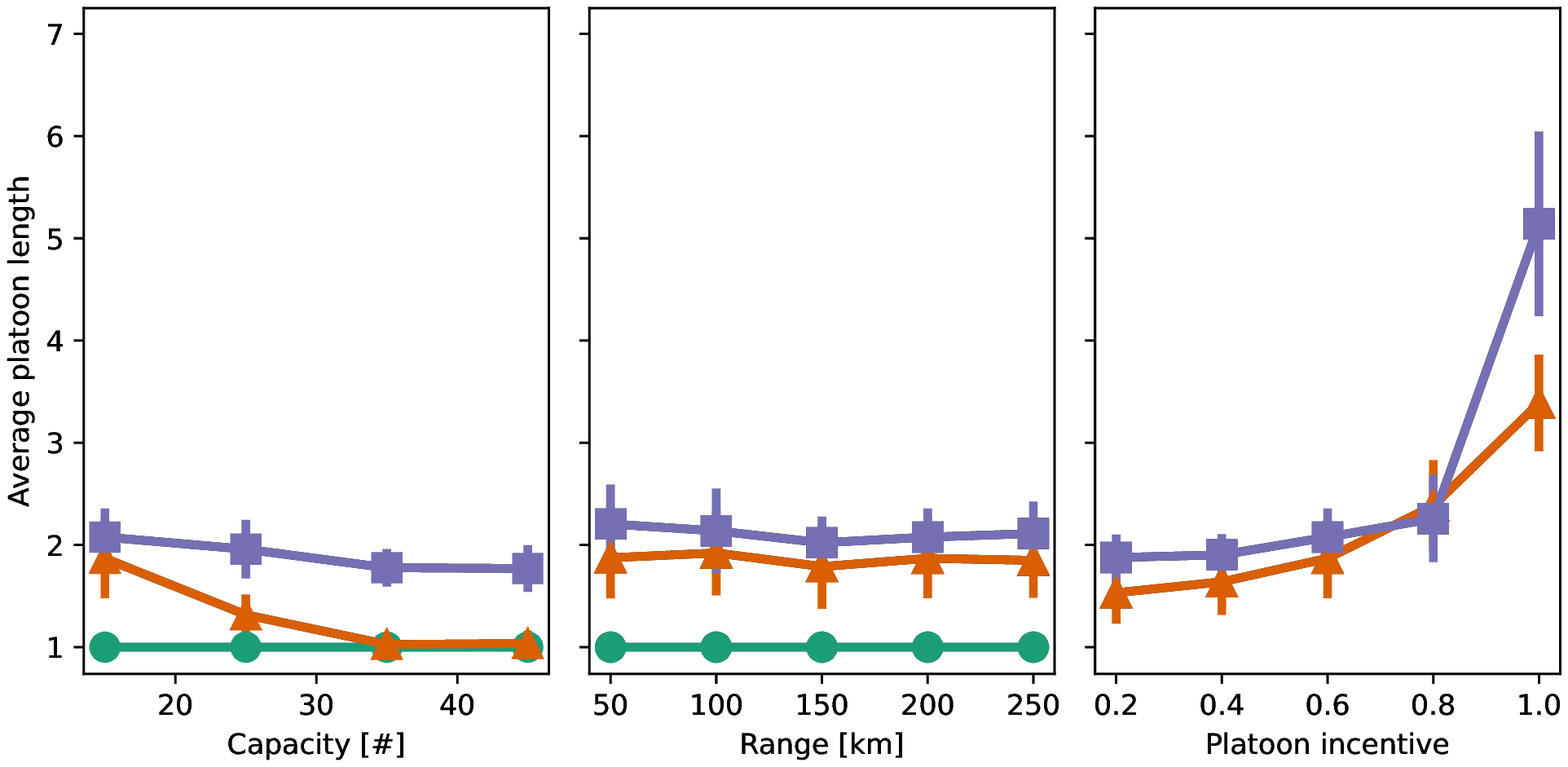}
    	\caption{Impact on Average train length}
    	\label{fig:parameter_study_Average_train_length}
	\end{subfigure}
	\hfill    
	\begin{subfigure}[b]{0.49\textwidth}
		\centering
    	\includegraphics[width=\textwidth]{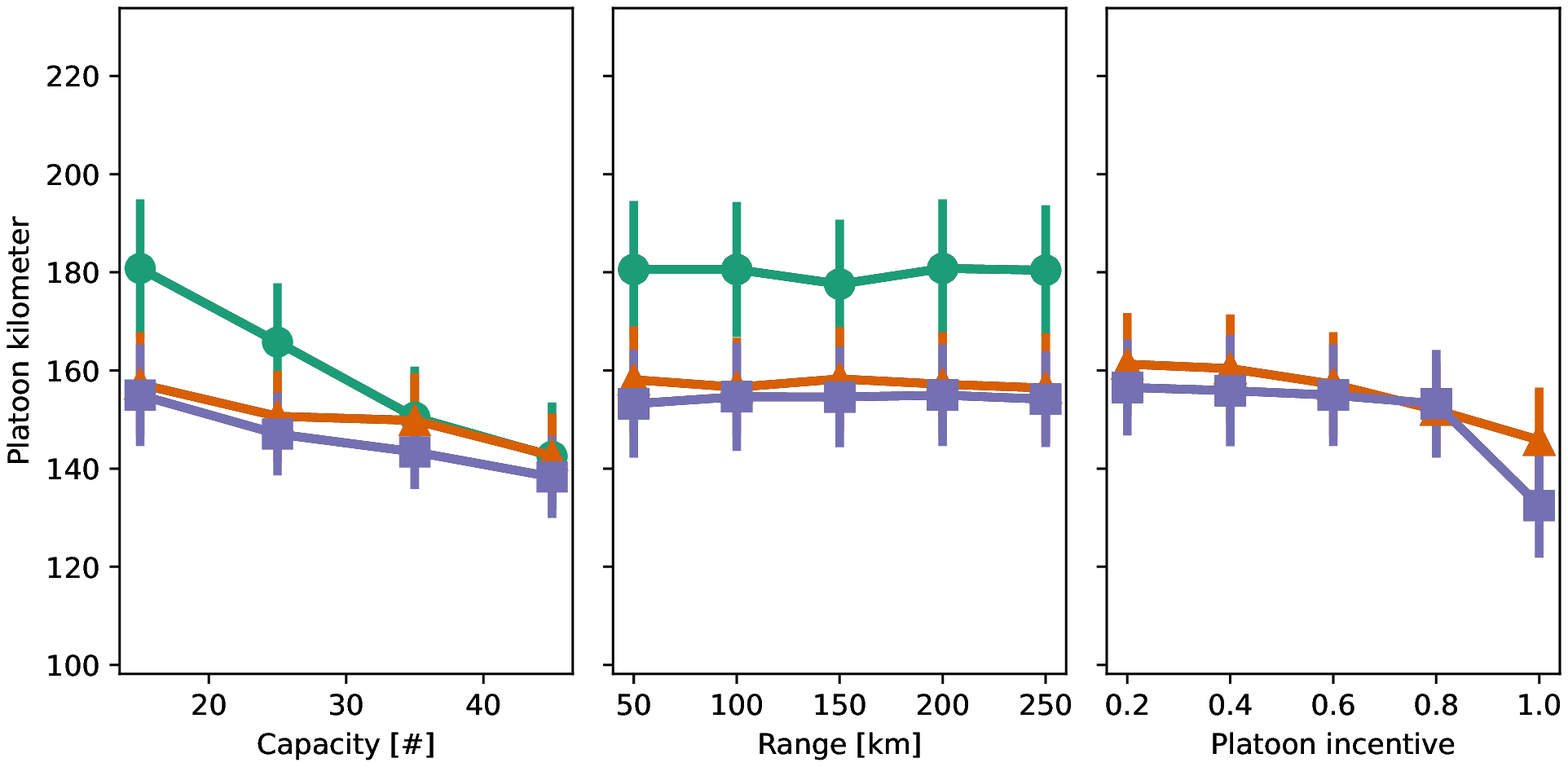}
    	\caption{Impact on Platoon kilometer}
    	\label{fig:parameter_study_Platoon_kilometer}
	\end{subfigure}
	\hspace{0pt}    
	\begin{subfigure}[b]{0.49\textwidth}
		\centering
    	\includegraphics[width=\textwidth]{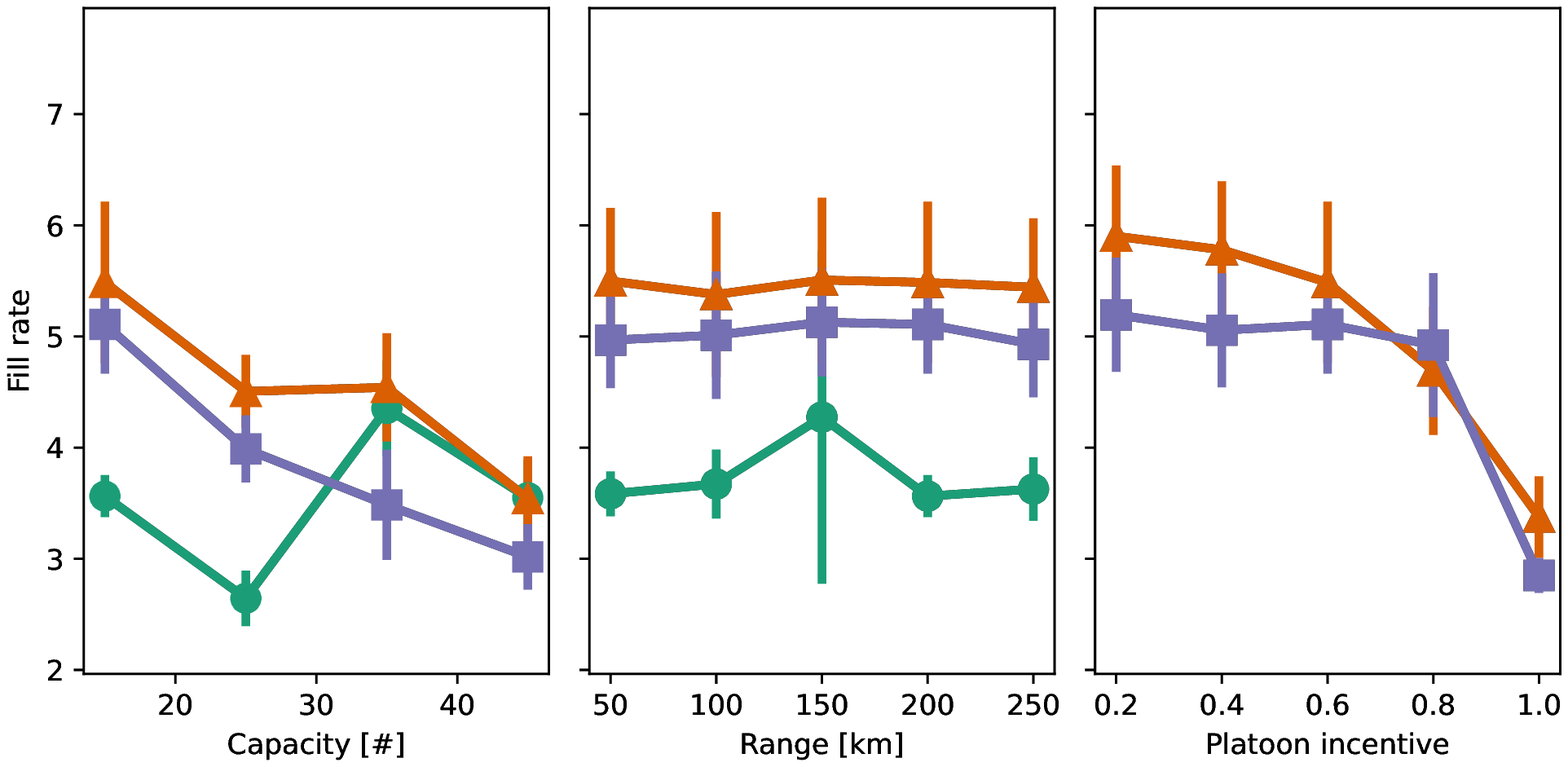}
    	\caption{Impact on Fill rate}
    	\label{fig:parameter_study_Fill_rate}
	\end{subfigure}
	\hfill    
	\begin{subfigure}[b]{0.49\textwidth}
		\centering
    	\includegraphics[width=\textwidth]{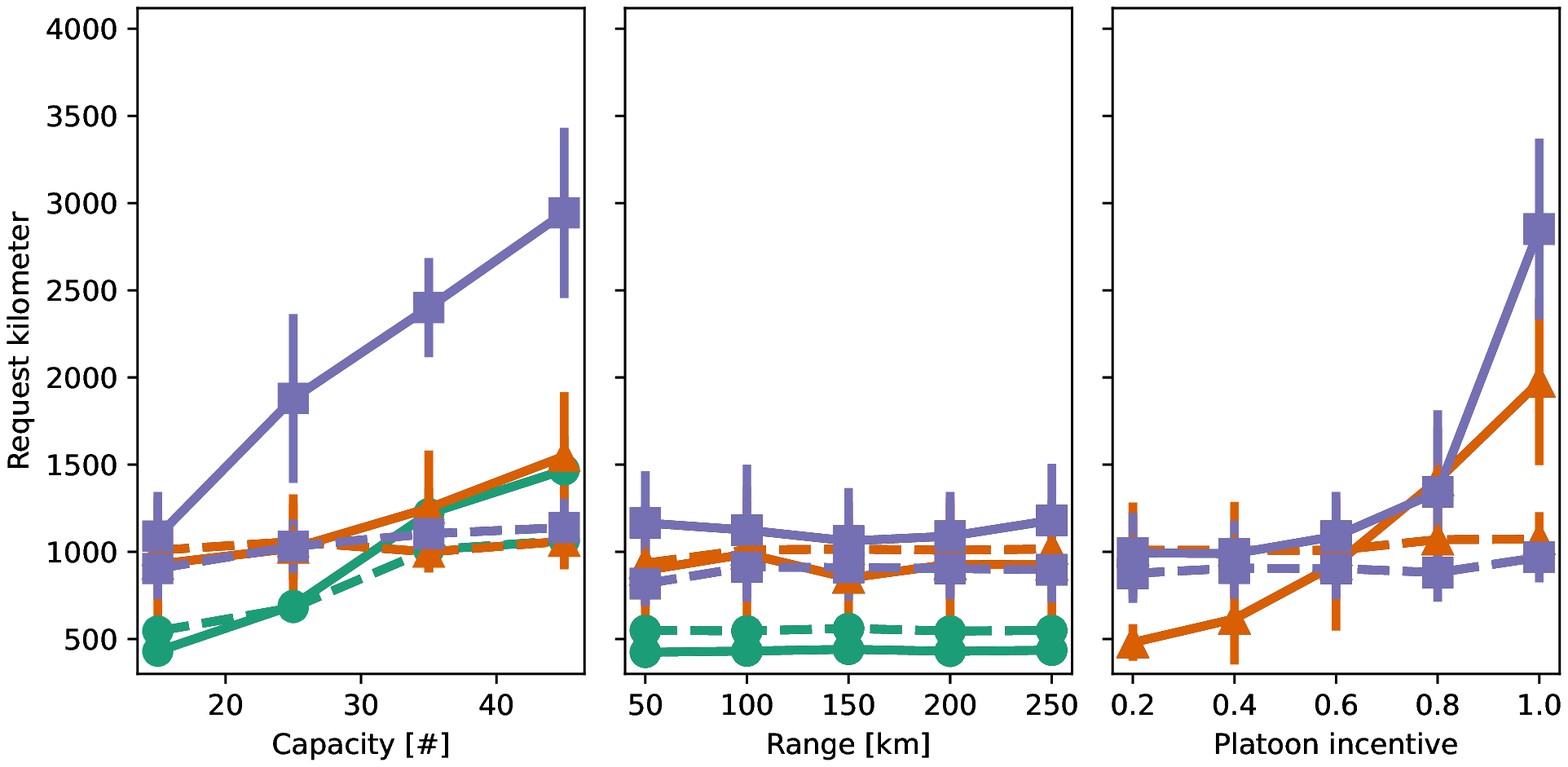}
    	\caption{Impact on Platoon kilometer}
    	\label{fig:unit_parameter_study_Platoon_kilometer}
	\end{subfigure}
	\hspace{0pt}    
	\begin{subfigure}[b]{0.49\textwidth}
		\centering
    	\includegraphics[width=\textwidth]{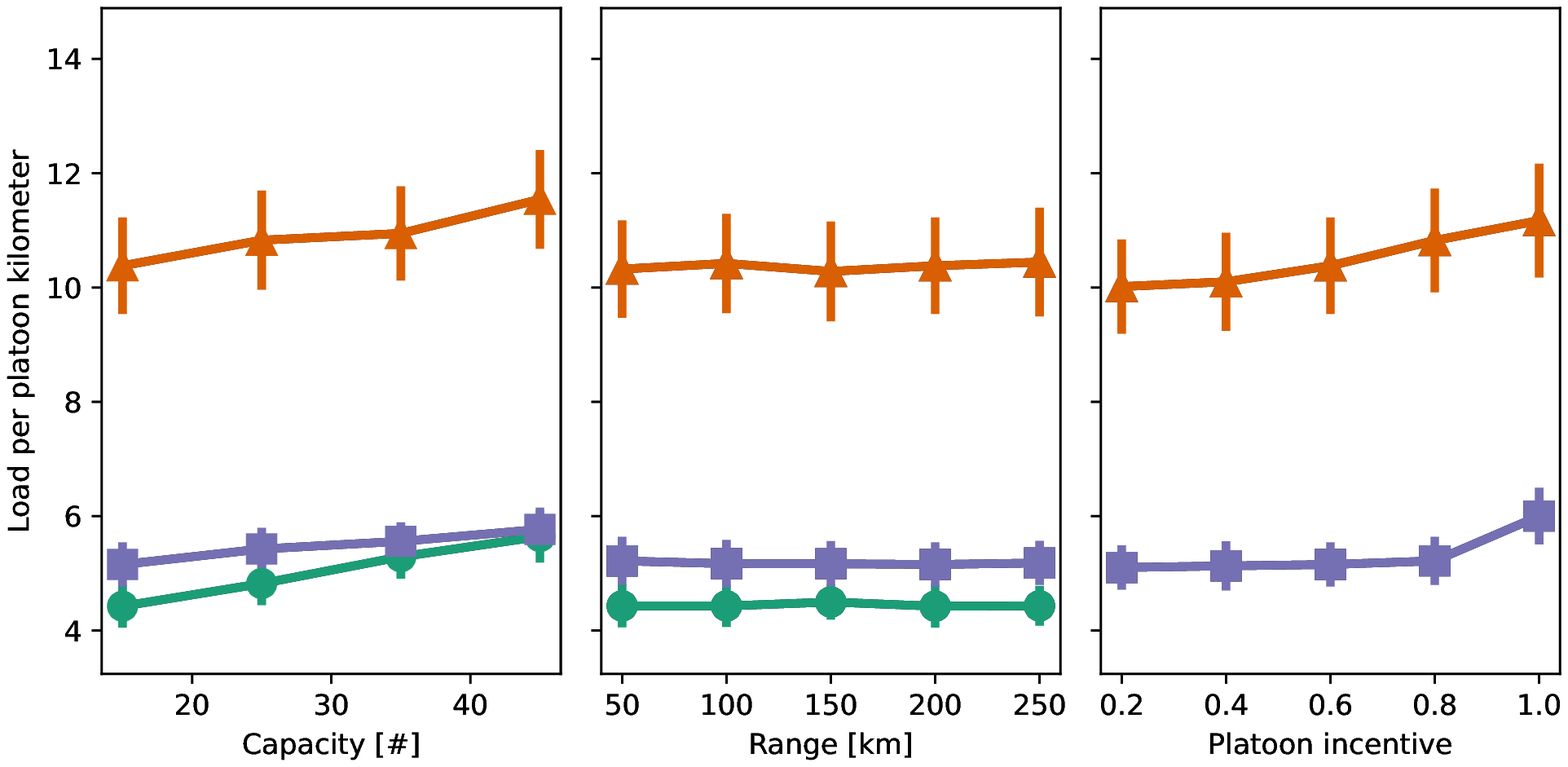}
    	\caption{Impact on Load per platoon kilometer}
    	\label{fig:parameter_study_Profit}
	\end{subfigure}
	\hfill    
	\begin{subfigure}[b]{0.49\textwidth}
		\centering
        \includegraphics[width=\textwidth]{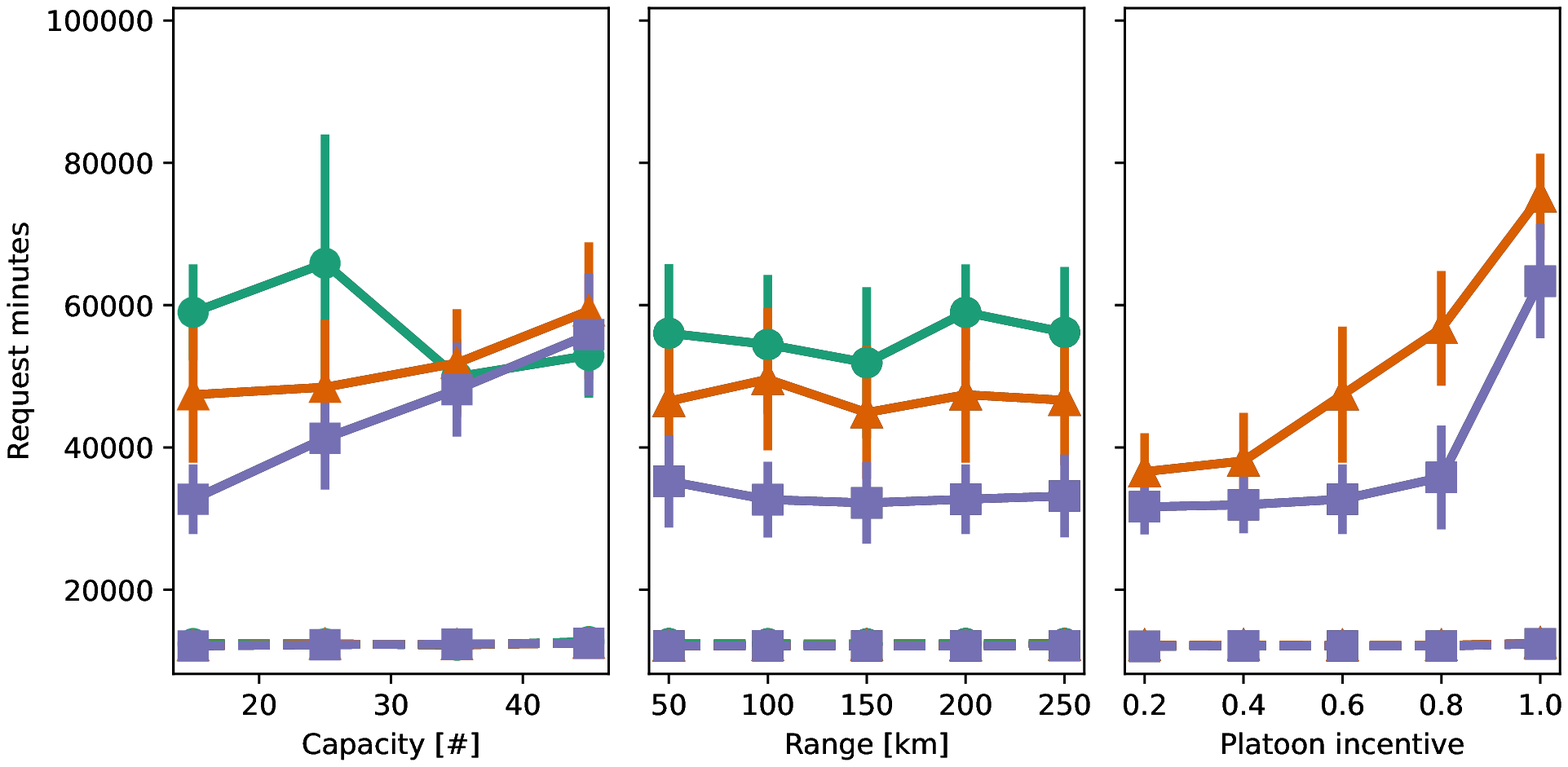}
        \caption{Impact on Request minutes}
        \label{fig:unit_parameter_study_Request_minutes}
	\end{subfigure}
	\caption{Results for the parameter study. The different operations are marked with different colors and symbols, conventional operations (\protect\includegraphics[height=1.5ex]{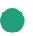}), modular without consolidation (\protect\includegraphics[height=1.5ex]{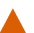}), and modular with consolidation (\protect\includegraphics[height=1.5ex]{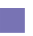}). The demand type is illustrated using different line patterns, freight demand (\protect\includegraphics[height=0.5ex]{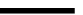}) and passenger demand (\protect\includegraphics[height=0.5ex]{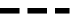}).}
	\label{fig:parameter_analysis}
\end{figure*}

The sensitivity analysis is performed using the distributed and peak demand distribution. By analyzing the graphs in Figures \ref{fig:parameter_study_Objective_} and \ref{fig:parameter_analysis}, it can be seen that the solutions do not depend on the module range settings. We therefore conclude that, first the range setting of \SI{50}{km} is sufficient for the study area chosen and second when providing a sufficiently large fleet of vehicles the range for urban deliveries is not a limiting factor. Different observations can be made for the module capacity and the platoon incentive. 

As can be seen in Figure \ref{fig:parameter_study_Objective_} larger module capacities in general lead to a reduction in overall costs. However, conventional vehicles are more sensitive to vehicle size since the benefits of larger capacities are greater for these vehicles, and a total cost reduction of ~50\% can be observed when increasing the capacity from 15 to 45. The overall cost for modular vehicle operations can also be reduced by ~15\% when increasing the module capacity, however the cost reduction saturates at a module capacity of 35. Indicating a turning point of optimal fleet characteristics. The cost reduction is mainly due to the reduced number of modules needed to serve the same demand (see Figure \ref{fig:unit_parameter_study_Request_Modules_}). For capacities above 35 the disadvantages of fewer models, i.e. longer traveled distances and longer travel times, outweigh the cost savings gained by further reducing the number of modules (see Figures \ref{fig:unit_parameter_study_Platoon_kilometer} and \ref{fig:unit_parameter_study_Request_minutes}). The distances traveled and the travel time increase by 13\% and 7\%, respectively, while the number of modules is reduced by only 5\% from vehicle capacity 35 to 45. The unserved demand stays unchanged throughout all parameter settings. 

A similar observation can be made for the empty kilometers and average platoon length, both decreasing with larger modules. The empty kilometers can be drastically reduced for conventional operations whereas the modular vehicle operations only reduce by a small margin. As expected the optimal platoon length becomes shorter with larger module capacities and as described in the total cost discussions this trend saturates at a module capacity of 35. 

The average fill rate per module is shown in Figure \ref{fig:parameter_study_Fill_rate}. The fill rate for modular operations is higher than for conventional vehicles and is less affected by larger module capacities. Furthermore, with larger module capacities, the fill rate gets lower, indicating less efficient use of the available supply. The jump for conventional operations indicates a shift in operation for larger capacities conventional operations and modular operations with separated demand operate identical. Interestingly, the fill rate for the consolidated operation of passenger and freight demand is lower than if both demand flows are transported separately. 

As shown in Figure \ref{fig:unit_parameter_study_Request_Modules_} and Figure \ref{fig:unit_parameter_study_Request_minutes} the passenger transportation contributes significantly less than the transportation of freight items to the overall cost. The accumulated minutes needed to transport all passengers and the number of passenger modules are lower for all parameter settings and for all types of modular vehicle operations. As discussed in the previous section, this shows that passenger requests can be transported next to the freight transportation without adding much additional trip time. As the module capacity increases, the change in request minutes shows that the total delivery time for each request is increasing. Additionally, the same number of passenger modules is needed in the conventional and modular operation for higher values of vehicle capacity, which minimizes the benefits of modular operations.

Looking at the load transported per platoon kilometer (see Figure \ref{fig:parameter_study_Profit}), the modular and consolidated operations outperform the conventional mode of operation. The load transported per platoon kilometer increases slightly with an increase in module capacity. 

The influence of the platoon incentive in Figures \ref{fig:parameter_study_Objective_} and \ref{fig:parameter_analysis} mainly shows two effects. First, the higher the incentive to form platoons, the longer the platoons are. This generally results in more used modules, since the number of trips is not significantly reduced. Since the parameter for platoon incentive is directly related to the objective value, the total cost for operating such systems decreases. Second, however, for larger platoon incentives, the request kilometers and request minutes increase, leading to a longer average trip distance and average trip duration. A similar negative effect can be seen in relation to the fill rate, which is slightly reduced for higher platoon incentives. This indicates that the reduction of operational costs leads to more supply which in turn is not used efficiently to serve the demand.

\section{Case Study}
\label{sec:case_study}

\begin{figure*}[ht]
	\centering
	\includegraphics[width=\textwidth]{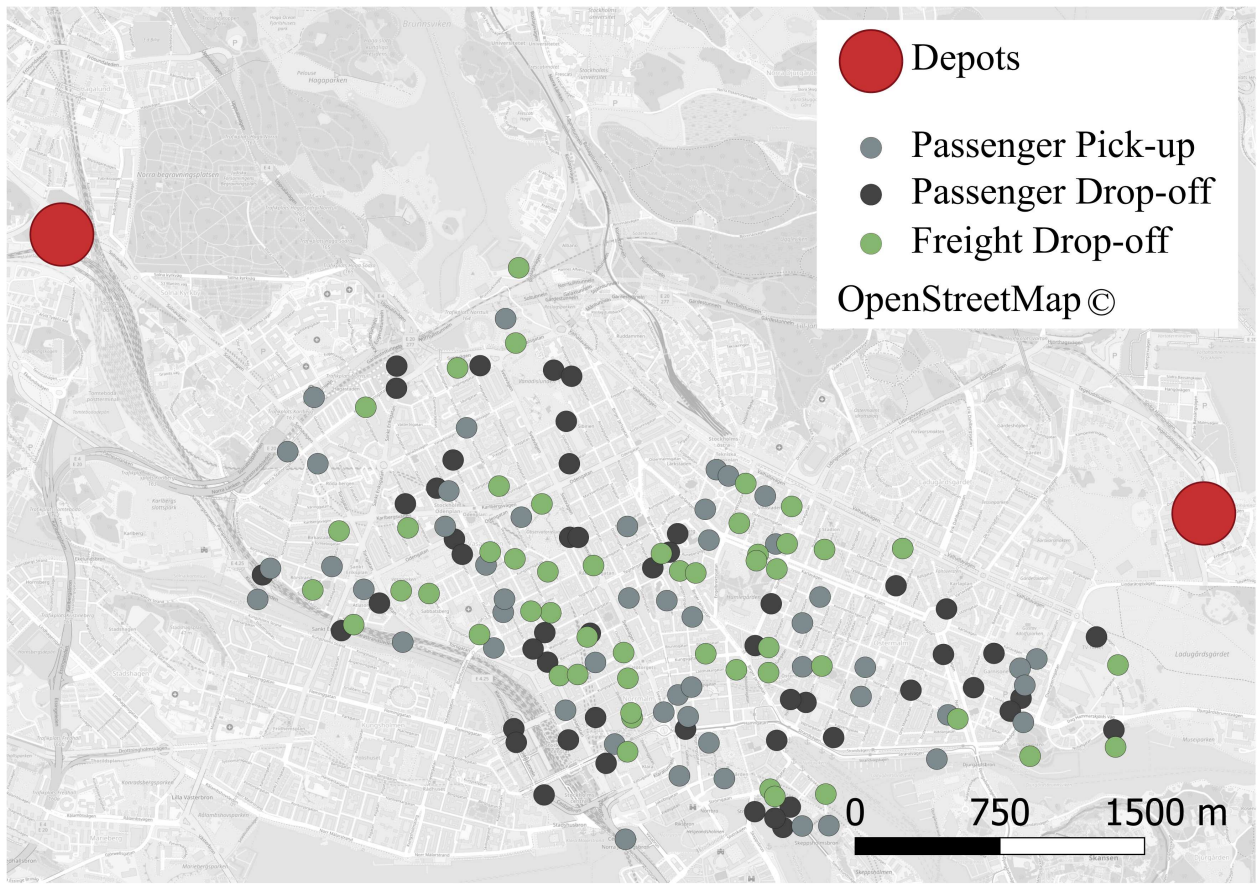}
	\caption{Spatial representation of the case study in Stockholm, Sweden.}
	\label{fig:stockholm_scenario}
\end{figure*}

The case study is situated in Stockholm, Sweden. In Figure \ref{fig:stockholm_scenario} the different pick-up and drop-off locations including the two depots are visualized. Each customer request has a certain demand level and time-windows assigned. The spatial passenger demand pattern and level is representing expected daily passenger movement patterns in Stockholm, while the freight demand is evenly distributed over the case study area. The freight pick-up locations are at the depots. The time windows follow a peak-hour distribution as exemplified in Figure \ref{fig:normal_scenario}. The ALNS and model parameters of this case study are in line with the values as shown in Table \ref{tab:ALNS_parameter} and Table \ref{tab:model_parameter}. The travel distances and travel times as used in the graph representation of the MMP-PDP are based on the underlying road network, for each node pair the shortest path is computed, and its distance recorded. The road network information is extracted from \cite{OpenStreetMap}.

\begin{figure*}[ht]
	\centering
	\includegraphics[width=\textwidth]{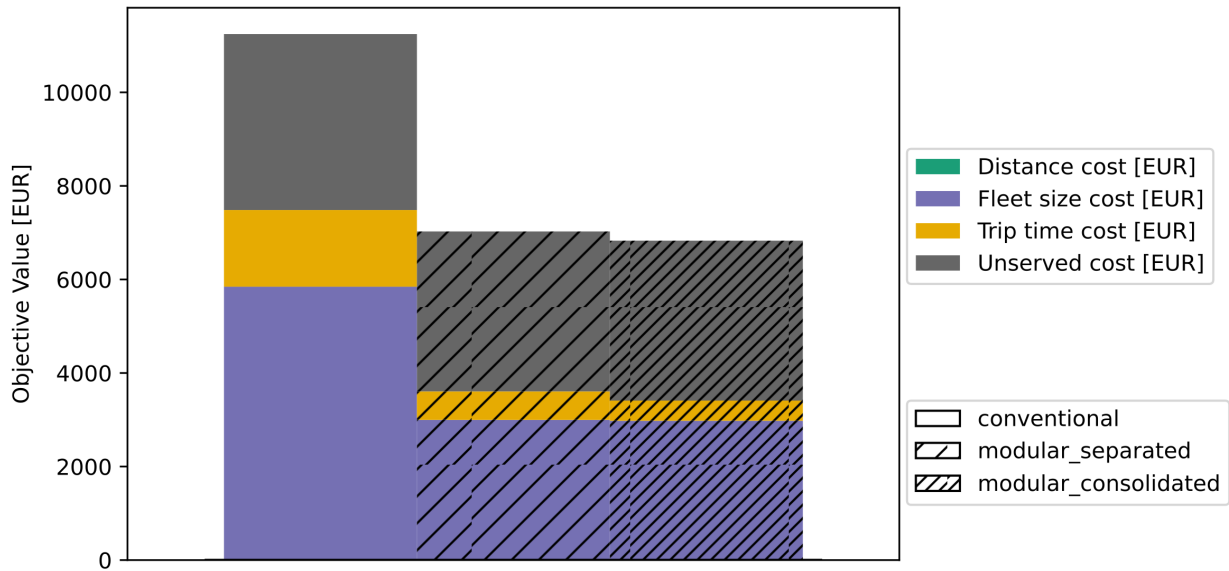}
	\caption{Objective value for different transport systems in the Stockholm case study.}
	\label{fig:stockholm_scenario_study_Objective_Value_}
\end{figure*}

\begin{figure*}[ht]
	\centering
	\begin{subfigure}[b]{0.49\textwidth}
		\centering
    	\includegraphics[width=\textwidth]{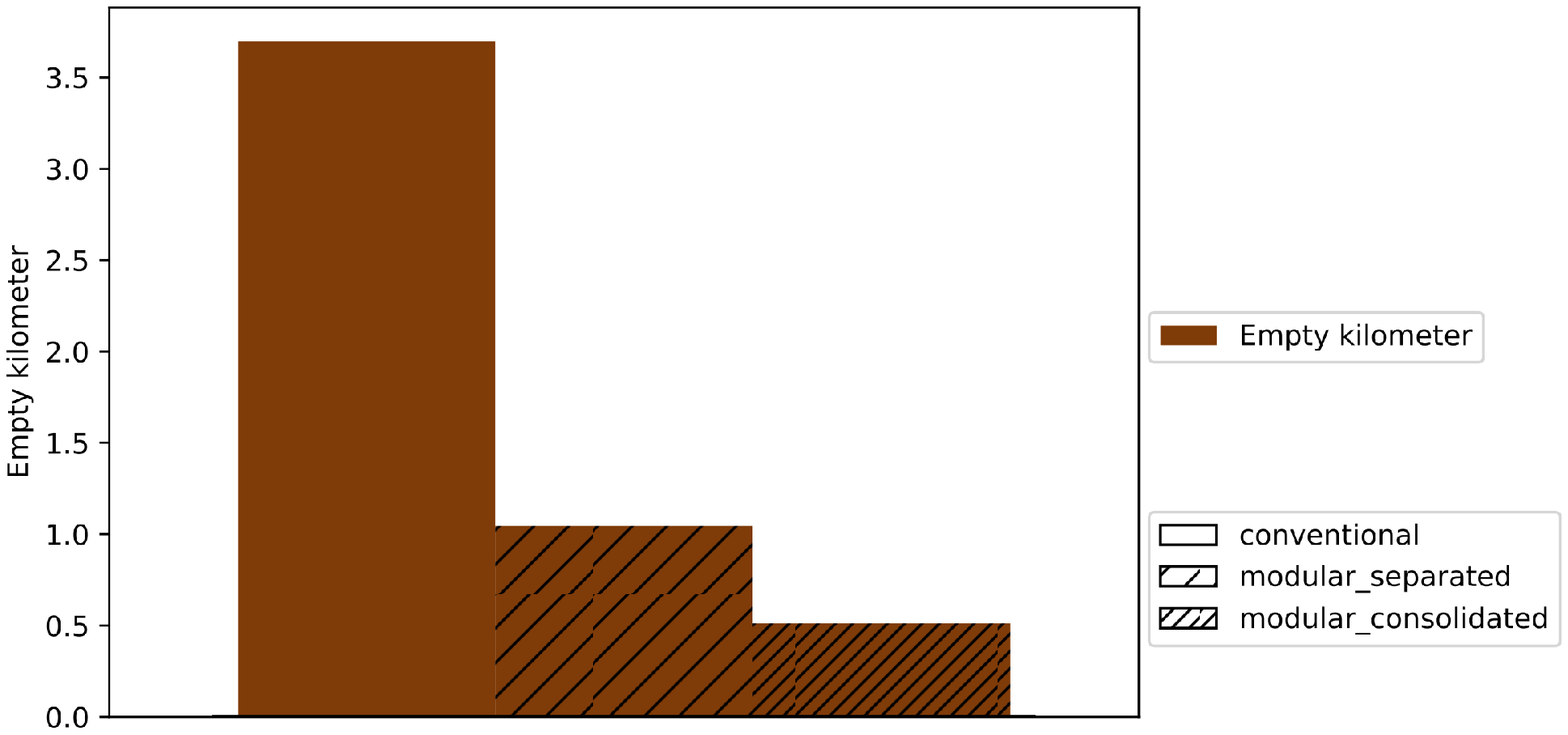}
    	\caption{Effect on empty kilometers traveled}
    	\label{fig:stockholm_scenario_study_Empty_kilometer}
	\end{subfigure}
	\hfill
	\begin{subfigure}[b]{0.49\textwidth}
		\centering
    	\includegraphics[width=\textwidth]{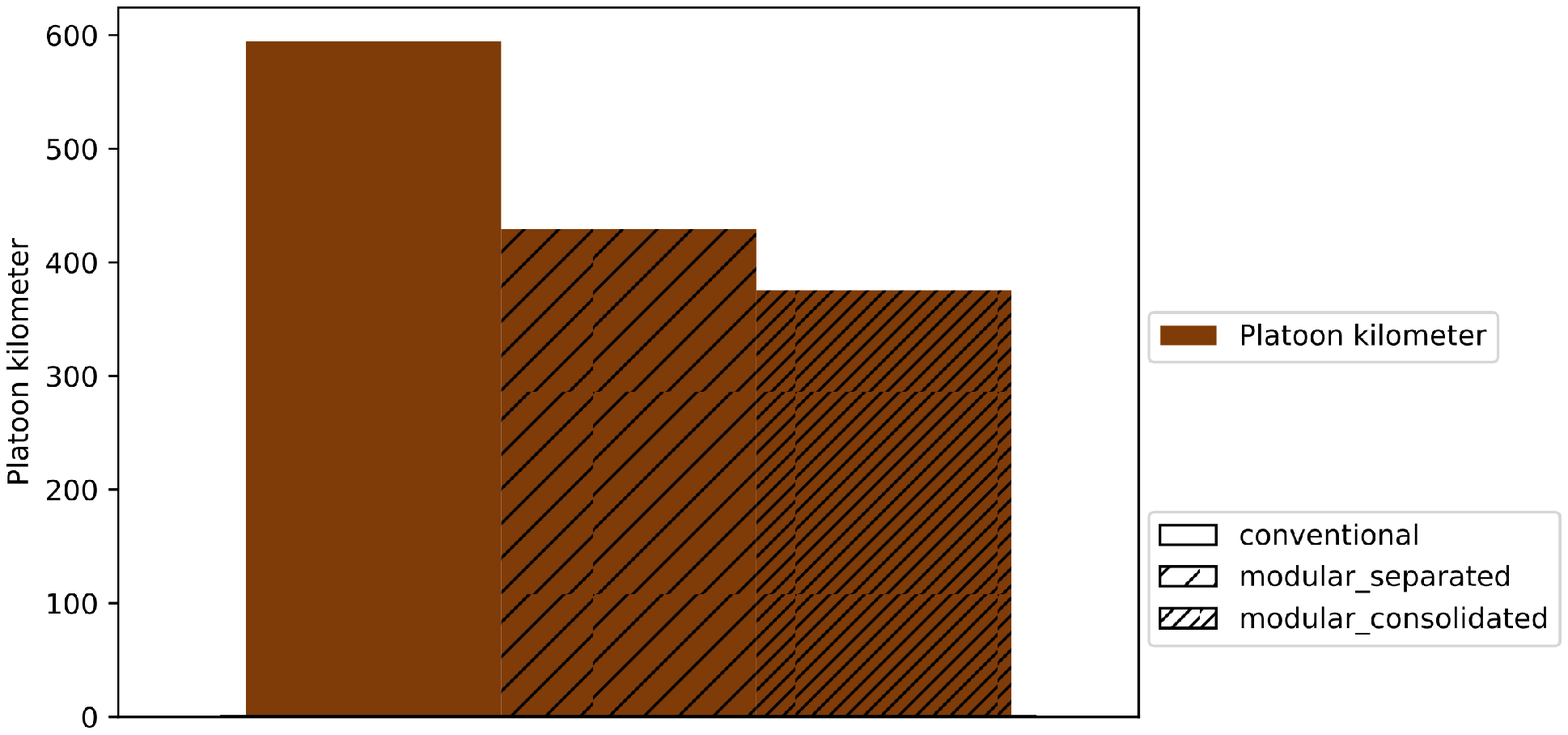}
    	\caption{Effect on Platoon kilometer}
    	\label{fig:stockholm_scenario_study_Platoon_kilometer}
	\end{subfigure}
	\hspace*{0.5\fill}
	\begin{subfigure}[b]{0.49\textwidth}
		\centering
	    \includegraphics[width=\textwidth]{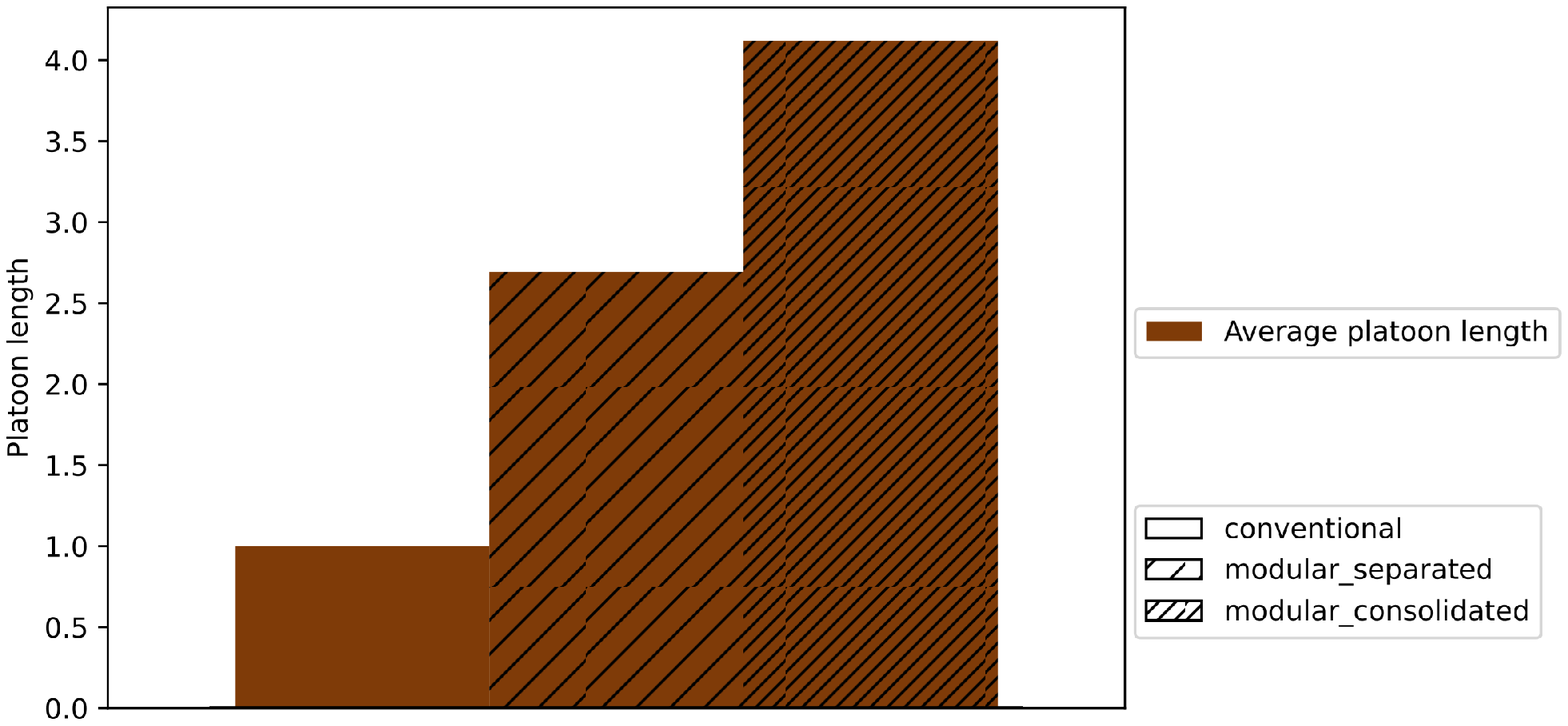}
	    \caption{Effect on Platoon length}
	    \label{fig:stockholm_scenario_study_Platoon_length}
	\end{subfigure}
	\caption{Results for the case study in Stockholm, Sweden. The vehicle operations are indicated by different patterns, conventional (\protect\includegraphics[height=1.5ex]{images/conventional.eps}), modular without consolidation (\protect\includegraphics[height=1.5ex]{images/modular_separated.eps}), and modular with consolidation (\protect\includegraphics[height=1.5ex]{images/modular_collaboration.eps}).}
	\label{fig:stockholm_scenario_analysis}
\end{figure*}

The results of the case study show similar trends as the analysis in Section \ref{sec:results}. The overall objective value can be substantially reduced when operating modular vehicles (see Figure \ref{fig:stockholm_scenario_study_Objective_Value_}). The objective value can be reduced further when consolidating passenger and freight demand. The main cost reduction is due to the reduced fleet size cost and the reduction in trip times. Multiple modules form platoons (see Figure \ref{fig:stockholm_scenario_study_Platoon_length}) and in some scenarios less modules are needed to serve the demand. This does not lead to a higher number of unserved requests compared to the conventional operations, indicating an unchanged level of service. 

In Figure \ref{fig:stockholm_scenario_study_Empty_kilometer} the effect on empty kilometers is shown. The large reduction in empty kilometers for modular operations highlights the higher utilization of available capacity. In modular operations fewer empty trips towards the depot need to be made since the platoons can serve more requests in one trip. Besides the reduction of empty kilometers traveled, this also results in a changed delivery order. The changes in platoon kilometer as shown in Figure \ref{fig:stockholm_scenario_study_Platoon_kilometer} follow a similar trend as discussed in Section \ref{sec:results} and shown in Figure \ref{fig:scenario_study_Platoon_kilometer}. The overall platoon kilometers can drastically be reduced when operating modular vehicles, and additional savings are achieved when consolidating passenger and freight demand. This trend reflects the changed delivery orders and the more efficient serving of the requests, and it indicates the reduced number of modules needed to serve the demand.

In Figure \ref{fig:stockholm_scenario_study_Platoon_length} the average platoon length is shown. The average platoon length is approximately 2.7 for modular and separated operations and approximately 4.1 for modular and consolidated operations. When comparing this to the scenario analysis in Section \ref{sec:results}two observations can be made. First, the general trend of longer platoons for modular and consolidated operations is persistent in all experiments. Second, the higher demand level and different demand patterns in the Stockholm case study led to longer platoons as well as a larger fleet size to serve the demand. 

The results of the case study highlight the overall system improvement when using modular vehicles compared to conventional vehicles. Similar effects as in the generated scenarios can be identified for the empty kilometers, platoon kilometers and platoon length. Including travel distance and travel time based on the underlying road network increases the practical relevance of the study and suggests significant cost savings for a variety of transportation scenarios.

\newpage
\section{Conclusion}
\label{sec:conclusion}
In this work, we study the impact of modular vehicle operations on the urban freight and passenger transportation system. The vehicle concept allows the operator to form platoons by connecting multiple modules. Additionally, different modules can be used to transport different types of demand. Furthermore, we analyze potential efficiency gains by simultaneously optimizing multiple demand flows in the transportation system. A novel modular multi-purpose pickup and delivery problem (MMP-PDP) is proposed to analyze the effects of each vehicle operation. The model extends existing research work by adding the necessary modularity and consolidating formulations to existing PDP formulations. The MMP-PDP is solved using an ALNS algorithm which is validated using CPLEX. We perform a series of experiments to study the impacts of the different vehicle operation types under different temporal and spatial demand configurations. Additionally, a parameter study is performed in order to understand the sensitivity of the results to vehicle capacity, vehicle range, and cost savings due to platoon formation.

The results suggest that cost savings of 48\% can be achieved with modular vehicles, while an additional 9\% can be saved when consolidating passenger and freight demand. These cost reductions are mainly due to the reduction in operating costs for the platoons and reduction of trip duration. The cost savings are insensitive to the demand distribution in space (clustered or distributed) and time (even or peak). The average platoon length for modular operations is 1.8 without consolidation and 2.1 with consolidation. In this study we show that the distance traveled as well as the number of modules required for different scenarios are slightly reduced with modular operations. Furthermore, empty vehicle kilometers can be significantly reduced, indicating a more efficient use of available capacity. The total served requests stay unchanged, representing an unchanged level of service. The results do not suggest additional benefits under uneven temporal demand distributions compared to even demand distributions. Finally, the results indicate that the delivery of freight requests contributes significantly more to the costs than passenger requests, and most cost savings can be achieved through the change of delivery order for freight requests. The combination of requiring fewer modules, lower operational costs, less platoon-kilometer, and the reduction of empty kilometers imply an increase in efficiency for the transport system. Moreover, a reduction in emissions is expected as a direct consequence of the reduced platoon kilometers.

The results of the parameter study show that the range parameter does not have a significant effect on the objective value. However, the total cost can be further reduced when increasing the module capacity. This is mainly due to lower operating costs, while the number of modules, the length of the platoon, and the empty kilometers all decrease. Furthermore, the fill rate is reduced for high platoon incentives, indicating a less efficient use of the supply. An increase in module capacity leads in general to a reduction in total costs for all vehicle operations. However, conventional operations are more sensitive to changes in module capacity. For modular operations, the benefits of larger modules are smaller and saturate at a module capacity of approximately 35 units for the scenarios studied. This cost reduction stems mainly from fewer modules required to serve the same number of requests.

The results of the case study show similar results to the scenario analysis. The objective value can be reduced when operating modular vehicles and is further reduced when consolidating passenger and freight demand. This is mainly due to the reduced fleet size cost and the reduction in trip times. Due to the higher demand and larger scale of this scenario modules form larger platoons and the platoon kilometer are higher however this does not affect the improved performance of the modular transport system. 

The proposed model and the ALNS algorithm allow for an efficient study of different transportation systems; by adjusting the input data and parameter settings, the proposed approach can be used by practitioners and decision makers to evaluate different scenarios specific for their transportation problem. Additionally, the formulation makes it possible to study different policy settings by adjusting the assumed cost parameters and the travel costs on the model graph. 

The presented results show the benefits modular vehicles can have on the consolidated transportation of passenger and freight in urban environments. To enable such vehicle operations additional challenges must be overcome. The formation of platoons at depots requires more time and space compared to conventional vehicle operations. Hence potentially larger and more central depots are required which requires approval of city authorities. Additionally, longer platoons might require larger curbside space when un-/loading requests which can obstruct surrounding traffic and pedestrian movements. As the results show the modular vehicle concept is beneficial for all spatial and temporal demand distributions, hence potential practical shortcomings can be solved by operating such systems during off-peak hours or in separated areas. This can be regulated by transport authorities applying policies supporting the benefits of the modular vehicles while minimizing the obstructions of ongoing traffic. The main benefits of the modular vehicle operations are the reduction of platoon-km by transporting passenger along with freight requests.

Additional consequences of the proposed transportation system affect the design and operation of conventional public transport stops, which could be combined with freight delivery stations to multimodal terminals. Furthermore, the combined optimization of freight and passenger demand requires the integrated planning of multiple entities, such as public transport operators and freight carriers, which traditionally have not been connected. All in all, there is need for planning processes and policies that support the benefits of the modular vehicles while minimizing the obstructions of ongoing traffic. The main benefits of the modular vehicle operations being the reduction of platoon-km by transporting passenger along with freight requests. Finally, the proposed model can be used to study different types of logistics and material flows, that is, forward and reverse flows, including goods, parcels, bulk deliveries and recycling.

The modularity of the proposed model is limited mainly by one assumption made in the problem description. The platoon configuration is required to stay unchanged throughout the route; hence, en-route platoon reconfiguration is not modeled. On one hand this limits the theoretical benefits of the vehicle technology, one the other hand this assumption increases the practical relevance of the model. Since en-route reconfiguration is not possible with current technologies such systems cannot be operated in a real-world scenarios. Vehicle operations as proposed in this study can be operated in a real-world scenario since there is the technology needed and no changes in infrastructure or legislation are required.

However, the proposed results should be considered as scenario specific. Mainly, the values for the optimal capacity and range-related impacts can be different for different input data. Furthermore, the presented results depend on the assumed cost terms for vehicle fixed costs and operational costs, which are uncertain and likely to change in the future, thereby affecting cost savings and efficiency gains. If the vehicle fixed costs can be reduced in the future, then the cost savings might be further increased. However, if we assume a reduction in vehicle operating costs, the benefits of the proposed modular system will be reduced and will mainly stem from the consolidation of demand. 

Analyzing data from pilot projects which utilize the assumed modules and operations will allow to refine these parameter settings in the future. Furthermore, we assume that all demand is known in advance. For the practical application of the model additional online planning steps which consider e.g. stochastic demand or new requests during the day would have to be added to the framework. Since most of freight and passenger demand is known or can be estimated before the planning period, we expect the effect of considering such unknown demand on the benefits of the modular system to be minimal. Nevertheless this assumption does influence the solutions generated and reduce their practical relevance. An extension of the proposed framework considering online planning as a post-processing step should therefore be considered in future studies.

Future research may extend the proposed model by allowing for dynamic transshipments of requests between different modules and by allowing for en-route platoon reconfiguration. Both of these features would allow further exploiting the flexibility and modularity of the modules, which is likely to further reduce the total costs. The proposed model can be further extended by adding multi-echelon delivery. By replicating consolidation centers and the possibility to model the interaction between different shipping operators. Another future research direction is to formulate the problem using time-expanded networks. Such a formulation might be helpful in solving larger instances with complex time window constraints and enables to model time-dependent travel times. Finally, in future studies demand patterns including return flows, and additional types of demand (e.g. waste, recycling) can be added to the scenario definitions to study other future urban transport systems.

\section{Acknowledgement}
This work was in part supported by the Sustainable and Integrated Transport Systems (Hållbara och Integrerade Transport System HITS) project under grant number 2020-00565 (Vinnova). The computations were enabled by resources provided by the Swedish National Infrastructure for Computing (SNIC) at HPC2N partially funded by the Swedish Research Council through grant agreement no. 2018-05973.

\hspace{1cm}

\bibliography{paper_4}

\end{document}